\documentclass[11pt]{article}

\usepackage{amsmath,amssymb,amsfonts}
\usepackage{pstricks,pst-node}

\numberwithin{equation}{section}
\newtheorem{thm}{Theorem}[section] 
\newtheorem{prp}[thm]{Proposition}
\newtheorem{lmm}[thm]{Lemma}   
\newtheorem{crl}[thm]{Corollary} 
\newtheorem{dfn}[thm]{Definition}

\renewcommand{\Bbb}{\mathbb}

\def\A{\mathcal A}
\def\C{\mathbb C}
\def\cC{\mathcal C}

\def\cD{\mathcal D}
\def\E{\mathbb E}
\def\F{\mathfrak F}
\def\cF{\mathcal F}
\def\G{\mathfrak G}
\def\H{\mathcal H}
\def\I{\mathfrak i}
\def\M{\mathfrak M}
\def\cM{\mathcal M}
\def\O{\mathcal O}
\def\Q{\mathbb Q}
\def\R{\mathbb R}
\def\T{\mathcal T}
\def\U{\mathcal U}

\def\W{\mathcal W}
\def\X{\mathfrak X}
\def\Z{\mathbb Z}
\def\cZ{\mathcal Z}

\def\e_ref#1{(\ref{#1})}
\def\under#1{\underline{#1}}

\def\ov#1{\overline{#1}}
\def\wt#1{\widetilde{#1}}
\def\ti#1{\tilde{#1}}
\def\sf#1{\textsf{#1}}

\def\lan{\langle}
\def\ran{\rangle}
\def\lr#1{\lan#1\ran}
\def\blr#1{\big\lan#1\big\ran}
\def\llrr#1{\lan\!\lan#1\ran\!\ran}
\def\bllrr#1{\big\lan\!\big\lan#1\big\ran\!\big\ran}
\def\lra{\longrightarrow}
\def\Lra{\Longrightarrow}

\def\al{\alpha}

\def\de{\delta}
\def\ep{\epsilon}
\def\ga{\gamma}
\def\io{\iota}
\def\ka{\kappa}
\def\la{\lambda}
\def\na{\nabla}
\def\om{\omega}
\def\si{\sigma}
\def\ups{\upsilon}
\def\ze{\zeta}
\def\ve{\varepsilon}

\def\vt{\vartheta}

\def\De{\Delta}
\def\Ga{\Gamma}
\def\La{\Lambda}
\def\Om{\Omega}
\def\Si{\Sigma}

\def\Aut{\textnormal{Aut}}
\def\const{\textnormal{const}}
\def\ev{\textnormal{ev}}
\def\GW{\textnormal{GW}}
\def\gd{\textnormal{gd}}
\def\id{\textnormal{id}}
\def\Id{\textnormal{Id}}

\def\End{\textnormal{End}}
\def\mn{\textnormal{mn}}
\def\Hol{\textnormal{Hol}}
\def\Hom{\textnormal{Hom}}
\def\Im{\textnormal{Im}~\!}
\def\es{\textnormal{es}}
\def\PD{\textnormal{PD}}
\def\ses{\textnormal{ses}}
\def\rk{\textnormal{rk}}
\def\st{\textnormal{s.t.}}

\def\P{\Bbb{P}^n}
\def\bP{\mathbb P}
\def\i{\infty}
\def\eset{\emptyset}
\def\bpar{\bar\partial}

\begin{document}

\title{Reduced Genus-One Gromov-Witten Invariants} 

\author{Aleksey Zinger\thanks{Partially supported by an NSF Postdoctoral Fellowship}}

\date{\it Revised: March 5, 2008}
\maketitle

\begin{abstract}
\noindent
In a previous paper we described a natural closed subset $\ov\M_{1,k}^0(X,A;J)$ of 
the moduli space $\ov\M_{1,k}(X,A;J)$ of stable genus-one $J$-holomorphic 
maps into a symplectic manifold~$X$.
In this paper we generalize the definition of the main component to moduli spaces 
of perturbed, in a restricted way, $J$-holomorphic maps.
This generalization implies that $\ov\M_{1,k}^0(X,A;J)$, just like $\ov\M_{1,k}(X,A;J)$, 
carries a virtual fundamental class and can be used to define symplectic invariants.
These truly genus-one invariants constitute part of the standard 
genus-one Gromov-Witten invariants, 
which arise from the entire moduli space $\ov\M_{1,k}(X,A;J)$.
The new invariants are more geometric  and can be used to compute 
the genus-one GW-invariants of complete intersections, as shown in a separate paper.
\end{abstract}

\thispagestyle{empty}

\tableofcontents

\section{Introduction}
\label{intro_sec}

\subsection{Background and Motivation}
\label{back_subs}

\noindent
Let $(X,\om,J)$ be a compact almost Kahler manifold.
In other words, $(X,\om)$ is a symplectic manifold and 
$J$ is an almost complex structure on $X$ tamed by $\om$,~i.e.
$$\om(v,Jv)>0 \qquad\forall v\in TX-X.$$
If $g,k$ are nonnegative integers and $A\!\in\!H_2(X;\Z)$,
we denote by $\ov\M_{g,k}(X,A;J)$ the moduli space of (equivalence classes of) 
stable $J$-holomorphic maps from genus-$g$ Riemann surfaces with $k$ marked points 
in the homology class~$A$.
Let $\M_{g,k}^0(X,A;J)$ be the subspace of $\ov\M_{g,k}(X,A;J)$
consisting of the stable maps $[\cC,u]$ such that the domain $\cC$
is a smooth Riemann surface.
The compact moduli space $\ov\M_{g,k}(X,A;J)$ was constructed 
in order to ``compactify" $\M_{g,k}^0(X,A;J)$ and to define invariants of $(X,\om)$
enumerating $J$-holomorphic curves of genus~$g$ in~$X$.
If $g\!=\!0$, $(X,\om;A)$ is positive in a certain sense, and $J$ is generic, 
then $\M_{g,k}^0(X,A;J)$ is a dense open subset of $\ov\M_{g,k}(X,A;J)$ and
the corresponding Gromov-Witten invariants do indeed count genus-zero 
$J$-holomorphic curves in~$X$;
see \cite[Chapter~7]{McSa} and \cite[Sections~1,9]{RT}, for example.
However, if $g\!\ge\!1$, it is usually the case that $\M_{g,k}^0(X,A;J)$ is not dense in 
$\ov\M_{g,k}(X,A;J)$ and the genus-$g$ GW-counts include 
$J$-holomorphic curves of lower genera.\\

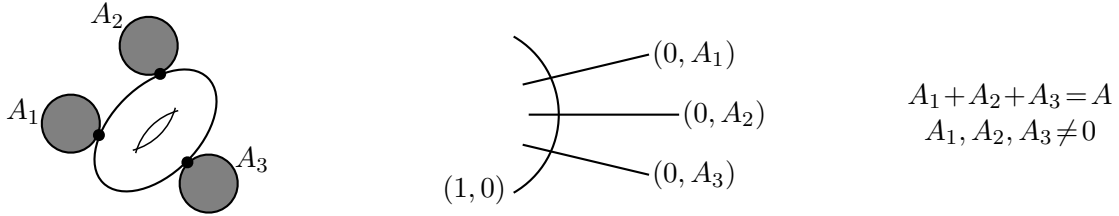
\begin{figure}
\begin{pspicture}(-1.1,-1.8)(10,1.25)
\psset{unit=.4cm}
\rput{45}(0,-4){\psellipse(5,-1.5)(2.5,1.5)
\psarc[linewidth=.05](5,-3.3){2}{60}{120}\psarc[linewidth=.05](5,0.3){2}{240}{300}
\pscircle[fillstyle=solid,fillcolor=gray](5,-4){1}\pscircle*(5,-3){.2}
\pscircle[fillstyle=solid,fillcolor=gray](6.83,.65){1}\pscircle*(6.44,-.28){.2}
\pscircle[fillstyle=solid,fillcolor=gray](3.17,.65){1}\pscircle*(3.56,-.28){.2}}
\rput(.2,-.9){$A_1$}\rput(3,2.3){$A_2$}\rput(7.8,-2.5){$A_3$}
\psarc(15,-1){3}{-60}{60}\psline(17,-1)(22,-1)\psline(16.8,-2)(21,-3)\psline(16.8,0)(21,1)
\rput(15.2,-3.5){$(1,0)$}\rput(22.5,1){$(0,A_1)$}
\rput(23.5,-1){$(0,A_2)$}\rput(22.5,-3){$(0,A_3)$}
\rput(33,-1){\begin{tabular}{c}$A_1\!+\!A_2\!+\!A_3\!=\!A$\\
$A_1,A_2,A_3\!\neq\!0$\end{tabular}}
\end{pspicture}
\caption{The domain of an element of $\M_{1,k}^3(X,A;J)$}
\label{m3_fig}
\end{figure}

\noindent
If $g\!=\!1$ and $(X,\om;A)$ is positive, the above deficiencies are due exclusively 
to the presence of large subspaces of stable maps $[\cC,u]$ in $\ov\M_{1,k}(X,A;J)$
such that $u$ is constant on the principal components of $\cC$,
i.e.~the irreducible components that carry the genus of~$\cC$.
More precisely, if $m$ is a positive integer, let $\M_{1,k}^m(X,A;J)$ be the subset
of $\ov\M_{1,k}(X,A;J)$ consisting of the stable maps $[\cC,u]$
such that $\cC$ is a smooth genus-one curve~$\cC_P$ with $m$~rational 
components attached directly to~$\cC_P$, $u|_{\cC_P}$ is constant, and
the restriction of~$u$ to each rational component is non-constant.
Figure~\ref{m3_fig} shows the domain of an element of $\M_{1,k}^3(X,A;J)$,
from the points of view of symplectic topology and of algebraic geometry.
In the first diagram, each shaded disc represents a sphere;
the homology class next to each rational component 
$\cC_i$ indicates the degree  of~$u|_{\cC_i}$.
In the second diagram, the components of $\cC$ are represented by curves,
and the pair of indices next to each component $\cC_i$ shows 
the genus of $\cC_i$ and the degree of~$u|_{\cC_i}$.
We denote by $\ov\M_{1,k}^m(X,A;J)$ the closure of $\M_{1,k}^m(X,A;J)$
in $\ov\M_{1,k}(X,A;J)$.
The image $u(\cC)$ of an element of $\ov\M_{1,k}^m(X,A;J)$ is a genus-zero,
instead of genus-one, $J$-holomorphic curve in~$X$.
We note that if $J$ is sufficiently regular, then
\begin{alignat*}{2}
&\dim\M_{1,k}^0(X,A;J)  ~&=&~
2\big( \lr{c_1(TX),A}+k\big)\equiv\dim_{1,k}(X,A) \qquad\hbox{and}\\
&\dim\M_{1,k}^m(X,A;J) ~&=&~ 
\dim_{1,k}(X,A)+2(n\!-\!m),
\end{alignat*}
where $2n$ is the real dimension of $X$.
Thus, the complement of $\M_{1,k}^0(X,A;J)$ in $\ov\M_{1,k}(X,A;J)$
contains subspaces of dimension at least as large as the dimension of
$\M_{1,k}^0(X,A;J)$, as long as $n\!\ge\!1$, i.e.~$X$ is not a finite collection
of points.\\

\noindent
In \cite[Definition~1.1]{g1comp}, 
we describe a subset $\ov\M_{1,k}^0(X,A;J)$ of $\ov\M_{1,k}(X,A;J)$,
for an arbitrary compact almost Kahler manifold $(X,\om,J)$,
obtained from $\ov\M_{1,k}(X,A;J)$ by discarding most elements of
the spaces $\ov\M_{1,k}^m(X,A;J)$ with $m\!\le\!n$.
In particular, $\ov\M_{1,k}^0(X,A;J)$ contains $\M_{1,k}^0(X,A;J)$.
By \cite[Theorem~1.2]{g1comp}, $\ov\M_{1,k}^0(X,A;J)$ 
is a closed subset of $\ov\M_{1,k}(X,A;J)$ and thus is compact.
If $(X,\om;A)$ is positive in the same sense as in the genus-zero case
and $J$ is generic, $\M_{1,k}^0(X,A;J)$ is a dense open subset of $\ov\M_{1,k}^0(X,A;J)$.
In addition, $\ov\M_{1,k}^0(X,A;J)$ carries a rational fundamental class,
which can be used to define a symplectic invariant of $(X,\om)$
counting genus-one $J$-holomorphic curves in $X$, 
without any genus-zero contribution in contrast to the standard Gromov-Witten invariants;
see \cite[Subsection~1.3]{g1comp}.
Unlike the genus-zero case, $\ov\M_{1,k}^0(X,A;J)$ has the topological structure
of a singular, instead of smooth, orbivariety.\\

\noindent
A $J$-holomorphic map into $X$ is a smooth map $u$ from a Riemann surface $(\Si,j)$ 
that satisfies the Cauchy-Riemann equation corresponding to~$(J,j)$:
$$\bar{\partial}_{J,j}u\equiv \frac{1}{2}\big(du+J\circ du\circ j\big) =0.$$
The Riemann surface $(\Si,j)$ may have simple nodes.
In this paper we generalize the results of~\cite{g1comp} to smooth maps $u$,  from 
genus-one Riemann surfaces, that satisfy a family of perturbed Cauchy-Riemann equations:
$$\bar{\partial}_{J,j}u +\nu(u)=0.$$
The perturbation term $\nu(u)$ is a section of the vector bundle
$$\La^{0,1}_{J,j}T^*\Si\!\otimes\!u^*TX
\equiv\big\{\eta\!\in\!\hbox{Hom}_{\R}(T\Si,u^*TX)\!: J\circ\eta=-\eta\circ j\big\}
\lra \Si.$$
We will study the moduli space $\ov\M_{1,k}(X,A;J,\nu)$ of 
\sf{$(J,\nu)$-holomorphic maps}, i.e.~of solutions to the perturbed 
Cauchy-Riemann equations, for a continuous family $\nu\!=\!\nu(u)$
chosen from a proper linear subspace of the space of all such families;
see Definition~\ref{pert_dfn}.
A key condition on $\nu$ will be that if the degree of $u$ restricted to
the principal components of $\Si$ is zero, then the restriction of $\nu(u)$ 
to the principal components and all nearby degree-zero bubble components is also zero.
Such a family $\nu\!=\!\nu(u)$ will be called \sf{effectively supported}.\\

\noindent
We will show that if $\nu$ is sufficiently small and effectively supported,
then the moduli space $\ov\M_{1,k}(X,A;J,\nu)$ contains a natural closed
subspace $\ov\M_{1,k}^0(X,A;J,\nu)$ containing $\M_{1,k}^0(X,A;J,\nu)$,
i.e.~the subspace of maps with smooth domains;
see Definition~\ref{degen_dfn} and Theorem~\ref{comp_thm}.
For a generic choice of~$\nu$, the ``boundary" of $\ov\M_{1,k}^0(X,A;J,\nu)$
is of real codimension two and thus
$\ov\M_{1,k}^0(X,A;J,\nu)$ determines a rational homology class.
This \sf{virtual fundamental class} (VFC) for $\ov\M_{1,k}^0(X,A;J)$ 
does not change under small changes in~$\nu$ and is an invariant of~$(X,\om)$.
It can be used to define new GW-style invariants, which we denote by~$\GW_{1,k}^0$.
These invariants differ from the standard GW-invariants by a combination
of the genus-zero GW-invariants of~$X$; see~Subsection~\ref{gwdiff_subs} below
for some special cases.\\

\noindent
We note that effectively supported families $\nu\!=\!\nu(u)$ are in no sense generic
in the space of all families.
If fact, for a generic $\nu$, $\M_{1,k}^0(X,A;J,\nu)$ is dense in $\ov\M_{1,k}(X,A;J,\nu)$,
and the latter space determines the standard GW-invariants of $(X,\om)$.
In particular, the statements of the previous paragraph do not hold for 
a generic family $\nu$ of perturbations.\\

\noindent
An algebraic approach to reduced genus-one GW-invariants is suggested 
by Vakil and the author at the end of~\cite{VZ2}.
It still remains to verify that the resulting algebraic invariants agree
with the symplectic ones defined in this paper (whenever the target space
is a smooth algebraic variety), but this should be deducible 
from the desingularizations for certain natural sheaves constructed by
Vakil and the author in~\cite[Section~5]{VZ}.\\

\noindent
Since the symplectic invariants arising from the moduli space $\ov\M_{1,k}^0(X,A;J,\nu)$,
with $\nu$ effectively supported, are closely related to the standard GW-invariants,
they do not in principle carry any new information.
In practice, they behave better geometrically.
In particular, Li and the author show in~\cite{LZ} that there is a simple
relation between {\it reduced} genus-one GW-invariants of a projective complete
intersection and twisted {\it reduced} genus-one GW-invariants of the ambient space.
This relation mimics the corresponding well-known relation in genus zero
(see \cite[(1.2)]{LZ}, for example), but no relation in positive genera
had been even conjectured until~\cite{LZ}.
Combining \cite[Theorem~1.1]{LZ} for the reduced genus-one invariants 
constructed in this paper with the desingularization 
of Vakil and the author in~\cite{VZ} and Theorem~\ref{gwdiff_thm} below, 
the author finally confirms the 1993 mirror symmetry of 
Bershadsky-Cecotti-Ooguri-Vafa \cite{BCOV}
for the genus-one GW-invariants of a quintic threefold\footnotemark
in~\cite{bcov1}.
\footnotetext{This is 
the genus-one analogue of the 1991 genus-zero mirror symmetry prediction 
of Candelas-de la Ossa-Green-Parkes \cite{CDGP},
which was proved in several different ways in the mid to late 90s.}\\

\noindent
In Subsection~\ref{pertmaps_subs}, we describe a geometric reinterpretation of 
the VFC constructions of Fukaya-Ono~\cite{FuO} and Li-Tian~\cite{LT} 
which is well suited for defining a VFC for $\ov\M_{1,k}^0(X,A;J)$.
We state the main results of this paper in Subsection~\ref{res_subs}.
In Subsections~\ref{notation0_subs} and~\ref{notation1_subs},
we generalize the setup of~\cite{g1comp} for $J$-holomorphic maps to 
$(J,\nu)$-holomorphic maps.
In Subsection~\ref{str_subs}, we state three propositions that together are equivalent
to Theorem~\ref{comp_thm}.
They are proved in Subsections~\ref{g0prp_subs} and~\ref{comp1prp_subs}
by extending some of the analytic arguments of \cite[Sections 3,4]{g1comp} 
to the present situation.
The difference between the standard and reduced GW-invariants is analyzed in
Section~\ref{gwdiff_sec}; see also the next subsection.

\subsection{Standard vs.~Reduced Gromov-Witten Invariants}
\label{gwdiff_subs}

\noindent
From the construction of VFC for $\ov\M_{1,k}^0(X,A;J)$
in Subsection~\ref{res_subs}, it is immediate that 
the difference between the standard and reduced genus-one GW-invariants of $X$
must be a combination of the genus-zero GW-invariants of~$X$.
The exact form of this combination can be determined in each specific case
from Proposition~\ref{bdcontr_prp}.
In this subsection, we give an explicit expression for the difference
between the standard and reduced genus-one GW-invariants in the two simplest cases.\\

\noindent
For each $l\!=\!1,\ldots,k$, let
$$\ev_l\!: \ov\M_{g,k}(X,A;J)\lra X, \qquad 
\big[\Si,y_1,\ldots,y_k;u\big] \lra  u(y_l),$$
 the evaluation map at the $l$th marked point.
We will call a cohomology class $\psi$  on $\ov\M_{g,k}(X,A;J)$ \sf{geometric}
if $\psi$ is a product of the classes $\ev_l^*\mu_l$ for $\mu_l\!\in\!H^*(X;\Z)$. 
We denote by $\bar\Z^+$ the set of nonnegative integers.

\begin{thm}
\label{gwdiff_thm}
Suppose $(X,\om)$ is a compact symplectic manifold,
$A\!\in\!H_2(X;\Z)^*$, $k\!\in\!\bar\Z^+$.
If $J$ is an $\om$-compatible almost complex structure on~$X$ and
$\psi$ is a geometric cohomology class on $\ov\M_{1,k}(X,A;J)$, then
$$\GW_{1,k}^X(A;\psi)-\GW_{1,k}^{0;X}(A;\psi)
=\begin{cases}
0,& \hbox{if}~\dim_{\R}X\!=\!4;\\
\frac{2-\lr{c_1(TX),A}}{24}\GW_{0,k}^X(A;\psi),& \hbox{if}~\dim_{\R}X\!=\!6.
\end{cases}$$\\
\end{thm}

\noindent
Theorem~\ref{gwdiff_thm} is proved in Section~\ref{gwdiff_sec}
by studying the obstruction theory along each stratum of the moduli
space $\ov\M_{1,k}(X,A;J,\nu)$, after capping it with~$\psi$.
This proof generalizes to higher-dimensional manifolds~$X$ and 
more general cohomology classes; an explicit formula is obtained by the author
in~\cite{g1diff}.
In fact, for geometric cohomology classes and higher-dimensional manifolds~$X$,
the difference is given by an expression similar to the correction term 
in \cite[Theorem~1.1]{g1}; this can be  seen a priori 
from Proposition~\ref{bdcontr_prp}  and \cite[Subsection~3.2]{g1}.
A special case of Theorem~\ref{gwdiff_thm} is \cite[Theorem~A]{G};
its proof has not yet appeared.\\

\noindent
Theorem~\ref{gwdiff_thm} has a natural, but rather speculative,
generalization to higher-genus invariants.
Suppose that the main component
$$\ov\M_{g,k}^0(X,A;J)\subset \ov\M_{g,k}(X,A;J)$$
is well-defined, as is the main component
$$\ov\M_{g,k}^0(X,A;J,\nu)\subset \ov\M_{g,k}(X,A;J,\nu)$$
for a sufficiently large subspace of perturbations $\nu$ of the $\bpar_J$-operator so that 
$\ov\M_{g,k}(X,A;J,\nu)$ has a regular structure for a generic $\nu$ in this subspace; 
see Subsection~\ref{res_subs} for the $g\!=\!1$ case.
If so, $\ov\M_{g,k}^0(X,A;J)$ carries a virtual fundamental class and
determines \sf{reduced} genus-$g$ GW-invariants $\GW_{g,k}^{0;X}(A;\psi)$.
Theorem~\ref{gwdiff_thm} and its proof should then generalize to higher-genus invariants.
If $\dim_{\R}X\!=\!6$, the expected relationship~is
$$\GW_{g,k}^X(A;\psi)-\GW_{g,k}^{0;X}(A;\psi)
=\sum_{g'=0}^{g-1}C_g^{g'}(\lr{c_1(TX),A})\, \GW_{g',k}^{0;X}(A;\psi),$$
where $\GW_{0,k}^0\!\equiv\!\GW_{0,k}$.
The coefficients $C_g^{g'}(\lr{c_1(TX),A})$ are given by Hodge integrals, 
i.e.~integrals 
of natural cohomology classes on the moduli spaces $\ov\cM_{*,*}$ of curves.
They are of the form expected from the usual obstruction bundle approach.
For example,
\begin{gather*}
C_2^1\big((5\!-\!a)d\big)=-\frac{d(a\!-\!5)}{24},  \\
C_2^0\big((5\!-\!a)d\big)=\frac{1}{2}\bigg(\frac{2\!+\!d(a\!-\!5)}{24}\bigg)^2
+\blr{c(\E^*\!\otimes\!TX)c(L_{2,1}\!\otimes\!T\Bbb{P}^1)^{-1},
\big[\ov\cM_{2,1}\big]\!\times\![\Bbb{P}^1]},
\end{gather*}
where $L_{2,1}\!\lra\!\ov\cM_{2,1}$ is the universal tangent line bundle,
$\E\!\lra\!\ov\cM_{2,1}$ is the rank-two Hodge bundle,
and $\bP^1$ is viewed as a smooth degree-$d$ curve in~$Y$.
The coefficients $C_g^{g'}(\lr{c_1(TX),A})$ can be expressed in terms of 
the numbers $C_{g'}(g\!-\!g',X,A)$ of~\cite{P} and vice versa.

\subsection{Configuration Spaces}
\label{pertmaps_subs}

\noindent
In this subsection we recall certain configuration spaces that are standard
in the theory of Gromov-Witten invariants.
We then define what we mean by effectively supported perturbations
of the $\bar{\partial}_J$-operator that are central to this paper.\\

\noindent
Fix $p\!>\!2$.
Suppose $X$ is a compact manifold, $A\!\in\!H_2(X;\Z)$, and $g,k\!\in\!\bar\Z^+$.
We denote by $\X_{g,k}(X,A)$  the space of equivalence classes
of stable $L^p_1$-maps $u\!:\Si\!\lra\!X$ from genus-$g$ Riemann surfaces
with $k$~marked points, which may have simple nodes, to~$X$ of degree~$A$, i.e. 
$$u_*[\Si]=A\in H_2(X;\Z).$$
Let $\X_{g,k}^0(X,A)$ be the subset of $\X_{g,k}(X,A)$
consisting of the stable maps with smooth domains. 
The spaces $\X_{g,k}(X,A)$ are topologized using $L^p_1$-convergence on compact subsets 
of smooth points of the domain and certain convergence requirements near the nodes;
see  \cite[Section~3]{LT}.
The spaces $\X_{g,k}(X,A)$ can be stratified by the smooth infinite-dimensional orbifolds
$\X_{\T}(X)$ of stable maps from domains of the same geometric type and with
the same degree distribution between the components of the domain;
see Subsections~\ref{notation0_subs} and~\ref{notation1_subs}.
The closure of the main stratum, $\X_{g,k}^0(X,A)$, is $\X_{g,k}(X,A)$.\\

\noindent
If $J$ is an almost complex structure on $X$, let 
$$\Ga_{g,k}^{0,1}(X,A;J)\!\lra\!\X_{g,k}(X,A)$$
be the bundle of $(TX,J)$-valued $(0,1)$ $L^p$-forms. 
In other words, the fiber of $\Ga_{g,k}^{0,1}(X,A;J)$ over a point
$[b]\!=\![\Si,j;u]$ in $\X_{g,k}(X,A)$ is the space
$$\Ga_{g,k}^{0,1}(X,A;J)\big|_{[b]}=\Ga^{0,1}(b;J)\big/\hbox{Aut}(b),
\quad\hbox{where}\quad
\Ga^{0,1}(b;J)=L^p\big(\Si;\La_{J,j}^{0,1}T^*\Si\!\otimes\!u^*TX\big).$$
Here $j$ is the complex structure on $\Si$, the domain of the smooth map~$u$.
The bundle \hbox{$\La_{J,j}^{0,1}T^*\Si\!\otimes\!u^*TX$} over $\Si$
consists of $(J,j)$-antilinear homomorphisms:
$$\La_{J,j}^{0,1}T^*\Si\!\otimes\!u^*TX=\big\{
\eta\!\in\!\hbox{Hom}(T\Si,u^*TX)\!:J\!\circ\!\eta\!=\!-\eta\!\circ\!j\big\}.$$
The total space of the bundle $\Ga_{g,k}^{0,1}(X,A;J)\!\lra\!\X_{g,k}(X,A)$ 
is topologized using $L^p$-convergence on compact subsets of smooth points of the domain
and certain convergence requirements near the nodes.
The restriction of $\Ga_{g,k}^{0,1}(X,A;J)$ to each stratum
$\X_{\T}(X)$ is a smooth vector orbibundle of infinite rank.
Let
$$\G_{g,k}^{0,1}(X,A;J)=
\Ga\big(\X_{g,k}(X,A),\Ga_{g,k}^{0,1}(X,A;J)\big)$$
denote the space of all continuous multisections\footnote{Our term 
\sf{multisection} corresponds to \sf{locally liftable multisection} described by
\cite[Definition 3.5]{FuO}.}
$\nu$ of $\Ga_{g,k}^{0,1}(X,A;J)$
such that the restriction of $\nu$ to each stratum $\X_{\T}(X)$ is smooth.\\

\noindent
We define a continuous section of the bundle 
$\Ga_{g,k}^{0,1}(X,A;J)\!\lra\!\X_{g,k}(X,A)$  by
$$\bar{\partial}_J\big([\Si,j;u]\big) = \bar{\partial}_{J,j}u
= \frac{1}{2}\big(du+J\!\circ\!du\!\circ\!j\big).$$
By definition, the zero set of this section
is the moduli space $\ov\M_{g,k}(X,A;J)$ 
of equivalence classes of stable $J$-holomorphic degree-$A$ maps 
from genus-$g$ curves with $k$~marked points into~$X$.
The restriction of~$\bar\partial_J$ to each stratum of $\X_{g,k}(X,A)$ is smooth.
The section $\bar\partial_J$ of $\Ga_{g,k}^{0,1}(X,A;J)$ 
is Fredholm, i.e.~the linearization of its restriction to every stratum $\X_{\T}(X)$
has finite-dimensional kernel and cokernel
at every point of $\bar{\partial}_J^{-1}(0)\!\cap\!\X_{\T}(X)$.
The index of the linearization of~$\bar{\partial}_J$ at an element of 
$\M_{g,k}^0(X,A;J)$ is the expected dimension $\dim_{g,k}(X,A)$ of 
the moduli space $\ov\M_{g,k}(X,A;J)$.\\

\noindent
If $\nu$ is a sufficiently small element of $\Ga_{g,k}^{0,1}(X,A;J)$, the space
$$\ov\M_{g,k}(X,A;J,\nu)\equiv 
\big\{\bar{\partial}_J\!+\!\nu\big\}^{-1}(0)\subset \X_{g,k}(X,A)$$
is compact, since $\ov\M_{g,k}(X,A;J)$ is.
For a small generic choice of $\nu$, $\ov\M_{g,k}(X,A;J,\nu)$ admits a stratification
by orbifolds of even dimensions; see the first remark below.
The main stratum, 
$$\M_{g,k}^0(X,A;J,\nu)=\ov\M_{g,k}(X,A;J,\nu)\cap\X_{g,k}^0(X,A),$$
is a smooth orbifold of dimension $\dim_{g,k}(X,A)$.
Since $\X_{g,k}(X,A)$ is locally a Banach space,
there exist arbitrary small neighborhoods $U$~of 
$$\ov\M_{g,k}(X,A;J,\nu)-\M_{g,k}^0(X,A;J,\nu)$$
in $\X_{g,k}(X,A)$ such that
$$H_l\big(U;\Q)=\{0\}\qquad\forall~l\ge \dim_{g,k}(X,A)-1.$$
Since $\ov\M_{g,k}(X,A;J,\nu)\!-\!U$ is compact,
via the pseudocycle construction of \cite[Chapter~7]{McSa} and
 \cite[Section~1]{RT}, $\M_{g,k}^0(X,A;J,\nu)$ determines a homology class
\begin{equation*}\begin{split}
\big[\ov\M_{g,k}(X,A;J,\nu)\big]&\in 
H_{\dim_{g,k}(X,A)}(W,U;\Q)\\
&\qquad\qquad \approx H_{\dim_{g,k}(X,A)}(W;\Q),
\end{split}\end{equation*}
for any small neighborhood $W$ of $\ov\M_{g,k}(X,A;J,\nu)$ in $\X_{g,k}(X,A)$.
The isomorphism between the two homology groups is induced by inclusion.
Since $\nu$ can be chosen to be arbitrarily small, this procedure defines a rational 
homology class in an arbitrary small neighborhood of $\ov\M_{g,k}(X,A;J)$ in $\X_{g,k}(X,A)$.
This topological reinterpretation of the VFC constructions 
done in~\cite{FuO} and~\cite{LT}
turns out to be very suitable for constructing a VFC for
the moduli space $\ov\M_{1,k}^0(X,A;J)$.\\

\noindent
{\it Remark 1:} The strata of $\ov\M_{1,k}(X,A;J,\nu)$ locally are 
unions of finitely many smooth suborbifolds of a smooth orbifold.
The branches of the strata correspond to the branches of~$\nu$.
We will call such objects orbifolds, nevertheless, as
these generalized orbifolds are just as suitable for the topological purposes
of~\cite{FuO}, \cite{LT}, and this paper; see  in \cite[Sections 3,4]{FuO}
for details.\\

\noindent
{\it Remark 2:} The above construction defines a homology class
$$\Om_W\in H_{\dim_{g,k}(X,A)}(W;\Q)$$
for every neighborhood $W$ of $\ov\M_{g,k}(X,A;J)$ in $\X_{g,k}(X,A)$.
Furthermore, if 
$$\io_{W',W}\!: W\lra W'$$
is the inclusion map of a neighborhood $W$ into a larger neighborhood $W'$, then
$$\io_{W',W*}\, \Om_W=\Om_{W'}.$$
Thus, the above construction defines VFC for $\ov\M_{g,k}(X,A;J)$ as an element
of the inverse limit of the homology groups $H_*(W;\Q)$ under inclusion, 
taken over all neighborhoods of $\ov\M_{g,k}(X,A;J)$ in $\X_{g,k}(X,A)$.
If $(X,J)$ is algebraic, 
$\ov\M_{g,k}(X,A;J)$ is a deformation retract of a neighborhood~$W$,
and one can then define VFC for $\ov\M_{g,k}(X,A;J)$
as a homology class in such a neighborhood~$W$.
However, these formalities are not essential for defining GW-invariants as intersection
numbers of $\ov\M_{g,k}(X,A;J,\nu)$ with certain natural classes on $\X_{g,k}(X,A)$.\\

\noindent 
For a small {\it generic} perturbation $\nu$ of $\bar{\partial}_J$, 
the closure~of $\M_{g,k}^0(X,A;J,\nu)$
is the entire moduli space $\ov\M_{g,k}(X,A;J,\nu)$.
In particular, the results of~\cite{g1comp}, that are summarized in Subsection~\ref{back_subs},
cannot possibly generalize to $\M_{g,k}^0(X,A;J,\nu)$, even with $g\!=\!1$, for a generic~$\nu$.
Instead, for $g\!=\!1$, we consider {\it non-generic} perturbations~$\nu$
of $\bar{\partial}_J$, which we now describe.\\

\noindent
An element $[\Si;u]$ of $\X_{1,k}(X,A)$ is an equivalence class
of pairs consisting of a prestable genus-one Riemann surface $\Si$
and a smooth map $u\!:\Si\!\lra\!X$.
The prestable surface~$\Si$ is a union of the principal component(s)~$\Si_P$,
which is either a smooth torus or a circle of spheres, and
trees of rational bubble components, which together will be denoted by~$\Si_B$.
Let
$$\X_{1,k}^{\{0\}}(X,A)=
\big\{[\Si;u]\!\in\!\X_{1,k}(X,A)\!: u_*[\Si_P]\neq0\in H_2(X;\Z)\big\}.$$\\

\noindent
Suppose 
\begin{equation}\label{compsubs_e1}
[\Si;u]\in\X_{1,k}(X,A)-\X_{1,k}^{\{0\}}(X,A),
\end{equation}
i.e.~the degree of $u|_{\Si_P}$ is zero.
Let $\chi^0(\Si;u)$ be the set of components $\Si_i$ of $\Si$ such that 
for every bubble component $\Si_h$ that lies between $\Si_i$ and $\Si_P$,
including $\Si_i$ itself, the degree of $u|_{\Si_h}$ is zero.
The set $\chi^0(\Si;u)$ includes the principal component(s) of~$\Si$.
We give an example of the set  $\chi^0(\Si;u)$ in Figure~\ref{chi_fig}.
In this figure, as in Figure~\ref{m3_fig},
we show the domain~$\Si$ of the stable map $(\Si;u)$
and shade the components of the domain on which the degree of the map~$u$ is not zero.
Let
$$\Si_u^0=\bigcup_{i\in\chi^0(\Si;u)}\!\!\!\!\!\!\Si_i.$$

\begin{dfn}
\label{pert_dfn}
Suppose $(X,\om)$ is a compact symplectic manifold, $\under{J}\!\equiv\!(J_t)_{t\in[0,1]}$ 
is a continuous family of $\om$-tamed almost structures on~$X$, $A\!\in\!H_2(X;\Z)$,
and $k\!\in\!\bar\Z^+$.
A continuous family of multisections $\under{\nu}\!\equiv\!(\nu_t)_{t\in[0,1]}$,
with $\nu_t\!\in\!\G_{1,k}^{0,1}(X,A;J_t)$ for all $t\!\in\![0,1]$, 
is  \sf{effectively supported} if for every element 
$$b\!\equiv\![\Si,u] \in \X_{1,k}(X,A)\!-\!\X_{1,k}^{\{0\}}(X,A)$$
there exists a neighborhood $\W_b$ of $\Si_u^0$ in a semi-universal family of
deformations for~$b$ such that
$$\nu_t(\Si';u')\big|_{\Si'\cap\W_b}=0 \qquad
\forall~~  [\Si';u']\in\X_{1,k}(X,A),~t\!\in\![0,1].$$\\
\end{dfn}

\noindent
We use the $C^1$-topology on the space of all almost complex structures on $X$, 
in this definition and throughout the rest of the paper.
The bundles $\Ga_{1,k}^{0,1}(X,A;J_t)$ are contained in the bundle
$\Ga_{1,k}^1(X,A)$ over $\X_{1,k}(X,A)$
with the fibers
$$\Ga_{1,k}^1(X,A)\big|_{[b]}=\Ga^1(b)\big/\Aut(b),
\quad\hbox{where}\quad
\Ga^1(b)=L^p\big(\Si;T^*\Si\!\otimes_{\R}\!u^*TX\big),$$
and with the topology constructed as for $\Ga_{g,k}^{0,1}(X,A;J)$.
Finally, let $b\!=\![\Si;u]$ be an element of $\X_{1,k}(X,A)$.
A \sf{semi-universal universal family of deformations} for~$b$ is 
a fibration
$$\si_b\!: \ti\U_b\lra \De_b$$
such that $\De_b/\Aut(b)$ is a neighborhood of $b$ in $\X_{1,k}(X,A)$
and the fiber of $\si_b$ over a point $[\Si';u']$ is~$\Si'$.\\

\noindent
If $\under\nu$ is effectively supported and $[\Si;u]$ is as in~\e_ref{compsubs_e1}, 
then the restriction of $\nu_t$ to a neighborhood $W$ of $\Si_u^0$ in $\Si$ is zero
for all~$t$.
Furthermore, if $\{[\Si_k;u_k]\}$ is a sequence converging to $[\Si;u]$ in 
$\X_{1,k}(X,A)$, then for all $k$ sufficiently large and a choice of representatives
$(\Si_k;u_k)$ there is an open subset $W_k$ of $\Si_k$ such that 
$\nu_t(\Si_k;u_k)|_{W_k}\!=\!0$ and the open sets $W_k$ converge to~$W$;
see the beginning of Section~3 in~\cite{LT} for a detailed setting.\\

\noindent
For example, if $[\Si;u]$ is as indicated in Figure~\ref{chi_fig} and
$\under\nu$ is effectively supported, then $\nu_t(\Si;u)$ vanishes on a neighborhood
of $\Si_P\!\cup\!\Si_{h_3}$ in $\Si$.
On the other hand, even if $\Si_{h_2}$ had not been shaded, 
i.e.~the degree of $u|_{\Si_{h_2}}$ were zero, 
there still would have been no condition on $\nu_t(\Si;u)|_{\Si_{h_2}}$
because the degree of $u|_{\Si_{h_1}}$ is not~zero.\\

\noindent
If $\under{J}\!\equiv\!(J_t)_{t\in[0,1]}$ is a continuous family of $\om$-tamed almost 
structures on~$X$, we denote the space of effectively supported families
$\under{\nu}$ as in Definition~\ref{pert_dfn} by $\G^{\es}_{1,k}(X,A;\under{J})$.
Similarly, if $J$ is an almost complex structure on $X$,
we denote by $\G^{\es}_{1,k}(X,A;J)$ the subspace of elements $\nu$ of $\G^{0,1}_{1,k}(X,A;J)$ 
such that the family $\nu_t\!=\!\nu$ is effectively supported.\\

\noindent
{\it Remark:} Since $\nu$ is a multisection of $\Ga^{0,1}_{g,k}(X,A;J)$,
which is a union of orbi-vector spaces 
$$\Ga_{g,k}^{0,1}(X,A;J)\big|_{[b]}=\Ga^{0,1}(b;J)\big/\Aut(b),$$
$\nu$ is a family of equivalence classes of elements of $\Ga_{g,k}^{0,1}(X,A;J)$
and can be locally represented by a family of elements of $\Ga^{0,1}(\cdot;J)$.
In order to simplify notation, we will use the same symbol for both,
as the exact meaning will be determined by the context.

\subsection{Main Results}
\label{res_subs}

\noindent
In this subsection we state the main results of this paper.
We begin by describing the subspace $\ov\M_{1,k}^0(X,A;J,\nu)$ of $\ov\M_{1,k}(X,A;J,\nu)$.
We then state the main compactness result, i.e.~Theorem~\ref{comp_thm}.
One of its consequences is that for a small generic choice of~$\nu$
the moduli space $\ov\M_{1,k}^0(X,A;J,\nu)$ determines a virtual fundamental
class for $\ov\M_{1,k}^0(X,A;J)$, which is independent of~$J$;
see Theorem~\ref{reg_thm} and Corollary~\ref{fundclass_crl}.\\

\noindent
Suppose $[\Si;u]$ is an element of $\X_{1,k}(X,A)$. 
Every bubble component $\Si_i\!\subset\!\Si_B$ is a sphere and has 
a distinguished singular point, which will be called the {\it attaching node of~$\Si_i$}.
This is the node of $\Si_i$ that lies either on~$\Si_P$
or on a bubble $\Si_h$ that lies between $\Si_i$ and~$\Si_P$.
For example, if $\Si$ is as shown in Figure~\ref{chi_fig},
the attaching node of $\Si_{h_3}$ is the node $\Si_{h_3}$ shares with the torus.
If $[\Si;u]$ is as in~\e_ref{compsubs_e1},
we denote by $\chi(\Si;u)$ the set of bubble components $\Si_i$ such that
the attaching node of $\Si_i$ lies on $\Si_u^0$ and the degree of $u|_{\Si_i}$ is not zero,
i.e.~$\Si_i$ is not an element of $\chi^0(\Si;u)$; see Figure~\ref{chi_fig}.
These components are called {\it first-level $(\Si;u)$-effective} 
in \cite[Subsection~1.2]{g1comp}.\\

\noindent
Suppose $\nu\!\in\!\G^{\es}_{1,k}(X,A;J)$ and $[\Si;u]$ is an element
of $\ov\M_{1,k}(X,A;J,\nu)$ as in~\e_ref{compsubs_e1}.
Since $\Si_i\!\subset\!\Si_B$ is a sphere, we can represent every element of
$\X_{1,k}(X,A)$ by a pair $(\Si;u)$ such that
the attaching node of every bubble component $\Si_i\!\subset\!\Si_B$
is the south pole, or the point $\i\!=\!(0,0,-1)$, of $S^2\!\subset\!\Bbb{R}^3$.
Let $e_{\i}\!=\!(1,0,0)$ be a nonzero tangent vector to $S^2$ at the south pole.
If $i\!\in\!\chi(\Si;u)$, we put
$$\cD_i(\Si;u)= d\big\{u|_{\Si_i}\big\}\big|_{\i}e_{\i}
\in T_{u|_{\Si_i}(\i)}X.$$
Since $u|_{\Si_i}$ is $J$-holomorphic on a neighborhood of $\i$ in~$\Si_i$, 
the linear subspace $\C\!\cdot\!\cD_i(\Si;u)$ is determined 
by~$[\Si;u]$, just as in the $\nu\!=\!0$ case, which is considered
 \cite[Subsection~1.2]{g1comp}.
We also note that $u|_{\Si_u^0}$ is a degree-zero holomorphic map and thus constant.
Thus, $u$  maps the attaching nodes of all elements of $\chi(\Si;u)$
to the same point in~$X$,
just as in the $\nu\!=\!0$ case of \cite[Subsection~1.2]{g1comp}.

\begin{figure}
\begin{pspicture}(-1.1,-1.8)(10,1.25)
\psset{unit=.4cm}
\psellipse(8,-1.5)(1.5,2.5)
\psarc[linewidth=.05](6.2,-1.5){2}{-30}{30}\psarc[linewidth=.05](9.8,-1.5){2}{150}{210}
\pscircle[fillstyle=solid,fillcolor=gray](5.5,-1.5){1}\pscircle*(6.5,-1.5){.2}
\pscircle[fillstyle=solid,fillcolor=gray](3.5,-1.5){1}\pscircle*(4.5,-1.5){.2}
\pscircle(10.5,-1.5){1}\pscircle*(9.5,-1.5){.2}
\pscircle[fillstyle=solid,fillcolor=gray](11.91,-.09){1}\pscircle*(11.21,-.79){.2}
\pscircle[fillstyle=solid,fillcolor=gray](11.91,-2.91){1}\pscircle*(11.21,-2.21){.2}
\rput(8.5,1.5){$h_0$}\rput(5.5,0){$h_1$}\rput(3.5,0){$h_2$}
\rput(10.3,0){$h_3$}\rput(13.5,0.1){$h_4$}\rput(13.5,-2.9){$h_5$}
\rput(8,-5){\small ``tacnode"}
\pnode(8,-5){A1}\pnode(6.5,-1.5){B1}
\ncarc[nodesep=.35,arcangleA=-25,arcangleB=-15,ncurv=1]{->}{A1}{B1}
\pnode(8,-4.65){A2}\pnode(10.3,-1.5){B2}
\ncarc[nodesep=0,arcangleA=40,arcangleB=30,ncurv=1]{-}{A2}{B2}
\pnode(11,-.95){B2a}\pnode(11.02,-2.02){B2b}
\ncarc[nodesep=0,arcangleA=0,arcangleB=10,ncurv=1]{->}{B2}{B2a}
\ncarc[nodesep=0,arcangleA=0,arcangleB=10,ncurv=1]{->}{B2}{B2b}
\rput(25,-1.5){\begin{tabular}{l}$\chi^0(\Si;u)\!=\!\{h_0,h_3\}$\\
${}$\\ $\chi(\Si;u)\!=\!\{h_1,h_4,h_5\}$\end{tabular}}
\end{pspicture}
\caption{An illustration of Definition~\ref{degen_dfn}} 
\label{chi_fig}
\end{figure}

\begin{dfn}
\label{degen_dfn}
Suppose $(X,\om,J)$ is a compact almost Kahler manifold, $A\!\in\!H_2(X;\Z)^*$, and 
\hbox{$k\!\in\!\bar{\Z}^+$}.
If $\nu\!\in\!\G^{\es}_{1,k}(X,A;J)$ is an effectively supported perturbation
of the $\bar{\partial}_J$-operator,
the \sf{main component} of the space 
$\ov\M_{1,k}(X,A;J,\nu)$ is the subset $\ov\M_{1,k}^0(X,A;J,\nu)$ 
consisting of the elements $[\Si;u]$ of $\ov\M_{1,k}(X,A;J,\nu)$
such~that\\
${}\quad$ (a) the degree of $u|_{\Si_P}$ is not zero, or\\
${}\quad$ (b) the degree of $u|_{\Si_P}$ is zero and 
$\dim_{\Bbb{C}}\text{Span}_{(\Bbb{C},J)}\{\cD_i(\Si;u)\!:
i\!\in\!\chi(\Si;u)\}<|\chi(\Si;u)|$.
\end{dfn}

\noindent
This definition generalizes  \cite[Definition~1.1]{g1comp}.
As in~\cite{g1comp}, we let
$$H_2(X;\Z)^*= H_2(X;\Z)-\{0\}.$$
If $[\Si;u]$ is as in~\e_ref{compsubs_e1}, 
$[\Si;u]$ belongs to $\ov\M_{1,k}^0(X,A;J,\nu)$ if and only if
the branches of $u(\Si)$ corresponding to the attaching nodes
of the first-level effective bubbles of $[\Si;u]$ form a {\it generalized tacnode}.
In the case of Figure~\ref{chi_fig}, this means that the complex dimension of
the span of the images of $du$ at the attaching nodes 
of the bubbles $h_1$, $h_4$, and $h_5$ is at most two.

\begin{thm}
\label{comp_thm}
Suppose $(X,\om)$ is a compact symplectic manifold, $\under{J}\!\equiv\!(J_t)_{t\in[0,1]}$
is a continuous family of $\om$-tamed almost complex structures on $X$,
$A\!\in\!H_2(X;\Z)^*$, and $k\!\in\!\bar{\Z}^+$.
If 
$$\under{\nu}\!\equiv\!(\nu_t)_{t\in[0,1]} \in \G^{\es}_{1,k}(X,A;\under{J})$$
is a family of sufficiently small perturbations of the $\bar{\partial}_{J_t}$-operators
on $\X_{1,k}(X,A)$, then
$$\ov\M_{1,k}^0(X,A;\under{J},\under\nu)\equiv\!
\bigcup_{t\in[0,1]}\!\!\ov\M_{1,k}^0(X,A;J_t,\nu_t)$$
is compact.
\end{thm}

\noindent
The requirement that $\nu_t$ be sufficiently small means that it lies in 
a neighborhood of the zero section with respect to a $C^1$-type of topology,
with appropriate interpretations of the rate of change in the normal directions 
to the boundary strata of~$\X_{1,k}(X,A)$.
This topology will be made apparent in the proof.\\

\noindent
Theorem~\ref{comp_thm} follows immediately from 
Propositions~\ref{cuspmap_prp}-\ref{torus_prp};
see also the beginning of Subsection~\ref{str_subs}.
These propositions generalize 
\cite[Propositions 5.1-5.3]{g1comp}.

\begin{thm}
\label{reg_thm}
Suppose $(X,\om,J)$ is a compact almost Kahler manifold,
$A\!\in\!H_2(X;\Z)^*$, $k\!\in\!\bar\Z^+$, and $W$ is a neighborhood of 
$\ov\M_{1,k}^0(X,A;J)$ in $\X_{1,k}(X,A)$.
If $\nu\!\in\!\G_{1,k}^{\es}(X,A;J)$ is a sufficiently small generic perturbation 
of the $\bar{\partial}_J$-operator on $\X_{1,k}(X,A)$, 
then $\ov\M_{1,k}^0(X,A;J,\nu)$ determines a rational homology class in~$W$.
Furthermore, if $\under{J}\!=\!(J_t)_{t\in[0,1]}$ is a family of $\om$-tamed
almost complex structures on $X$, such that $J_0\!=\!J$ and $J_t$ is sufficiently close
to~$J$ for all~$t$, and $\nu_0$ and $\nu_1$ are sufficiently small generic perturbations
of $\bar{\partial}_{J_0}$ and $\bar{\partial}_{J_1}$ on $\X_{1,k}(X,A)$, 
then there exists  a homotopy 
$$\under\nu\!=\!(\nu_t)_{t\in[0,1]} \in \G^{\es}_{1,k}(X,A;\under{J})$$
between $\nu_0$ and $\nu_1$ 
such that $\ov\M_{1,k}^0(X,A;\under{J},\under\nu)$ determines a chain in $W$~and
$$\partial\ov\M_{1,k}^0(X,A;\under{J},\under\nu)
=\ov\M_{1,k}^0(X,A;J_1,\nu_1)-\ov\M_{1,k}^0(X,A;J_0,\nu_0).$$
\end{thm}

\begin{crl}
\label{fundclass_crl}
If $(X,\om,J)$ is a compact almost Kahler manifold,
$A\!\in\!H_2(X;\Z)^*$, and $k\!\in\!\bar\Z^+$, 
the moduli space  $\ov\M_{1,k}^0(X,A;J)$ carries 
a well-defined virtual fundamental class of the expected dimension.
This class is an invariant of~$(X,\om)$.
\end{crl}

\noindent
It is straightforward to see that for a generic $\nu\!\in\!\G_{1,k}^{\es}(X,A;J)$
the space $\ov\M_{1,k}(X,A;J,\nu)$ is stratified by smooth orbifolds of even dimensions.
The strata of 
$$\M_{1,k}^{\{0\}}(X,A;J,\nu) \equiv 
\ov\M_{1,k}(X,A;J,\nu)\cap  \X_{1,k}^{\{0\}}(X,A)$$ 
have the expected dimension,
based on the index of a certain elliptic operator.
In particular, the dimension of the main stratum of $\M_{1,k}^{\{0\}}(X,A;J,\nu)$ 
is $\dim_{1,k}(X,A)$, while the dimensions of all other strata of 
$\M_{1,k}^{\{0\}}(X,A;J,\nu)$ are smaller than $\dim_{1,k}(X,A)$.\\

\noindent
On the other hand, suppose $\U_{\T,\nu}(X;J)$ is a stratum of the complement of
$\M_{1,k}^{\{0\}}(X,A;J,\nu)$ in $\ov\M_{1,k}(X,A;J,\nu)$;
see Subsection~\ref{notation1_subs} for more details.
The sets $\chi^0(\Si;u)$ and $\chi(\Si;u)$ are independent of the choice
of  $[\Si;u]$ in~$\U_{\T,\nu}(X;J)$.
We denote them by $\chi^0(\T)$ and $\chi(\T)$, respectively.
By Definition~\ref{degen_dfn}, for every
$[\Si,u]\!\in\!\U_{\T,\nu}(X;J)$ and $i\!\in\!\chi^0(\T)$,
$u|_{\Si_i}$ is constant.
Thus,
$$\U_{\T,\nu}(X;J)\subset \cM_{\T}\times\X_{\bar{\T}}(X),$$
where $\cM_{\T}$ is a product of $|\chi^0(\T)|$ moduli spaces of smooth genus-zero and genus-one
curves and $\X_{\bar{\T}}(X)$ is a certain collection of $|\chi(\T)|$-tuples
of stable smooth genus-zero bubble maps. 
For example, if the elements of $\U_{\T,\nu}(X;J)$ are described by Figure~\ref{chi_fig},
$$\cM_{\T} = \cM_{1,2}\times \cM_{0,3}.$$
In this case, $\X_{\bar{\T}}(X)$ consists of triples of stable genus-zero bubble 
maps each with a special marked point, corresponding to the attaching nodes of
the elements of $\chi(\T)$, such that the values of three maps at the special marked 
points are the same.
If $\nu\!\in\!\G_{1,k}^{\es}(X,A;J)$ is generic, we have a fiber bundle
$$\pi_0\!:\U_{\T,\nu}(X;J)\lra\cM_{\T},$$
with fibers of the expected dimension.
An index computation then shows that
\begin{equation}\label{bddim_e1}
\dim \U_{\T,\nu}(X;J)\le\dim_{1,k}(X,A)+2\big(n\!-\!|\chi(\T)|),
\end{equation}
where $2n$ is the dimension of $X$ as before.\\

\noindent
We denote by $E_{\T}\!\lra\!\U_{\T,\nu}(X;J)$
the direct sum of the $|\chi(\T)|$ universal tangent line bundles 
for the special marked points of 
the elements of each $\chi(\T)$-tuple in $\X_{\bar{\T}}(X)$ and~by 
$$\ev_P\!:\U_{\T,\nu}(X;J)\lra X$$
the map sending an element $[\Si;u]$ of $\U_{\T,\nu}(X;J)$ to the value
of $u$ on~$\Si_P$.
Let $\ga_{E_{\T}}\!\lra\!\bP E_{\T}$ be the tautological line bundle.
By Definition~\ref{degen_dfn},
$$\U_{\T,\nu}(X;J) \cap \ov\M_{1,k}^0(X,A;J,\nu) = \pi_{\T}(\cZ_{\T}),$$
where 
$$\pi_{\T}\!: \bP E_{\T}\lra \U_{\T,\nu}(X;J)$$
is the bundle projection map and $\cZ_{\T}$ is the zero set of the section
of the vector bundle
$$\ga_{E_{\T}}^*\otimes\ev_P^*TX\lra\bP E_{\T}$$
induced by the differentials $\cD_i$, with $i\!\in\!\chi(\T)$, defined above.
It is straightforward to see that this section is transverse to the zero set
if $\nu\!\in\!\G_{1,k}^{\es}(X,A;J)$ is generic.
Thus,
\begin{equation*}\begin{split}
\dim\, \U_{\T,\nu}(X;J) \cap \ov\M_{1,k}^0(X,A;J,\nu)
&\le \dim \cZ_{\T}\\
&=\dim \U_{\T,\nu}(X;J)+2\big(\rk_{\C}E_{\T}\!-\!1)
-2\,\rk_{\C}\big(\ga_{E_{\T}}^*\!\otimes\!\ev_P^*TX\big)\\
&\le \dim_{1,k}(X,A)-2,
\end{split}\end{equation*}
by~\e_ref{bddim_e1}.\\

\noindent
By the above, for a generic $\nu\!\in\!\G_{1,k}^{\es}(X,A;J)$,
$\ov\M_{1,k}^0(X,A;J,\nu)$ is stratified by smooth orbifolds of even dimensions,
such that the main stratum is of dimension $\dim_{1,k}(X,A)$,
while all other strata have smaller dimensions.
Thus, the first claim of Theorem~\ref{reg_thm} follows from Theorem~\ref{comp_thm}
by the same topological construction as in Subsection~\ref{pertmaps_subs}.
The second claim of Theorem~\ref{reg_thm} is obtained by a similar argument.\\

\noindent
By the first claim of Theorem~\ref{reg_thm}, we can define a homology class
for $\ov\M_{1,k}^0(X,A;J)$, which is induced by $\ov\M_{1,k}^0(X,A;J,\nu)$, for any~$J$. 
By the last statement of Theorem~\ref{reg_thm}, this class is independent of the choice~$\nu$
and does not change under small changes in~$J$.
Since the space of $\om$-tamed almost complex structures on $X$ is path-connected,
it follows that the virtual fundamental class of $\ov\M_{1,k}^0(X,A;J)$ is an invariant of
$(X,\om)$.\\

\noindent
{\it Remark:} It is simplest to view the last statement above as 
the independence of all numbers $\GW_{1,k}^{0;X}(\mu)$ obtained by evaluating
natural cohomology classes on $\ov\M_{1,k}^0(X,A;J,\nu)$.

\section{Proof of Theorem~\ref{comp_thm}}
\label{gluing_sec}

\subsection{Notation: Genus-Zero Maps}
\label{notation0_subs}

\noindent
We now describe our notation for bubble maps from genus-zero Riemann surfaces and
for the spaces of such bubble maps that form
the standard stratifications of moduli spaces of stable maps.
We also state analogues of Definition~\ref{pert_dfn} for genus-zero maps
with one and two special marked points.\\

\noindent
In general, moduli spaces of stable maps can stratified by the dual graph.
However, in the present situation, it is more convenient to make use
of {\it linearly ordered sets}:

\begin{dfn}
\label{index_set_dfn1}
(1) A finite nonempty partially ordered set $I$ is a 
\sf{linearly ordered set} if 
for all \hbox{$i_1,i_2,h\!\in\!I$} such that $i_1,i_2\!<\!h$, 
either $i_1\!\le\!i_2$ or $i_2\!\le\!i_1$.\\
(2) A linearly ordered set $I$ is a \sf{rooted tree} if
$I$ has a unique minimal element, 
i.e.~there exists \hbox{$\hat{0}\!\in\!I$} such that $\hat{0}\!\le\!i$ 
for {all $i\!\in\!I$}.
\end{dfn}

\noindent
If $I$ is a linearly ordered set, let $\hat{I}$ be 
the subset of the non-minimal elements of~$I$.
For every $h\!\in\!\hat{I}$,  denote by $\io_h\!\in\!I$
the largest element of $I$ which is smaller than~$h$, i.e.
$\io_h\!=\!\max\big\{i\!\in\!I:i\!<\!h\big\}$.\\

\noindent
We identify $\C$ with $S^2\!-\!\{\i\}$ via 
the stereographic projection mapping the origin in $\C$ 
to the north pole, or the point $(0,0,1)$, in~$S^2$.
If $M$ is a finite set,
a \sf{genus-zero $X$-valued bubble map with $M$-marked points} is a~tuple
$$b=\big(M,I;x,(j,y),u\big),$$
where $I$ is a rooted tree, and
\begin{equation}\label{stablemap_e1}
x\!:\hat{I}\!\lra\!\C\!=\!S^2\!-\!\{\i\},\quad  j\!:M\!\lra\!I,\quad
y\!:M\!\lra\!\C,        \hbox{~~~and~~~} 
u\!:I\!\lra\!C^{\i}(S^2;X)
\end{equation}
are maps such that $u_h(\i)\!=\!u_{\io_h}(x_h)$ for all $h\!\in\!\hat{I}$.
We associate such a tuple with Riemann surface
\begin{equation}\label{stablemap_e2}
\Si_b=
\Big(\bigsqcup_{i\in I}\Si_{b,i}\Big)\Big/\!\sim,
\hbox{~~where}\qquad \Si_{b,i}=\{i\}\!\times\!S^2
\quad\hbox{and}\quad
(h,\i)\sim (\io_h,x_h)
~~\forall h\!\in\!\hat{I},
\end{equation}
with marked points 
$$y_l(b)\!\equiv\!(j_l,y_l)\in\Si_{b,j_l} \qquad\hbox{and}\qquad
y_0(b)\!\equiv\!(\hat{0},\i)\in\Si_{b,\hat{0}},$$
and continuous map $u_b\!:\Si_b\!\lra\!X$,
given by $u_b|_{\Si_{b,i}}\!=\!u_i$ for \hbox{all $i\!\in\!I$}.
The general structure of bubble maps is described
by tuples ${\cal T}\!=\!(M,I;j,\under{A})$, where 
$$A_i=u_{i*}[S^2]\in H_2(X;\Z) \qquad\forall\, i\!\in\!I.$$ 
We call such tuples \sf{bubble types}.
Let $\X_{\T}(X)$ denote the subset of $\X_{0,\{\hat{0}\}\sqcup M}(X,A)$ 
consisting of stable maps $[{\cal C};u]$ such that
$$[{\cal C};u]=
\big[(\Si_b,(\hat{0},\i),(j_l,y_l)_{l\in M});u_b\big],$$
for some bubble map $b$ of type ${\cal T}$ as above,
where $\hat{0}$ is the minimal element of~$I$; 
see \cite[Section~2]{gluing} for details.
For $l\!\in\!\{0\}\!\sqcup\!M$, let 
$$\ev_l\!:\X_{\T}(X)\lra X$$ 
be the evaluation map corresponding to the marked point $y_l$.\\

\noindent
With notation as above, suppose
$$[b]\!\equiv\!\big[M,I;x,(j,y),u\big] \in \X_{0,\{0\}\sqcup M}(X,A).$$
Let $\chi^0(b)$ be the set of components $\Si_{b,i}$ of $\Si_b$ such that 
for every component $\Si_{b,h}$ that lies between $\Si_i$ and  $\Si_{b,\hat{0}}$,
including $\Si_{b,i}$ and $\Si_{b,\hat{0}}$, the degree of $u|_{\Si_{b,h}}$ is zero.
For example, if $b$ is as indicated by Figure~\ref{derivest_fig} on page~\pageref{derivest_fig},
the set $\chi^0(b)$ consists of the two components that are not shaded.
The set $\chi^0(b)$ is empty if and only if the degree of the restriction of $u_b$ 
to the component containing the special marked point is not zero.
Let
$$\Si_b^0=\big\{(\hat{0},\i)\big\} \cup 
\bigcup_{i\in\chi^0(b)}\!\!\!\Si_{b,i}.$$
We denote by $\chi(b)$ the set of components $\Si_{b,i}$ such that
the attaching node of $\Si_{b,i}$ lies on $\Si_b^0$ and 
the degree of $u_b|_{\Si_{b,i}}$ is not zero,
i.e.~$\Si_{b,i}$ is not an element of $\chi^0(b)$.
If the degree of $u_b|_{\Si_{b,\hat{0}}}$ is not zero,
$\chi(b)\!=\!\{\hat{0}\}$.
If $A\!\neq\!0$ and  the degree of $u_b|_{\Si_{b,\hat{0}}}$ is zero,
the set $\chi(b)$ is not empty, but does not contain~$\hat{0}$.

\begin{dfn}
\label{pert0_dfn1}
Suppose $(X,\om)$ is a compact symplectic manifold, $\under{J}\!\equiv\!(J_t)_{t\in[0,1]}$ 
is a continuous family of $\om$-tamed almost structures on~$X$,
$A\!\in\!H_2(X;\Z)^*$, and $M$ is a finite set.
A continuous family of multisections $\under{\nu}\!\equiv\!(\nu_t)_{t\in[0,1]}$,
with $\nu_t\!\in\!\G_{0,\{0\}\sqcup M}^{0,1}(X,A;J_t)$ for all $t\!\in\![0,1]$, 
is  \sf{effectively supported} if for every element 
$b$ of $\X_{0,\{0\}\sqcup M}(X,A)$ 
there exists a neighborhood $\W_b$ of $\Si_b^0$ in a semi-universal family of
deformations for~$b$ such that
$$\nu_t(b')\big|_{\Si_{b'}\cap\W_b}=0 \qquad
\forall~~  [b']\in\X_{0,\{0\}\sqcup M}(X,A),~t\!\in\![0,1].$$
\end{dfn}

\begin{dfn}
\label{pert0_dfn2}
Suppose $(X,\om)$, $\under{J}\!\equiv\!(J_t)_{t\in[0,1]}$, $A$, and $M$
are as in Definition~\ref{pert0_dfn1}.
A continuous family of multisections $\under{\nu}\!\equiv\!(\nu_t)_{t\in[0,1]}$,
with $\nu_t\!\in\!\G_{0,\{0,1\}\sqcup M}^{0,1}(X,A;J_t)$ for all $t\!\in\![0,1]$, 
is  \sf{semi-effectively supported} if for every element 
$b$ of $\X_{0,\{0,1\}\sqcup M}(X,A)$ such that the marked point $y_1(b)$ lies on $\Si_b^0$
there exists a neighborhood $\W_b$ of $\Si_b^0$ in a semi-universal family of
deformations for~$b$ such that
$$\nu_t(b')\big|_{\Si_{b'}\cap\W_b}=0 \qquad
\forall~~  [b']\in\X_{0,\{0,1\}\sqcup M}(X,A),~t\!\in\![0,1].$$\\
\end{dfn}

\noindent
We denote the spaces of effectively and semi-effectively supported families 
$\under{\nu}$ as in Definitions~\ref{pert0_dfn1} and~\ref{pert0_dfn2} by 
$$\G^{\es}_{0,\{0\}\sqcup M}(X,A;\under{J})   \qquad\hbox{and}\qquad 
\G^{\ses}_{0,\{0,1\}\sqcup M}(X,A;\under{J}),$$ 
respectively.
Similarly to the genus-one case, if $J$ is an almost complex structure on $X$,
we denote~by 
$$\G^{\es}_{0,\{0\}\sqcup M}(X,A;J) \qquad\hbox{and}\qquad
\G^{\ses}_{0,\{0,1\}\sqcup M}(X,A;J)$$
the subspaces of elements $\nu$ of $\G^{0,1}_{0,\{0\}\sqcup M}(X,A;J)$ and
$\G^{0,1}_{0,\{0,1\}\sqcup M}(X,A;J)$ such that the family $\nu_t\!=\!\nu$ is 
effectively supported or semi-effectively supported, respectively.\\

\noindent
If $[b]\!=\![\Si_b;u_b]$ is an element of $\X_{0,\{0\}\sqcup M}(X,A)$ is as above 
and $i\!\in\!\chi(b)$, we~put
$$\cD_ib= d\big\{u_b|_{\Si_{b,i}}\big\}\big|_{\i}e_{\i}
\in T_{u_b|_{\Si_{b,i}}(\i)}X.$$
If  $\nu\!\in\!\G^{\es}_{0,\{0\}\sqcup M}(X,A;J)$ and $b$ is an element of 
$$\ov\M_{0,\{0\}\sqcup M}(X,A;J,\nu) \equiv
\big\{\bar\partial_J\!+\!\nu\big\}^{-1}(0),$$
then $u_b|_{\Si_{b,i}}$ is $J$-holomorphic on a neighborhood of $\i$ in~$\Si_{b,i}$ and
$\C\!\cdot_J\!\cD_ib$ is determined by~$b$, just as in Subsection~\ref{res_subs}.
This is also the case if $\nu\!\in\!\G^{\ses}_{0,\{0,1\}\sqcup M}(X,A;J)$
and $[b]$ is an element of $\ov\M_{0,\{0,1\}\sqcup M}(X,A;J,\nu)$
such that $y_1(b)\!\in\!\Si_b^0$.
In both of these cases, $u_b|_{\Si_b^0}$ is a degree-zero holomorphic map and thus constant.
Thus, $u_b$  maps the attaching nodes of all elements of $\chi(b)$
to the same point in~$X$, as in the genus-one case of Subsection~\ref{res_subs}.

\subsection{Notation: Genus-One Maps}
\label{notation1_subs}

\noindent
We next set up analogous notation for maps from genus-one Riemann surfaces.
In this case, we also need to specify the structure of the principal component.
Thus, we index the strata of $\X_{1,M}(X,A)$
by {\it enhanced linearly ordered sets}:

\begin{dfn}
\label{index_set_dfn2}
An \sf{enhanced linearly ordered set} is a pair $(I,\aleph)$,
where $I$ is a linearly ordered set, $\aleph$ is a subset of $I_0\!\times\!I_0$,
and $I_0$ is the subset of minimal elements of~$I$,
such that if $|I_0|\!>\!1$, 
$$\aleph=\big\{(i_1,i_2),(i_2,i_3),\ldots,(i_{n-1},i_n),(i_n,i_1)\big\}$$
for some bijection $i\!:\{1,\ldots,n\}\!\lra\!I_0$.
\end{dfn}

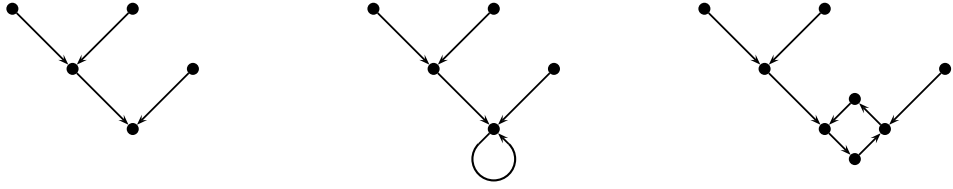
\begin{figure}
\begin{pspicture}(-1.1,-2)(10,1)
\psset{unit=.4cm}
\pscircle*(6,-3){.2}
\pscircle*(4,-1){.2}\psline[linewidth=.06]{->}(4.14,-1.14)(5.86,-2.86)
\pscircle*(8,-1){.2}\psline[linewidth=.06]{->}(7.86,-1.14)(6.14,-2.86)
\pscircle*(2,1){.2}\psline[linewidth=.06]{->}(2.14,.86)(3.86,-.86)
\pscircle*(6,1){.2}\psline[linewidth=.06]{->}(5.86,.86)(4.14,-.86)
% 2nd picture starts here
\pscircle*(18,-3){.2}\psline[linewidth=.06](17.86,-3.14)(17.5,-3.5)
\psarc(18,-4){.71}{135}{45}\psline[linewidth=.06]{->}(18.5,-3.5)(18.14,-3.14)
\pscircle*(16,-1){.2}\psline[linewidth=.06]{->}(16.14,-1.14)(17.86,-2.86)
\pscircle*(20,-1){.2}\psline[linewidth=.06]{->}(19.86,-1.14)(18.14,-2.86)
\pscircle*(14,1){.2}\psline[linewidth=.06]{->}(14.14,.86)(15.86,-.86)
\pscircle*(18,1){.2}\psline[linewidth=.06]{->}(17.86,.86)(16.14,-.86)
% 3rd picture starts here
\pscircle*(30,-2){.2}\pscircle*(30,-4){.2}\pscircle*(29,-3){.2}\pscircle*(31,-3){.2}
\psline[linewidth=.06]{->}(29.86,-2.14)(29.14,-2.86)
\psline[linewidth=.06]{->}(29.14,-3.14)(29.86,-3.86)
\psline[linewidth=.06]{->}(30.14,-3.86)(30.86,-3.14)
\psline[linewidth=.06]{->}(30.86,-2.86)(30.14,-2.14)
\pscircle*(27,-1){.2}\psline[linewidth=.06]{->}(27.14,-1.14)(28.86,-2.86)
\pscircle*(33,-1){.2}\psline[linewidth=.06]{->}(32.86,-1.14)(31.14,-2.86)
\pscircle*(25,1){.2}\psline[linewidth=.06]{->}(25.14,.86)(26.86,-.86)
\pscircle*(29,1){.2}\psline[linewidth=.06]{->}(28.86,.86)(27.14,-.86)
\end{pspicture}
\caption{Some enhanced linearly ordered sets}
\label{index_set_fig}
\end{figure}

\noindent
An enhanced linearly ordered set can be represented by an oriented connected graph.
In Figure~\ref{index_set_fig}, the dots denote the elements of~$I$.
The arrows outside the loop, if there are any, 
specify the partial ordering of the linearly ordered set~$I$.
In fact, every directed edge outside of the loop
connects a non-minimal element $h$ of $I$ with~$\io_h$.
Inside of the loop, there is a directed edge from $i_1$ to $i_2$
if and only if $(i_1,i_2)\!\in\!\aleph$.\\

\noindent
The subset $\aleph$ of $I_0\!\times\!I_0$ will be used to describe
the structure of the principal curve of the domain of stable maps in 
a stratum of~$\X_{1,M}(X,A)$.
If $\aleph\!=\!\eset$, and thus $|I_0|\!=\!1$,
the corresponding principal curve $\Si_P$ 
is a smooth torus, with some complex structure.
If $\aleph\!\neq\!\eset$, the principal components form a circle of spheres:
$$\Si_P=\Big(\bigsqcup_{i\in I_0}\{i\}\!\times\!S^2\Big)\Big/\sim,
\qquad\hbox{where}\qquad
(i_1,\i)\sim(i_2,0)~~\hbox{if}~~(i_1,i_2)\!\in\!\aleph.$$
A \sf{genus-one $X$-valued bubble map with $M$-marked points} is a tuple
$$b=\big(M,I,\aleph;S,x,(j,y),u\big),$$
where $S$ is a smooth Riemann surface of genus one if $\aleph\!=\!\eset$
and the circle of spheres $\Si_P$ otherwise.
The objects $x$, $j$, $y$, $u$, and $(\Si_b,u_b)$ are as in 
\e_ref{stablemap_e1} and \e_ref{stablemap_e2}, except 
the sphere $\Si_{b,\hat{0}}$ is replaced by the genus-one curve $\Si_{b;P}\!\equiv\!S$.
Furthermore, if $\aleph\!=\!\eset$, and thus $I_0\!=\!\{\hat{0}\}$ is a single-element set,
$u_{\hat{0}}\!\in\!C^{\i}(S;X)$ and $y_l\!\in\!S$ if $j_l\!=\!\hat{0}$.
In the genus-one case, the general structure of bubble maps is encoded by
the tuples of the form ${\cal T}\!=\!(M,I,\aleph;j,\under{A})$.
Similarly to the genus-zero case, we denote by $\X_{\T}(X)$
the subset of $\X_{1,M}(X,A)$ 
consisting of stable maps $[{\cal C};u]$ such that
$$[{\cal C};u]=\big[(\Si_b,(j_l,y_l)_{l\in M});u_b\big],$$
for some bubble map $b$ of type $\T$ as above.
If $\nu$ is an element of $\G_{1,M}^{\es}(X,A)$, we put
$$\U_{\T,\nu}(X;J)=\big\{[b]\!\in\!\X_{\T}(X)\!:
\{\bar{\partial}_J\!+\!\nu\}(b)=0\big\}.$$\\

\noindent
All vector orbi-bundles we encounter will be assumed to be normed.
Some will come with natural norms; for others, we choose a norm,
sometimes implicitly, once and for~all.
If \hbox{$\F\!\lra\!\X$} is a normed vector bundle
and $\de\!\!\in\!\R^+$, let
$$\F_{\de}=\big\{\ups\!\in\!\F\!:|\ups|\!<\!\de\big\}.$$
If $\Om$ is any subset of $\F$, we take 
$\Om_{\de}\!=\!\Om\cap\F_{\de}$.

\subsection{Outline of the Proof of Theorem~\ref{comp_thm}}
\label{str_subs}

\noindent
Suppose $(X,\om)$ is a compact symplectic manifold, $\under{J}\!\equiv\!(J_t)_{t\in[0,1]}$
is a continuous family of $\om$-tamed almost complex structures on $X$,
$A\!\in\!H_2(X;\Z)^*$, $M$ is a finite set, and 
$$\under\nu\!\equiv\!(\nu_t)_{t\in[0,1]} \in \G^{\es}_{1,M}(X,A;\under{J})$$
is a family of sufficiently small perturbations of the $\bpar_{J_t}$-operators
on $\X_{1,M}(X,A)$.
Let $t_r$ and $b_r$ be sequences of elements in $[0,1]$ and 
in $\ov\M_{1,M}^0(X,A;J_{t_r},\nu_{t_r})$ such that 
$$\lim_{r\lra\i}t_r=0  \quad\hbox{and}\quad
\lim_{r\lra\i}b_r=b\in\ov\M_{1,M}(X,A;J_0,\nu_0).$$
We need to show that $b\!\in\!\ov\M_{1,M}^0(X,A;J_0,\nu_0)$.
By Definition~\ref{degen_dfn}, it is sufficient to assume that $b$ is an
element of $\U_{\T,\nu_0}(X;J_0)$ for a bubble type
$$\T=\big(M,I,\aleph;j,\under{A})$$
such that $A_i\!=\!0$ for all minimal elements $i\!\in\!I$.\\

\noindent
We can also assume that for some bubble type
$$\T'=\big(M,I',\aleph';j',\under{A}')$$
$b_r\!\in\!\U_{\T',\nu_{t_r}}(X;J_{t_r})$ for all $r$.
If $A_i'\!=\!0$ for all minimal elements $i\!\in\!I'$,
the desired conclusion follows Proposition~\ref{cuspmap_prp} below,
as it implies that the second condition in Definition~\ref{degen_dfn} 
is closed with respect to the stable map topology.
If $A_i'\!\neq\!0$ for some minimal element $i\!\in\!I'$ and $\aleph'\!\neq\!\eset$,
i.e.~the principal component of $\Si_{b_r}$ is a circle of spheres,
Proposition~\ref{node_prp} implies that $b$ satisfies the second condition in
Definition~\ref{degen_dfn}.
Finally, if $\aleph'\!=\!\eset$ and $A_i'\!\neq\!0$ for the unique minimal element
$i$ of~$I'$,  the desired conclusion follows from Proposition~\ref{torus_prp}.
We note that the three propositions are applied with $b$ and $b_r$ 
that are components of the ones above.\\

\noindent
Let $[n]=\big\{1,\ldots,n\big\}$.

\begin{prp}
\label{cuspmap_prp}
Suppose $(X,\om)$ is a compact symplectic manifold, $\under{J}\!\equiv\!(J_t)_{t\in[0,1]}$
is a continuous family of $\om$-tamed almost complex structures on $X$,
$A\!\in\!H_2(X;\Z)^*$, $M$ is a finite set, and 
$$\under\nu\!\equiv\!(\nu_t)_{t\in[0,1]} \in 
\G^{\es}_{0,\{0\}\sqcup M}(X,A;\under{J})$$
is a family of sufficiently small perturbations of the $\bpar_{J_t}$-operators
on $\X_{0,\{0\}\sqcup M}(X,A)$.
If $t_r$ and $[b_r]$ are sequences of elements in $[0,1]$ and 
in $\M_{0,\{0\}\sqcup M}^0(X,A;J_{t_r},\nu_{t_r})$  such that
$$\lim_{r\lra\i}t_r=0 \quad\hbox{and}\quad
\lim_{r\lra\i}[b_r]=[b]\in \ov\M_{0,\{0\}\sqcup M}(X,A;J_0,\nu_0),$$
then either\\
${}\quad$ (a) $\dim_{\C}\hbox{Span}_{(\C,J_0)}\{\cD_ib\!:i\!\in\!\chi(b)\}<|\chi(b)|$, or\\
${}\quad$ (b) $\bigcap_{r=1}^{\i}\ov{\bigcup_{r'>r}\C\cdot_{J_{r'}}\!\cD_{\hat{0}}b_{r'}}
\subset\hbox{Span}_{(\C,J_0)}\{\cD_ib\!:i\!\in\!\chi(b)\}$.
\end{prp}

\begin{prp}
\label{node_prp}
Suppose $(X,\om)$ and $\under{J}$ are as in Proposition~\ref{cuspmap_prp},
$n\!\in\!\Z^+$, $A_1,\ldots,A_n\!\in\!H_2(X;\Z)^*$, $M_1,\ldots,M_n$ are finite sets, and
for each $k\!\in\![n]$
$$\under\nu_k\!\equiv\!(\nu_{k,t})_{t\in[0,1]} \in 
\G^{\ses}_{0,\{0,1\}\sqcup M_k}(X,A_k;\under{J})$$
is a family of sufficiently small perturbations of the $\bpar_{J_t}$-operators
on $\X_{0,\{0,1\}\sqcup M_k}(X,A)$.
Let $t_r$  and $[b_{k,r}]$ be sequences of elements in $[0,1]$ and
in $\M_{0,\{0,1\}\sqcup M_k}^0(X,A_k;J_{t_r},\nu_{k,t_r})$
for $k\!\in\![n]$  such~that
\begin{gather*}
\ev_1(b_{k,r})=\ev_0(b_{k+1,r}) \quad\forall\, k\!\in\![n\!-\!1], 
\qquad  \ev_1(b_{n,r})=\ev_0(b_{1,r}),\\
\lim_{r\lra\i}t_r=0, \quad\hbox{and}\quad
\lim_{r\lra\i}[b_{k,r}]=[b_k]\in
\ov\M_{0,\{0,1\}\sqcup M_k}(X,A_k;J_0,\nu_{k,0})
\quad\forall\, k\!\in\![n].
\end{gather*}
If $y_1(b_k)\!\in\!\Si_{b_k}^0$ for all $k\!\in\![n]$, then 
$$\dim_{\C}\hbox{Span}_{(\C,J_0)}\big\{\cD_ib_k\!:i\!\in\!\chi(b_k),~k\!\in\![n]\big\}
< \sum_{k=1}^{k=n}|\chi(b_k)|.$$
\end{prp}

\begin{prp}
\label{torus_prp}
Suppose $(X,\om)$, $\under{J}$, $A$, and $M$  are as in Proposition~\ref{cuspmap_prp} and 
$$\under\nu\!\equiv\!(\nu_t)_{t\in[0,1]} \in 
\G^{\es}_{1,M}(X,A;\under{J})$$
is a family of sufficiently small perturbations of the $\bpar_{J_t}$-operators
on $\X_{1,M}(X,A)$.
Let  $t_r$ and $[b_r]$ be sequences of elements in $[0,1]$ and 
in $\M_{1,M}^0(X,A;J_{t_r},\nu_{t_r})$  such that
$$\lim_{r\lra\i}t_r=0 \quad\hbox{and}\quad
\lim_{r\lra\i}[b_r]=[b]\in \ov\M_{1,M}(X,A;J_0,\nu_0).$$
If $b\!=\!(\Si;u)$ is such that the degree of $u|_{\Si_P}$ is zero, then
$$\dim_{\C}\hbox{Span}_{(\C,J_0)}\{\cD_ib\!:i\!\in\!\chi(b)\}<|\chi(b)|.$$\\
\end{prp}

\noindent
Propositions~\ref{cuspmap_prp}, \ref{node_prp}, and~\ref{torus_prp} follow immediately 
from the estimates \e_ref{cuspmap_e4}, \e_ref{node_e4}, and~\e_ref{torus_e10} below.
These estimates are obtained by combining the approach of \cite[Sections 3,4]{g1comp}
with some aspects of the local setting of \cite[Section~3]{LT}.
A key step is \cite[Lemma~3.5]{g1comp}
that gives power series expansions for the behavior of derivatives of 
$J$-holomorphic genus-zero maps under gluing.
They lead to estimates on obstructions to smoothing genus-one $J$-holomorphic
maps from singular domains in  \cite[Lemma~4.4]{g1comp}.
While the maps we encounter are not $J$-holomorphic on the entire domain,
they are $J$-holomorphic  around the part of the domain which is essential 
for the estimates of  \cite[Lemmas 3.5,4.4]{g1comp},
i.e.~$\Si_b^0$ in the notation of Subsections~\ref{pertmaps_subs} 
and~\ref{notation0_subs} above.
The argument in the next two subsections is in fact an extension of
\cite[Section~5]{g1comp}, but is far more detailed (as promised in~\cite{g1comp}).

\subsection{Proofs of Propositions~\ref{cuspmap_prp} and~\ref{node_prp}} 
\label{g0prp_subs}

\noindent
Let $(X,\om)$, $\under{J}$, $A$, $M$, $\under\nu$, 
$$b=(M,I;x,(j,y),u), \qquad\hbox{and}\qquad u_i\equiv u_b|_{\Si_{b,i}}$$ 
be as in the statement of Proposition~\ref{cuspmap_prp}.
For each $i\!\in\!I$, we~put
$$\Ga(b;i)=\big\{\xi\!\in\!L^p_1(\Si_{b,i};u_i^*TX\big)\!:
\xi(\i)\!=\!0\big\} \qquad\hbox{and}\qquad
\Ga^{0,1}(b;i)=
L^p(\Si_{b,i};\La^{0,1}_{J_0,j}T^*\Si_{b,i}\!\otimes\!u_i^*TX\big),$$
where $j$ is the complex structure on $\Si_b$.
We denote~by
$$D_{J_0,\nu_0;b,i}\!: \Ga(b;i) \lra \Ga^{0,1}(b;i)$$
the linear operator induced by the linearization $D_{J_0,\nu_0;b}$
of the section $\bar{\partial}_{J_0}\!+\!\nu_0$ at~$b$.\\

\noindent
We put
\begin{equation}\label{iplusdfn_e}
I^+=\big\{i\!\in\!I\!:A_i\!\neq\!0\big\}.
\end{equation}
For each $i\!\in\!I^+$, choose a finite-dimensional linear subspace 
$$\ti\Ga^{0,1}_-(b;i)\subset \Ga\big(\Si_{b,i}\!\times\!X;
\La^{0,1}_{J_0,j}\pi_1^*T^*\Si_{b,i}\!\otimes\!\pi_2^*TX\big)$$
such that
\begin{equation*}\begin{split}
\Ga^{0,1}(b;i) =\Im D_{J_0,\nu_0;b,i} \oplus
\big\{\{\id\!\times\!u_i\}^*\eta\!:\eta\!\in\!\ti\Ga^{0,1}_-(b;i)\big\}
\end{split}\end{equation*}
and every element of $\ti\Ga^{0,1}_-(b;i)$ vanishes on a neighborhood 
of $\i\!\in\!\Si_{b,i}$ and the nodes $x_h(b)\!\in\!\Si_{b,i}$ with $\io_h\!=\!i$.
If $i\!\in\!I\!-\!I^+$, we denote by $\ti\Ga^{0,1}_-(b;i)$ the zero vector space.
Let $\T$ be the bubble type of the map~$b$.
We put
\begin{equation*}\begin{split}
\ti\U= \big\{ b'\!\equiv\!(M,I;x',(j,y'),u')\!:~& [b']\!\in\!\X_{\T}(X),\\
& \pi_i\{\bar\partial_{J_0,j}\!+\!\nu_0\}u_{b'}\in
\{\id\!\times\!u_i'\}^*\ti\Ga^{0,1}_-(b;i)~\forall i\!\in\!I\big\},
\end{split}\end{equation*}
where
$$\pi_i\!: \Ga^{0,1}(b';J_0) \lra \Ga^{0,1}(b';i)$$
is the natural projection map.
By the Implicit Function Theorem, $\ti\U$ is a smooth manifold near~$b$.
Let
$$\ev_0\!: \ti\U\lra X,  \qquad b'\lra u_{b'}(\hat{0},\i),$$
be the evaluation map for the special marked point~$0$; 
see also Subsection~\ref{notation0_subs}.
Let
$$\wt{\cal F}\equiv \ti\U\times\C^{\hat{I}}$$
be the bundle of smoothing parameters.
We denote by $\wt{\cal F}^{\eset}$ the subset of  $\wt{\cal F}$ consisting
of the elements with all components nonzero.
For each $\ups\!=\!(b',v)$, where $b'\!\in\!\ti\U$
and $v\!=\!(v_i)_{i\in\hat{I}}$, and $i\!\in\!\chi(b)$, we~put
$$\rho_i(\ups)=\!\prod_{\hat{0}<h\le i}\!\!\!v_h\in\Bbb{C}
\qquad\hbox{and}\qquad
x_i(\ups)=\!\!\sum_{\hat{0}<i'\le i}\!\Big(x_{i'}(b')
\!\!\!\prod_{\hat{0}<h<i'}\!\!\!\!v_h\Big)~\in~\Bbb{C},$$
where $x_i(b')$ is the point of $\Si_{b',\io_i}$ to which the bubble $\Si_{b',i}$ is
attached; see \e_ref{stablemap_e2} and Figure~\ref{derivest_fig} 
on page~\pageref{derivest_fig}.\\

\noindent
For each sufficiently small element $\ups\!=\!(b',v)$ of $\wt{\cal F}^{\eset}$, let 
$$q_{\ups}\!:\Si_{\ups}\lra\Si_{b'}$$
be the basic gluing map constructed in  \cite[Subsection~2.2]{gluing}.
In this case, $\Si_{\ups}$ is the projective line $\Bbb{P}^1$ with $|M|\!+\!1$ marked points.
The map~$q_{\ups}$ collapses $|\hat{I}|$ circles on~$\Si_{\ups}$.
It induces a metric $g_{\ups}$ on~$\Si_{\ups}$ such that
$(\Si_{\ups},g_{\ups})$ is obtained from $\Si_{b'}$ by replacing 
the $|\hat{I}|$ nodes of $\Si_{b'}$ by thin necks.
Let
$$u_{\ups}=u_{b'}\circ q_{\ups}.$$
We fix a $J_0$-compatible metric $g$ on $X$ and denote 
the corresponding $J_0$-compatible connection by~$\na$.
The map~$q_{\ups}$ induces norms $\|\cdot\|_{\ups,p,1}$ and~$\|\cdot\|_{\ups,p}$
on the spaces 
$$\Ga(\Si_{\ups};u_{\ups}^*TX)     \qquad\hbox{and}\qquad
\Ga(\Si_{\ups};\La^{0,1}_{J_0,j}T^*\Si_{\ups}\!\otimes\!u_{\ups}^*TX),$$
respectively; see  \cite[Subsection~3.3]{gluing}.
We denote the corresponding completions by $\Ga(\ups)$ and~$\Ga^{0,1}(\ups)$.
The norms $\|\cdot\|_{\ups,p,1}$ and $\|\cdot\|_{\ups,p}$
are equivalent to the ones used in Section~3 of~\cite{LT}.\\

\noindent
Let $t_r$ and $b_r$ be as in Proposition~\ref{cuspmap_prp}.
Since the sequence $[b_r]$ converges to $[b]$, for all $r$ sufficiently large
there exist\\
$$b_r'\in\ti\U, \qquad
\ups_r=(b_r',v_r)\in\wt{\cal F}^{\eset}, \qquad\hbox{and}\qquad
\xi_r\in\Ga(\ups_r)$$
such that
\begin{gather}\label{cuspmap_e0}
\lim_{r\lra\i}b_r'=b, \qquad \lim_{r\lra\i}|v_r|=0, \qquad
\xi_r(\i)=0~~~\forall\, r,\qquad
\lim_{r\lra\i}\|\xi_r\|_{\ups_r,p,1}=0,\\
\hbox{and}\qquad b_r\!\equiv\!\big(\Si_{b_r};u_{b_r}\big)=
\big(\Si_{\ups_r};\exp_{u_{\ups_r}}\!\xi_r\big).\notag
\end{gather}
The last equality holds for a representative $b_r$ for $[b_r]$.\\

\noindent
{\it Remark:}
The existence of $b_r'$, $\ups_r$, and $\xi_r$ as above can be shown 
by an argument similar to the surjectivity argument in  \cite[Section~4]{gluing}, 
with significant simplifications.
In fact, the only facts about the bubble maps $b_r'$ we use below are that they
are constant on the degree-zero components and holomorphic on fixed neighborhoods
of the attaching nodes of the first-level effective bubbles.
Such bubble maps $b_r'$, along with $\ups_r$ and $\xi_r$, can be constructed directly from
the maps $b_r$; see the beginning of Subsection~4.4 in~\cite{gluing}.\\

\noindent
If $\de\!\in\!\Bbb{R}^+$, $b'\!\in\!\ti\U$, and 
$\ups\!=\!(b',v)\!\in\!\wt{\cal F}^{\eset}$ is sufficiently small, we put
\begin{gather}\label{mainpart_e1}
\Si_{b'}^0(\de)=\Si_{b'}^0 \cup \!\bigcup_{i\in\chi(b)}\!\!\! A_{b',i}(\de),
\quad\hbox{where}\quad
A_{b',i}(\de)=\big\{(i,z)\!:|z|\!\ge\!\de^{-1/2}/2\big\}\subset \Si_{b',i}\!\approx\!S^2,\\
\label{mainpart_e2}
\hbox{and}\qquad \Si_{\ups}^0(\de)=q_{\ups}^{-1}\big(\Si_{b'}^0(\de)\big).
\end{gather}
Choose $\de\!\in\!\Bbb{R}^+$ such that for all $i\!\in\!\chi(b)$
all elements of $\ti\Ga_-^{0,1}(b;i)$ vanish on $A_{b,i}(2\de)$
and for all $r$ sufficiently large 
$$\nu_t(b_r')\big|_{\Si_{b_r'}^0(2\de)}=0 \quad\hbox{and}\quad
\nu_t(b_r)\big|_{\Si_{\ups_r}^0(2\de)}=0 \qquad\forall~t\!\in\![0,1].$$
Such a positive number $\de$ exists by our assumptions on the spaces
$\ti\Ga_-^{0,1}(b;i)$ and the family of perturbations~$\under\nu$;
see Definition~\ref{pert0_dfn1}.\\

\noindent
For every element $b'\!=\!(\Si_{b'};u_{b'})$ of $\ti\U$
and every sufficiently small element $\ups\!=\!(b',v)$ of $\wt{\cal F}^{\eset}$,
we denote~by
$$\Hol_{J_0}\big(\Si_{b'}^0(\de);T_{\ev_0(b')}X\big)
\qquad\hbox{and}\qquad
\Hol_{J_0}\big(\Si_{\ups}^0(\de);T_{\ev_0(b')}X\big)$$
the spaces of holomorphic maps from $\Si_{b'}^0(\de)$ and $\Si_{\ups}^0(\de)$ 
into the complex vector space $(T_{\ev_0(b')}X;J_0)$.
Let $\exp$ be the $\na$-exponential map.
For every $b'\!\in\!\ti\U$ as above, $u_{b'}|_{\Si_{b'}^0}$ is constant 
and $u_{b'}|_{\Si_{b'}^0(2\de)}$ is $J_0$-holomorphic.
Thus, if $\de$ is sufficiently small, there exist continuous families of maps
$$\Phi_{b'}\in L^p_1\big(\Si_{b'}^0(\de);\End(T_{\ev_0(b')}X)\big)
\qquad\hbox{and}\qquad
\vt_{b'}\in\Hol_{J_0}\big(\Si_b^0(\de);T_{\ev_0(b')}X\big)$$
with $b'\!\in\!\ti\U$ such that for all $b'$ sufficiently close to~$b'$
\begin{gather*}
\Phi_{b'}|_{\Si_{b'}^0}=\Id, \qquad 
\big\|\Phi_{b'}\!-\!\Id\big\|_{b',p,1}\le\frac{1}{2}, \qquad\hbox{and}\\
\exp_{\ev_0(b')}\big(\Phi_{b'}(z)\vt_{b'}(z)\big)=u_{b'}(z) 
\qquad\forall z\in\Si_{b'}^0(\de).
\end{gather*}
This statement follows immediately from the proof of Theorem~2.2 in~\cite{FlHS}.
Similarly, for every
$$b_r\!\equiv\!\big(\Si_{b_r};u_{b_r}\big)=
\big(\Si_{\ups_r};\exp_{u_{\ups_r}}\!\xi_r\big)$$
with $r$ sufficiently large, $u_{b_r}|_{\Si_{\ups_r}^0(2\de)}$ is $J_{t_r}$-holomorphic.
Since $\|\xi_r\|_{\ups_r,p,1}$ tends to zero as $r$ approaches~$\i$,
if $\de$ is sufficiently small and $r$ is sufficiently large,
there exist
$$\Phi_{b_r}\in L^p_1\big(\Si_{\ups_r}^0(\de);\End(T_{\ev_0(b_r')}X)\big), \qquad
\vt_{b_r}\in\Hol_{J_0}\big(\Si_{\ups_r}^0(\de);T_{\ev_0(b_r')}X\big)$$
such that
\begin{gather*}
\big\|\Phi_{b_r}\!-\!\Phi_{b_r'}\!\circ\!q_{\ups_r}\big\|_{\ups_r,p,1}
\le  C\big(\|J_0\!-\!J_{t_r}\|_{C^1}\!+\!|\ups_r|^{1/p}\!+\!\|\xi_r\|_{\ups_r,p,1}\big)
\qquad\hbox{and}\\
\exp_{\ev_0(b_r')}\big(\Phi_{b_r}(z)\vt_{b_r}(z)\big) = u_{b_r}(z)
\qquad\forall z\in\!\Si_{\ups_r}^0(\de).
\end{gather*}
In the inequality above, both norms $\|\cdot\|_{\ups_r,p,1}$ are the norms
induced from the pregluing construction as in Subsection~3.3 of~\cite{gluing}.
With these norms, the existence of $\Phi_{b_r}$ and $\vt_{b_r}$ follows easily
from the proof of Theorem~2.2 in~\cite{FlHS}; see the paragraph following 
Lemma~3.3 in~\cite{g1comp}.\\

\noindent
If $i\!\in\!\chi(b)$ and $b'\!\in\!\ti\U$, let $w_i$ be
the standard holomorphic coordinate centered at the point $\i$ in $\Si_{b',i}\!=\!S^2$.
If $m\!\in\!\Z^+$, we~put
$$\cD_i^{(m)}\vt_{b'} = \frac{1}{m!}
\frac{d^m}{dw_i^m} \vt_{b',i}(w_i)\Big|_{w_i=0} \in T_{\ev_0(b')}X,
\qquad\hbox{where}\qquad \vt_{b',i}\!=\!\vt_{b'}|_{\Si_{b',i}}.$$
Similarly, for all $r$ sufficiently large, we put
$$\cD_{\hat{0}}^{(m)}\vt_{b_r} = \frac{1}{m!}
\frac{d^m}{dw^m} \vt_{b_r,i}(w)\Big|_{w=0} \in T_{\ev_0(b_r')}X,$$
where $w$ is the standard holomorphic coordinate centered 
at the point $\i$ in $\Si_{\ups_r}\!\approx\!S^2$.
The key step in the proof of Propositions~\ref{cuspmap_prp} and~\ref{node_prp} is
the power series expansion
\begin{equation}\label{derivexp_e1}
\cD_{\hat{0}}^{(m)}\vt_{b_r}=
\sum_{k=1}^{k=m} \binom{m\!-\!1}{k\!-\!1}\sum_{i\in\chi(b)} \!\!\!
x_i^{m-k}(\ups_r)\rho_i^k(\ups_r) \big\{ \cD_i^{(k)}\vt_{b_r'}\!+\!\ve_{i,r}^{(k)}\big\}
\in (T_{\ev_0(b_r')}X,J_0),
\end{equation}
for some $\ve_{i,r}^{(k)}\!\in\!T_{\ev_0(b_r')}X$ such that 
\begin{equation}\label{derivexp_e2}
\big|\ve_{i,r}^{(k)}\big| \le C\de^{-k/2}
\big(\|J_{t_r}\!-\!J_0\|_{C^1}\!+\!|\ups_r|^{1/p}\!+\!\|\xi_r\|_{\ups_r,p,1}\big).
\end{equation}
The expansion~\e_ref{derivexp_e1} is obtained by exactly the same integration-by-parts
argument as the expansion in~(2a) of Lemma~3.4 in~\cite{g1comp}.
We point out that $\ve_{i,r}^{(k)}$ is independent of~$m$.
The $m\!=\!1$ case of the estimates~\e_ref{derivexp_e1} 
and~\e_ref{derivexp_e2} is illustrated in Figure~\ref{derivest_fig}.\\

\begin{figure}
\begin{pspicture}(-1.1,-1.8)(10,1.25)
\psset{unit=.4cm}
\pscircle(5,-1.5){1.5}\pscircle*(5,-3){.24}\rput(5,-3.8){$(\hat{0},\i)$}
\pscircle[fillstyle=solid,fillcolor=gray](3.23,.27){1}\pscircle*(3.94,-.44){.2}
\pscircle[fillstyle=solid,fillcolor=gray](1.23,.27){1}\pscircle*(2.23,.27){.2}
\pscircle(6.77,.27){1}\pscircle*(6.06,-.44){.2}
\pscircle[fillstyle=solid,fillcolor=gray](6.77,2.27){1}\pscircle*(6.77,1.27){.2}
\pscircle[fillstyle=solid,fillcolor=gray](8.77,.27){1}\pscircle*(7.77,.27){.2}
\rput(3.3,1.8){$h_1$}\rput(1.3,1.8){$h_2$}\rput(5.3,1){$h_3$}
\rput(8.4,2.4){$h_4$}\rput(10.4,.4){$h_5$}
\pnode(2,-3){A1}\rput(.75,-3){$x_{h_1}(b_r')$}\pnode(4,-.46){B1}
\ncarc[nodesep=.35,arcangleA=-30,arcangleB=-50,ncurv=.8]{->}{A1}{B1}
\pnode(7.3,-3){A2}\rput(8.5,-2.9){$x_{h_3}(b_r')$}\pnode(6,-.45){B2}
\ncarc[nodesep=.35,arcangleA=30,arcangleB=40,ncurv=.8]{->}{A2}{B2}
\pnode(8.8,-1.55){A3}\rput(10,-1.4){$x_{h_5}(b_r')$}\pnode(7.85,.25){B3}
\ncarc[nodesep=.35,arcangleA=45,arcangleB=110,ncurv=1]{->}{A3}{B3}
\rput(24.9,0){\begin{small}\begin{tabular}{l}
$\chi(b)\!=\!\{h_1,h_4,h_5\}$\\
$\rho(\ups)=(v_{h_1},v_{h_3}v_{h_4},v_{h_3}v_{h_5})$\\
$x_{h_5}(\ups)=x_{h_3}(b_r')+v_{h_3}x_{h_5}(b_r')$\\
\\
$\cD_{\hat{0}}^{(1)}\vt_{b_r}
\cong v_{h_1}\big(\cD_{h_1}^{(1)}\vt_{b_r'}\big)\!+\!
v_{h_3}v_{h_4}\big(\cD_{h_4}^{(1)}\vt_{b_r'}\big)\!+
\!v_{h_3}v_{h_5}\big(\cD_{h_5}^{(1)}\vt_{b_r'}\big)$
\end{tabular}\end{small}}
\end{pspicture}
\caption{An example of the estimates \e_ref{derivexp_e1} and~\e_ref{derivexp_e2}}
\label{derivest_fig}
\end{figure}

\noindent
We now complete the proof of Proposition~\ref{cuspmap_prp}.
By the $m\!=\!1$ case of~\e_ref{derivexp_e1} and~\e_ref{derivexp_e2},
\begin{equation}\label{cuspmap_e1}
\Big| \cD_{\hat{0}}^{(1)}\vt_{b_r} -
\sum_{i\in\chi(b)}\!\!\rho_i(\ups_r)\big(\cD_i^{(1)}\vt_{b_r'}\big)  \Big|
\le C
\big(\|J_{t_r}\!-\!J_0\|_{C^1}\!+\!|\ups_r|^{1/p}\!+\!\|\xi_r\|_{\ups_r,p,1}\big)
\sum_{i\in\chi(b)}\!\! \big|\rho_i(\ups_r)\big|.
\end{equation}
On the other hand, since $\Phi_{b_r'}(i,\i)\!=\!\Id$ and $\vt_{b_r'}(i,\i)\!=0$
for all $i\!\in\!\chi(b)$,
\begin{equation}\label{cuspmap_e2}
\cD_ib_r'=\cD_i^{(1)}\vt_{b_r'} \qquad\forall~i\!\in\!\chi(b). 
\end{equation}
Furthermore, since 
$$\big|\Phi_{b_r}(\hat{0},\i)\!-\!\Id\big|
\le C\big\|\Phi_{b_r}\!-\!\Phi_{b_r'}\!\circ\!q_{\ups_r}\big\|_{\ups_r,p,1}
\le  C'\big(\|J_0\!-\!J_{t_r}\|_{C^1}\!+\!|\ups_r|^{1/p}\!+\!\|\xi_r\|_{\ups_r,p,1}\big)$$
and $\vt_{b_r}(\hat{0},\i)\!=0$,
\begin{equation}\label{cuspmap_e3}
\big|\cD_{\hat{0}}b_r-\cD_{\hat{0}}^{(1)}\vt_{b_r}\big|
\le C'\big(\|J_0\!-\!J_{t_r}\|_{C^1}\!+\!|\ups_r|^{1/p}\!+\!\|\xi_r\|_{\ups_r,p,1}\big)
\big|\cD_{\hat{0}}^{(1)}\vt_{b_r}\big|.
\end{equation}
By \e_ref{cuspmap_e1}-\e_ref{cuspmap_e3},
\begin{equation}\label{cuspmap_e4}
\Big|\cD_{\hat{0}}b_r-
\sum_{i\in\chi(b)}\!\!\rho_i(v_r)\big(\cD_ib_r'\big)\Big|
\le C\big(\|J_0\!-\!J_{t_r}\|_{C^1}\!+\!|\ups_r|^{1/p}\!+\!\|\xi_r\|_{\ups_r,p,1}\big)
\sum_{i\in\chi(b)}\!\!\big|\rho_i(v_r)\big|
\end{equation}
for all $r$ sufficiently large.
Since 
$$\lim_{r\lra\i}
\big(\|J_0\!-\!J_{t_r}\|_{C^1}\!+\!|\ups_r|^{1/p}\!+\!\|\xi_r\|_{\ups_r,p,1}\big) =0$$
and $\cD_ib_r'\!\lra\!\cD_ib$ for all $i\!\in\chi(b)$, 
\e_ref{cuspmap_e4} implies that $b$ must satisfy one of
the two conditions in the statement of Proposition~\ref{cuspmap_prp}.\\

\noindent
We next complete the proof of Proposition~\ref{node_prp}.
By the assumption on the bubble maps $b_k$ made in Proposition~\ref{node_prp}
and by Definition~\ref{pert0_dfn2},
$\ev_0(b_k)\!=\!\ev_1(b_k)$ for all~$k$.
Thus,
$$\ev_1(b_k)=\ev_0(b_k)=\ev_1(b_l) \qquad\forall~ k,l\in[n].$$
Let $q$ denote the point $\ev_0(b_1)$.
We identify a small neighborhood of $q$ in $X$ with a small neighborhood 
of $q$ in $T_qX$ via the exponential map $\exp$ and 
the tangent space to $X$ at a point close to $q$
with $T_qX$ via the parallel transport with respect to the $J_0$-linear
connection~$\na$.\\

\noindent
For each pair $(k,r)$, with $r$ sufficiently large, 
let $(b_{k,r}',\ups_{k,r},\xi_{k,r})$ be an analogue of $(b_r',\ups_r,\xi_r)$ 
for~$b_{k,r}$.
We~put
\begin{gather*}
\ze_{k,r}=\ev_0(b_{k,r}')\in T_qX  \qquad\hbox{and}\\
\ti\ze_{k,r}=\ev_1(b_{k,r})-\ev_0(b_{k,r})
=\ev_1(b_{k,r})-\ev_0(b_{k,r}') \in T_qX.
\end{gather*}
By the assumption on the maps $b_{k,r}$ made in the statement of Proposition~\ref{node_prp},
\begin{gather}\begin{split}
&\big|\ze_{k,r}+\ti\ze_{k,r}-\ze_{k+1,r}\big| 
\le C\big|\ze_{k,r}\big|\cdot\big|\ti\ze_{k,r}\big|
\qquad\forall~ k\!\in[n\!-\!1],\\
&\big|\ze_{n,r}+\ti\ze_{n,r}-\ze_{1,r}\big| 
\le C\big|\ze_{n,r}\big|\cdot\big|\ti\ze_{n,r}\big|;
\end{split}\notag\\
\label{node_e1}
\Lra\qquad 
\big|\ti\ze_{1,r}+\ldots+\ti\ze_{n,r}\big|
\le \ep_r\sum_{k=1}^{k=n}\big|\ti\ze_{k,r}\big|,
\end{gather}
for a sequence $\ep_r$ converging to $0$.\\

\noindent
On the other hand, the marked point $y_1(b_{k,r})\!=\!y_1(\ups_{k,r})$ of the bubble map
$b_{k,r}$ lies in $\Si_{\ups_{k,r}}^0(\de)$.
Thus,
\begin{equation}\label{node_e2}
\ti\ze_{k,r}=\vt_{b_{k,r}}\big(y_1(\ups_{k,r})\big)
=\sum_{m=1}^{\i} y_1(\ups_{k,r})^{-m}\big(\cD_{\hat{0}}^{(m)}\vt_{b_{k,r}}\big),
\end{equation}
where $y_1(\ups_{k,r})\!\in\!S^2\!-\!\{\i\}$ is viewed as a complex number.
Combining \e_ref{node_e2} with~\e_ref{derivexp_e1} and~\e_ref{derivexp_e2}
and then taking the lowest-order terms, we obtain an expression of the form
\begin{equation}\label{node_e3}\begin{split}
&\Big|\ti\ze_{k,r}-
\sum_{i\in\chi(b_k)}\!\!(y_{1;i}(b_{k,r}')\!-\!x_{i;1}(b_{k,r}'))^{-1}
\rho_{i;1}(\ups_{k,r}) \big(\cD_i^{(1)}\vt_{b_{k,r}'}\big) \Big| \\
&\qquad\qquad\qquad\qquad\qquad\qquad
\le C\big(\|J_0\!-\!J_{t_r}\|_{C^1}\!+\!|\ups_{k,r}|^{1/p}\!+\!\|\xi_{k,r}\|_{\ups_{k,r},p,1}\big)
\sum_{i\in\chi(b_k)}\!\!\big|\rho_{i;1}(\ups_{k,r})\big|;
\end{split}\end{equation}
see the proof of Corollary~3.7 in~\cite{g1comp}
for a derivation and the notation involved.
For the present purposes, the only fact we need to know about~\e_ref{node_e3} is that
\begin{equation}\label{node_e3b}
0<\big|\rho_{i;1}(\ups_{k,r})\big| \le  \big|\ups_{k,r}\big|
\qquad\forall~i\!\in\!\chi(b_k), ~ k\!\in\![n].
\end{equation}
In particular,  $\rho_{i;1}(\ups_{k,r})$ is a sequence of nonzero complex numbers
that approaches zero as $r$ tends to infinity.
By \e_ref{node_e1} and \e_ref{node_e3},
\begin{equation}\label{node_e4}
\Big|\sum_{k=1}^{k=n}\sum_{i\in\chi(b_k)}\!\!(y_{1;i}(b_{k,r}')\!-\!x_{i;1}(b_{k,r}'))^{-1}
\rho_{i;1}(\ups_{k,r}) \big(\cD_i^{(1)}\vt_{b_{k,r}'}\big) \Big|
\le\ti\ep_r \sum_{k=1}^{k=n}\sum_{i\in\chi(b_k)}\!\!\big|\rho_{i;1}(\ups_{k,r})\big|,
\end{equation}
for a sequence $\{\ti\ep_r\}$ converging to $0$.
Since 
$$\cD_i^{(1)}\vt_{b_{k,r}'}= \cD_i b_{k,r}' \lra \cD_ib_k
\qquad\hbox{as}\quad r\!\lra\!\i,$$
\e_ref{node_e3b} and \e_ref{node_e4} imply the conclusion of Proposition~\ref{node_prp}.

\subsection{Proof of Proposition~\ref{torus_prp}} 
\label{comp1prp_subs}

\noindent
We prove Proposition~\ref{torus_prp} by combining elements of the previous subsection
with a version of the two-stage gluing construction of \cite[Section~4]{g1comp}.
At the first stage, we smooth out all nodes of an element $[\Si';u']$ close to~$b$
that lie away from the principal component(s) $\Si_P'$ of~$\Si'$.
This stage will be unobstructed.
The objective of the second stage of the gluing construction is to smooth out
the remaining nodes of~$\Si'$.
We obtain Proposition~\ref{torus_prp} by estimating the obstruction to 
achieving this objective.\\

\noindent
Let $(X,\om)$, $\under{J}$, $A$, $M$, $\under\nu$, 
$$b=(M,I,\aleph;S,x,(j,y),u), \qquad\hbox{and}\qquad u_i\equiv u_b|_{\Si_{b,i}}$$ 
be as in the statement of Proposition~\ref{torus_prp}.
For each $i\!\in\!I$, we define
$$D_{J_0,\nu_0;b,i}\!: \Ga(b;i) \lra \Ga^{0,1}(b;i)$$
as at the beginning of the previous subsection.
With $I^+\!\subset\!I$ as in~\e_ref{iplusdfn_e}, choose
$$\ti\Ga^{0,1}_-(b;i)\subset \Ga\big(\Si_{b,i}\!\times\!X;
\La^{0,1}_{J_0,j}\pi_1^*T^*\Si_{b,i}\!\otimes\!\pi_2^*TX\big)$$
as in Subsection~\ref{g0prp_subs}.
Let $\T$ be the bubble type of the map~$b$.
We put
\begin{equation*}\begin{split}
\ti\U= \big\{ b'\!\equiv\!(M,I,\aleph;S',x',(j,y'),u')\!:~& [b']\!\in\!\X_{\T}(X),\\
& \pi_i\{\bar\partial_{J_0,j}\!+\!\nu_0\}u_{b'}\in
\{\id\!\times\!u_i'\}^*\ti\Ga^{0,1}_-(b;i)~\forall i\!\in\!I\big\},
\end{split}\end{equation*}
where
$$\pi_i\!: \Ga^{0,1}(b';J_0) \lra \Ga^{0,1}(b';i)$$
is the natural projection map.
By the Implicit Function Theorem, $\ti\U$ is a smooth manifold near~$b$.
If $b'\!\in\!\ti\U$, $u_{b'}|_{\Si_{b';P}}$ is a degree-zero $J_0$-holomorphic
map and thus is constant. Let
$$\ev_P\!: \ti\U\lra X, \qquad b'\lra u_{b'}({\Si_{b';P}}),$$
be the map sending each element $b'$ of $\ti\U$ to the image of the main component
of its domain.\\

\noindent
For each $b'\!\in\!\ti\U$, let 
$$\Ga_-(b') = \big\{\xi\!\in\!\Ga(b')\!: 
\pi_iD_{J_0,\nu_0;b'}\xi\in\{\id\!\times\!u_{b',i}\}^*\ti\Ga^{0,1}_-(b;i)
~\forall i\!\in\!I\big\},$$
where $u_{b',i}\!=\!u_{b'}|_{\Si_{b',i}}$.
We denote by
$$\ti\Ga^{0,1}_-(b';X) \subset \Ga\big(\Si_{b'}\!\times\!X;
\La^{0,1}_{J_0,j}\pi_1^*T^*\Si_{b'}\!\otimes\!\pi_2^*TX\big)$$
the subspace obtained by extending all elements of each of the spaces 
$$\ti\Ga_-^{0,1}(b';i) \equiv\ti\Ga_-^{0,1}(b;i)$$
with $i\!\in\!\hat{I}$ by zero outside of the component $\Si_{b',i}$ of~$\Si_{b'}$.\\

\noindent
We put 
$$I_1=\big\{h\!\in\!\hat{I}\!:\io_h\!\in\!I_0\big\},$$
where $I_0$ is the subset of minimal elements of $I$.
Let
$$\wt{\cal F}\lra \ti\U$$
be the bundle of gluing parameters.
In this case, $\wt{\cal F}$ has three distinguished components:
\begin{gather*}
\wt{\cal F}=\wt{\cal F}_{\aleph} \oplus \wt{\cal F}_0 \oplus \wt{\cal F}_1, 
\qquad\hbox{where}\\
\wt{\cal F}_{\aleph}=\ti\U\times\C^{\aleph}, \qquad
\wt{\cal F}_1=\ti\U\times\C^{\hat{I}-I_1}, \quad\hbox{and}\quad
 \wt{\cal F}_0\big|_{b'}=\bigoplus_{h\in I_1}T_{x_h(b')}\Si_{b';P}~~\forall~b'\!\in\!\ti\U.
\end{gather*}
The total space of $\wt{\cal F}_0$ has a natural topology;
see \cite[Subsection~2.2]{g1comp}.
We denote by $\wt{\cal F}^{\eset}$ the subset of  $\wt{\cal F}$ consisting
of the elements with all components nonzero.
If $i\!\in\!\hat{I}$, let $h(i)\!\in\!I_1$ be the unique element such that
$h(i)\!\le\!i$.
For each $\ups\!=\!(b',v)$, where $b'\!\in\!\ti\U$
and $v\!=\!(v_i)_{i\in\aleph\sqcup\hat{I}}$, and $i\!\in\!\chi(b)$, we~put
\begin{gather*}
\ups_0=\big(b',(v_i)_{i\in\aleph\sqcup I_1}\big), \qquad
\ups_1=\big(b',(v_i)_{i\in\hat{I}- I_1}\big),\\
\ti\rho_i(\ups)=\!\prod_{h(i)< h\le i}\!\!\!\!\!\! v_h \in\Bbb{C},
\qquad\hbox{and}\qquad  \rho_i(\ups)=\ti\rho_i(\ups)\cdot v_{h(i)}
\in T_{x_{h(i)}(b')}\Si_{b';P}.
\end{gather*}
The component $\ups_1$ of $\ups$ consists of the smoothings of the nodes 
of $\Si_b$ that lie away from the principal component.\\

\noindent
For each sufficiently small element $\ups\!=\!(b',v)$ of $\wt{\cal F}^{\eset}$, let 
$$q_{\ups_1}\!:\Si_{\ups_1}\lra\Si_{b'}$$
be the basic gluing map constructed in \cite[Subsection~2.2]{gluing}.
In this case, the  principal component $\Si_{\ups_1;P}$ of $\Si_{\ups_1}$
is the same as principal component $\Si_{b';P}$ of~$\Si_{b'}$,
and $\Si_{\ups_1}$ has $|I_1|$ bubble components $\Si_{\ups_1,h}$,
with $h\!\in\!I_1$, attached directly to~$\Si_{\ups_1;P}$.
The map~$q_{\ups_1}$ collapses $|\hat{I}\!-\!I_1|$ circles on the bubble components
of~$\Si_{\ups_1}$.
It induces a metric $g_{\ups_1}$ on~$\Si_{\ups_1}$ such that
$(\Si_{\ups_1},g_{\ups_1})$ is obtained from $\Si_{b'}$ by replacing 
$|\hat{I}\!-\!I_1|$ nodes by thin necks.
Let
$$u_{\ups_1}=u_{b'}\circ q_{\ups_1}.$$
The map~$q_{\ups_1}$ induces norms $\|\cdot\|_{\ups_1,p,1}$ and~$\|\cdot\|_{\ups_1,p}$
on the spaces 
$$\big\{\xi\!\in\!\Ga(\Si_{\ups_1};u_{\ups_1}^*TX)\!:\xi|_{\Si_{\ups_1;P}}\!=\!0\big\}
\qquad\hbox{and}\qquad
\big\{\eta\!\in\!\Ga(\Si_{\ups_1};\La^{0,1}_{J_0,j}T^*\Si_{\ups_1}\!\otimes\!u_{\ups_1}^*TX)\!:
\eta|_{\Si_{\ups_1;P}}\!=\!0\big\}$$
respectively; see \cite[Subsection~3.3]{gluing}.
We denote the corresponding completions by $\Ga_B(\ups_1)$ and $\Ga_B^{0,1}(\ups_1)$.\\

\noindent
{\it Remark:} The weights for the norms $\|\cdot\|_{\ups_1,p,1}$ and $\|\cdot\|_{\ups_1,p}$
are constructed as \cite[Subsection~3.3]{gluing}, but on each of the $|I_1|$
 bubbles separately.
The restrictions of these norms to each of the bubbles 
are equivalent to the norms used in \cite[Section~3]{LT}.\\

\noindent
Fix $\ep_b\!\in\!\R^+$ such that for every $h\!\in\!I_1$ the disk of radius of $8\ep_b$
in $\Si_{b;P}$ around the node $x_h(b)$ contains no other special, i.e.~singular or
marked, point of~$\Si_b$.
For each $(b',v)\!\in\!\wt{\cal F}^{\eset}$ with $b'\!\in\!\ti\U$ sufficiently close to~$b$ 
and  $v$ sufficiently small, let
$$q_{\ups_0;2}\!:\Si_{\ups}\lra\Si_{\ups_1}
\qquad\hbox{and}\qquad
\ti{q}_{\ups_0;2}\!:\Si_{\ups}\lra\Si_{\ups_1}$$
be the basic gluing map of \cite[Subsection~2.2]{gluing} corresponding 
to the gluing parameter~$v_0$ and the modified basic gluing map defined
in the middle of Subsection~4.2 in~\cite{g1comp}
with the collapsing radius~$\ep_b$.
In this case, $\Si_{\ups}$ is a smooth genus-one curve.
For each $h\!\in\!I_1$, the maps $q_{\ups_0;2}$ and $\ti{q}_{\ups_0;2}$
collapse the circles of radii $|v_h|^{1/2}$ and $\ep_b$, respectively, 
around the point $x_h(b')\!\in\!\Si_{\ups_1;P}$. 
Once again, the~map 
$$q_{\ups} \equiv q_{\ups_0;2}\!\circ q_{\ups_1}\!: \Si_{\ups}\lra\Si_{b'}$$ 
induces a metric $g_{\ups}$ on~$\Si_{\ups}$ such that $(\Si_{\ups},g_{\ups})$ is 
obtained from $\Si_{b'}$ by replacing all nodes by thin necks.\\

\noindent
For each $t\!\in\![0,1]$, let $g_t$ be the $J_t$-compatible symmetric 
two-tensor on $X$ given by
$$g_t(\cdot,\cdot)=\frac{1}{2}
\big( g(J_0\cdot,J_t\cdot)+g(J_t\cdot,J_0\cdot)\big).$$
If $t$ is sufficiently close to $0$, $g_t$ is positive-definite, i.e.~it is 
a $J_t$-compatible metric on~$X$.
We denote the $J_t$-compatible connection induced by the Levi-Civita connection 
of the metric~$g_t$ by~$\na^t$
and the corresponding exponential map by~$\exp^t$.\\

\noindent
If $W$ is a small neighborhood of $b$ in $\ti\U$ and $\de\!\in\!\R^+$ is
sufficiently small, let
$$\X_1(W,\de)=\big\{ (\Si_{\ups_1};u_{\ups_1,\xi})\!\equiv\!
\big(\Si_{\ups_1};\exp_{u_{\ups_1}}\!\xi\big)\!:
\ups\!=\!(b',v)\!\in\!\wt{\cal F}_{\de}^{\eset}|_W;~
 \xi\!\in\!\Ga_B(\ups_1),~\|\xi\|_{\ups_1,p,1}\!\le\!\de\big\}.$$
For each element $(\Si_{\ups_1};u_{\ups_1,\xi})$ of $\X_1(W,\de)$, we~put
\begin{equation}\label{approxmap2_e}
u_{\ups,\xi}= u_{\ups_1,\xi}\circ \ti{q}_{\ups_0;2}.
\end{equation}
The map~$q_{\ups_0;2}$ induces norms $\|\cdot\|_{\ups,p,1}$ and~$\|\cdot\|_{\ups,p}$
on the spaces 
$$\Ga(\Si_{\ups};u_{\ups,\xi}^*TX)  \qquad\hbox{and}\qquad
\Ga(\Si_{\ups};\La^{0,1}_{J_t,j}T^*\Si_{\ups}\!\otimes\!u_{\ups,\xi}^*TX),$$
respectively.
For $t$ sufficiently small, we use the metric $g_t$ on $X$ to define a norm
on the latter space.
Let $\Ga(\ups;\xi)$ and~$\Ga_t^{0,1}(\ups;\xi)$ be the corresponding completions.
If $\ups\!=\!(b',v)$, we~put
$$\Ga_-(\ups;\xi)=\big\{ \big(\Pi_{\xi}(\ze\!\circ\!q_{\ups_1})\big)
\circ\ti{q}_{\ups_0;2}\!: \ze\!\in\!\Ga_-(b')\big\}
\subset \Ga(\ups;\xi),$$
where $\Pi_{\xi}$ is the $\na$-parallel transport along the $\na$-geodesics
$\tau\!\lra\!\exp_{u_{\ups_1}}\!\tau\xi$.
Let 
$$\Ga_+(\ups;\xi) \subset \Ga(\ups;\xi)$$
be the $L^2$-orthogonal complement of $\Ga_-(\ups;\xi)$ defined with respect to
the metrics $g_{\ups}$ on $\Si_{\ups}$ and $g$ on~$X$.
For every $t\!\in\![0,1]$, we denote~by
$$\pi_t^{0,1}\!:\Ga(\Si_{\ups};T^*\Si_{\ups}\!\otimes_{\R}\!u_{\ups,\xi}^*TX)
\lra \Ga(\Si_{\ups};\La^{0,1}_{J_t,j}T^*\Si_{\ups}\!\otimes\!u_{\ups,\xi}^*TX)$$
the natural projection map.\\

\noindent
For each $b'\!\in\!\ti\U$, each sufficiently small element $\ups\!=\!(b',v)$
of $\wt{\cal F}^{\eset}$, and $\de\!\in\!\R^+$, we define
$$A_{b',i}(\de),\Si_{b'}^0(\de)\subset \Si_{b'},        \qquad
\Si_{\ups_1}^0(\de)\subset \Si_{\ups_1},  \quad\hbox{and}\quad
\Si_{\ups}^0(\de)\subset \Si_{\ups}$$
by \e_ref{mainpart_e1} and~\e_ref{mainpart_e2}.
Choose a neighborhood $W$ of $b$ in $\ti\U$, $\de_b\!\in\!(0,\ep_b)$, and
$\de\!\in\!(0,\de_b^2/2)$ such~that\\
${}~~$ (i) the maps $q_{\ups}$ and $\ti{q}_{\ups_0;2}$ are defined for all 
$\ups\!\in\!\wt{\cal F}_{\de}^{\eset}|_W$;\\
${}~\,$ (ii) for all $i\!\in\!\chi(b)$ all elements of $\ti\Ga_-^{0,1}(b;i)$ 
vanish on $A_{b,i}(\de_b)$;\\
${}~$ (iii) $\nu_t(b')\big|_{\Si_{b'}^0(2\de_b)}\!=\!0$  for every $t\!\in\![0,1]$
and $b'\!\in\!W$;\\
${}~$ (iv) $\nu_t(\ti{b})\big|_{\Si_{\ups}^0(2\de_b)}\!=\!0$ 
for every $t\!\in\![0,1]$ and every $\ti{b}\!=\!(\Si_{\ups};\exp_{u_{\ups_1,\xi}}\!\ze)$
such that\\ 
${}\qquad$ $\ups\!\in\!\wt{\cal F}_{\de}^{\eset}|_W$, 
$(\Si_{\ups_1};u_{\ups_1,\xi})\!\in\!\X_1(W,\de)$, $\ze\!\in\!\Ga(\ups;\xi)$,
and $\|\ze\|_{\ups,p,1}\!\le\!\de$.\\
Such a positive number $\de$ exists by our assumptions on the spaces
$\ti\Ga_-^{0,1}(b;i)$ and the family of perturbations~$\under\nu$;
see Definition~\ref{pert_dfn}.\\

\noindent
Suppose $\ups\!\in\!\wt{\cal F}_{\de}^{\eset}|_W$ and
$(\Si_{\ups_1};u_{\ups_1,\xi})\!\in\!\X_1(W,\de)$.
By the construction of the map~$\ti{q}_{\ups_0;2}$ in Subsection~4.2 
of~\cite{g1comp} and the assumptions that $\de\!\le\!\de_b^2/2$ and $\de_b\!\le\!\ep_b$,
$$\ti{q}_{\ups_0;2}\!: \Si_{\ups}-\Si_{\ups}^0(\de_b/2)
\lra \Si_{\ups_1}-\Si_{\ups_1}^0(\de_b/2)$$
is a biholomorphism.
Thus, by the assumption~(iv) above with $\ze\!=\!0$, for every $t\!\in\![0,1]$, 
we can define 
$$\nu_{t;\ups_1,\xi}\in\Ga^{0,1}\big(\Si_{\ups_1};
 \La_{J,j}^{0,1}T^*\Si_{\ups_1}\!\otimes\!u_{\ups_1,\xi}^*TX\big)
\quad\hbox{by}\quad
\ti{q}_{\ups_0;2}^*\nu_{t;\ups_1,\xi}=\nu_t(\Si_{\ups};u_{\ups,\xi}), ~~~
\nu_{t;\ups_1,\xi}\big|_{\Si_{\ups_1}^0(\de_b)}=0.$$
If $\ups\!=\!(b',v)$, we put
$$\Ga_{t,\de}(\ups)=\big\{\xi\!\in\!\Ga_B(\ups_1)\!: \|\xi\|_{\ups_1,p,1}\!\le\!\de;~
\big\{\bar\partial_J\!+\!\nu_{t;\ups_1,\xi}\big\}u_{\ups_1,\xi} \in
\pi_t\big\{q_{\ups_1}\!\times\!u_{\ups_1,\xi}\big\}^*\ti\Ga_-^{0,1}(b';X)\big\}.$$\\

\noindent
Let $t_r$ and $b_r$ be as in Proposition~\ref{torus_prp}.
Since the sequence $[b_r]$ converges to $[b]$, for all $r$ sufficiently large
there~exist 
$$b_r'\in W, \quad \ups_r=(b_r',v_r)\in\wt{\cal F}_{\de}^{\eset},  \qquad
\xi_r\in\Ga_{t,\de}(\ups_r),  \quad\hbox{and}\quad
\ze_r\in\Ga_+(\ups_r;\xi_r)$$
such that
\begin{gather}\label{torus_e1}
\lim_{r\lra\i}b_r'=b, \quad \lim_{r\lra\i}|v_r|=0, \quad
\lim_{r\lra\i}\|\xi_r\|_{\ups_{r;1},p,1}=0,  \quad
\lim_{r\lra\i}\|\ze_r\|_{\ups_r,p,1}=0,\\
\label{torus_e2}
\hbox{and}\qquad b_r\!\equiv\!\big(\Si_{b_r};u_{b_r}\big)=
\big(\Si_{\ups_r};\exp_{u_{\ups_r,\xi_r}}^t\!\ze_r\big).
\end{gather}
The last equality holds for a representative $b_r$ for $[b_r]$.
The existence of $(\ups_r,\xi_r,\ze_r)$ for 
\begin{equation}\label{torus_e3}
b_r\in\M_{1,M}^0(X,A;J_{t_r},\nu_{t_r})
\end{equation}
as above will imply that $b$ satisfies the second property in Definition~\ref{degen_dfn}.\\

\noindent
{\it Remark:} Similarly to the genus-zero case, the existence of 
$(\ups_r,\xi_r,\ze_r)$ can be shown by a variation on the surjectivity argument of 
\cite[Section~4]{gluing}; see also the paragraph following Lemma~4.4 in~\cite{g1comp}.
This is also the case if the map $q_{\ups_0;2}$ is used instead of $\ti{q}_{\ups_0;2}$
in~\e_ref{approxmap2_e} .
However, using the map $\ti{q}_{\ups_0;2}$ in~\e_ref{approxmap2_e} makes 
the maps $u_{\ups,\xi}$, with $\xi\!\in\!\Ga_{t,\de}(\ups)$,
closer to being $(J_t,\nu_t)$-holomorphic.
Since $u_{b'}|_{\Si_{b';P}}$ is constant, the choice of $\ti{q}_{\ups_0;2}$
for constructing approximately $(J_t,\nu_t)$-holomorphic maps 
is analogous to that of Section~3 in~\cite{LT}. \\

\noindent
For each $\ups\!\in\!\wt{\cal F}_{\de}^{\eset}|_W$ and $\xi\!\in\!\Ga_{t,\de}(\ups)$, 
let
$$D_{\ups;\xi}^t\!: \Ga(\ups;\xi) \lra \Ga_t^{0,1}(\ups;\xi)$$
be the linearization of section $\bar{\partial}_{J_t}\!+\!\nu_t$ at 
$(\Si_{\ups};u_{\ups,\xi})$ defined via the connection~$\na^t$.
We denote~by
$$\Ga_{t;+}^{0,1}(\ups;\xi) \subset \Ga_t^{0,1}(\ups;\xi)$$
the image $\Ga_+(\ups;\xi)$ under $D_{\ups;\xi}^t$.
By the same argument as in  \cite[Subsection~5.4]{gluing}, 
there exists $C\!\in\!\R^+$ such~that
\begin{equation}\label{torus_e4}
C^{-1}\|\ze\|_{\ups,p,1}\le \big\|D_{\ups;\xi}^t\ze\big\|_{\ups,p}
\le C\|\ze\|_{\ups,p,1} 
\quad \forall~ t\!\in\![0,\de],~
\ups\!\in\!\wt{\cal F}_{\de}^{\eset}|_W, ~\xi\!\in\!\Ga_{t,\de}(\ups),
~\ze\!\in\!\Ga_+(\ups;\xi),
\end{equation}
provided $W$ and $\de$ are sufficiently small.\\

\noindent
Put
$$\ti{\Ga}_{t;-}^{0,1}(\ups;\xi)=\big\{\pi_t\{q_{\ups}\!\times\!u_{\ups,\xi}\}^*\eta\!:
\eta\!\in\!\ti\Ga_-^{0,1}(b';X)\big\} \subset \Ga_t^{0,1}(\ups;\xi).$$
If $\de$ is sufficiently small, by the same argument as in \cite[Subsection~3.5]{gluing}
and our assumptions on the spaces $\ti\Ga_-(b;i)$,
\begin{equation}\label{ga01decomp_e}
\Ga_t^{0,1}(\ups;\xi)=\Ga_{t;+}^{0,1}(\ups;\xi)\oplus
\ti{\Ga}_{t;-}^{0,1}(\ups;\xi) \oplus \Ga_{t;-}^{0,1}(\ups;\xi),
\end{equation}
for some subspace $\Ga_{t;-}^{0,1}(\ups;\xi)$ of $\Ga_t^{0,1}(\ups;\xi)$ 
isomorphic to the cokernel of the composition:
$$\pi\circ D_{J_0,\nu_0;b'}\!:
\Ga(b') \lra \Ga^{0,1}(b';J_0)\lra 
\Ga^{0,1}(b';J_0)\big/  \big\{\id\!\times\!u_{b'}\big\}^*\ti{\Ga}^{0,1}_-(b';X).$$
This cokernel is naturally isomorphic to
$$\Ga_-^{0,1}(b';J_0) \equiv \H_{b'}^{0,1}\otimes_{J_0} T_{\ev_P(b')}X \subset 
\Ga^{0,1}(b';J_0),$$
where $\H_{b'}^{0,1}$ is the one-dimensional complex vector space of 
$(0,1)$-harmonic forms on the principal component~$\Si_{b';P}$ of~$\Si_{b'}$.
If $\Si_{b';P}$ is a circle of spheres,
the elements of $\H_{b'}^{0,1}$ have simple poles at the nodes of $\Si_{b';P}$
with the residues adding up to zero at each node.
Recall that $(\Si_{\ups},g_{\ups})$ is obtained from $\Si_{b'}$ by replacing 
the nodes of $\Si_{b'}$ with thin necks.
The~map
$$q_{\ups}\!: \Si_{\ups}\lra\Si_{b'}$$
collapses each neck at its thinnest position to the corresponding node.
For each element $\eta$ of $\Ga_-^{0,1}(b';J_t)$, we can construct an element
$R_{\ups,\xi}\eta$ of $\Ga_t^{0,1}(\ups;\xi)$ by parallel transporting $\eta$
along the restriction of $u_{\ups,\xi}$ to $q_{\ups}^{-1}(A_{b,h}(\de_b^2))$
for each $h\!\in\!I_1$ and cutting it off with a smooth function that drops
from $1$ to~$0$ over the annulus $q_{\ups}^{-1}(A_{b,h}(\de_b^2)\!-\!A_{b,h}(\de_b^2/4))$;
see the middle of Subsection~4.2 in~\cite{g1comp} for details.
We~denote by $\Ga_{t;-}^{0,1}(\ups;\xi)$ the image of $\Ga_-^{0,1}(b';J_t)$ 
under~$R_{\ups}$.\\

\noindent
{\it Remark:} In the construction of the map $\ti{q}_{\ups_0;2}$ in 
\cite[Subsection~4.2]{g1comp}, 
$\de_K$ corresponds to $\ep_b^2$ above.
In the construction of $R_{\ups,\xi}\eta$ in \cite[Subsection~4.2]{g1comp},
$\de_K$ corresponds to~$\de_b^2/4$ above.\\

\noindent
Let $\llrr{\cdot,\cdot}_t$ denote the Hermitian inner-product on $\Ga^{0,1}_t(\ups;\xi)$
induced by the $J_t$-compatible metric~$g_t$ on~$X$.
For each $\eta\!\in\!\Ga^{0,1}_-(b';J_t)$, let
$$\|\eta\|=\sum_{h\in I_1}\blr{\eta_{x_h(b')},\eta_{x_h(b')}}_t,$$
where $\lr{\cdot,\cdot}_t$ is the hermitian inner-product
on ${\cal H}_{b'}^{0,1}\!\otimes_{J_t}\!T_{\ev_P(b')}X$ 
defined via the metric $g$ on $X$ and the original metric $g_{b'}$ on~$\Si_{b'}$.
From the construction of $R_{\ups}$ in  \cite[Subsection~4.2]{g1comp}
and Holder's inequality, it is immediate that
\begin{equation}\label{g1reg_lmm1e6}
\big|\bllrr{\eta',R_{\ups,\xi}\eta}_t\big|
\le C\|\eta\|\|\eta'\|_{\ups,p}
\qquad\forall~\eta\!\in\!\Ga_-^{0,1}(b';J_t),~\eta'\!\in\!\Ga^{0,1}_t(\ups;\xi);
\end{equation}
see \cite[(4.11)]{g1comp}.
Another essential property of the above construction is that
\begin{equation}\label{torus_e5}
\big|\bllrr{D_{\ups;\xi}^t\ze,R_{\ups,\xi}\eta}_t\big|
\le C|\ups|^{1/2}\|\eta\|\|\ze\|_{\ups,p,1}
\quad \forall~ t\!\in\![0,\de],~
\ups\!\in\!\wt{\cal F}_{\de}^{\eset}|_W, ~\xi\!\in\!\Ga_{t,\de}(\ups),
~\ze\!\in\!\Ga(\ups;\xi);
\end{equation}
see part (7) of Lemma~4.4 in~\cite{g1comp}.\\

\noindent
Due to the assumption~\e_ref{torus_e2},
 the condition~\e_ref{torus_e3} is equivalent~to
\begin{equation}\label{torus_e6a}
\big\{\bar\partial_{J_{t_r}}\!+\!\nu_{t_r}\big\}u_{\ups_r,\xi_r} +
D_{\ups_r;\xi_r}^{t_r}\ze_r + N_{\ups_r;\xi_r}^{t_r}\ze_r=0.
\end{equation}
The quadratic term $N_{\ups;\xi}^t$ satisfies
\begin{equation}\label{torus_e6b}
\big\|N_{\ups;\xi}^t\ze\big\|_{\ups,p} \le C  \|\ze\|_{\ups,p,1}^2
\quad \forall~ t\!\in\![0,\de],~
\ups\!\in\!\wt{\cal F}_{\de}^{\eset}|_W, ~\xi\!\in\!\Ga_{t,\de}(\ups), ~
\ze\!\in\!\Ga(\ups;\xi) ~\st~ \|\ze\|_{\ups,p,1}\le\de.
\end{equation}
We will obtain Proposition~\ref{torus_prp} by estimating the  inner-product 
$\llrr{\cdot,\cdot}_{t_r}$
of each term in~\e_ref{torus_e6a} with each element of 
$\Ga_{t_r;-}^{0,1}(\ups_r;\xi_r)$.\\

\noindent
First, for every $\ups\!\in\!\wt{\cal F}_{\de}^{\eset}|_W$ and 
$\xi\!\in\!\Ga_{t,\de}(\ups)$, let
$$\ti\pi_{\ups;\xi}^{t;+}\!: 
\Ga_t^{0,1}(\ups;\xi) \lra \Ga_{t;+}^{0,1}(\ups;\xi)
\oplus \Ga_{t;-}^{0,1}(\ups;\xi)$$
be the projection map corresponding to the decomposition~\e_ref{ga01decomp_e}.
Then,
\begin{equation}\label{torus_e7a}\begin{split}
\ti\pi_{\ups;\xi}^{t;+}\big\{\bar\partial_{J_t}\!+\!\nu_t\big\}u_{\ups,\xi}
&= \ti\pi_{\ups;\xi}^{t;+}\big\{\bar\partial_{J_t}\!+\!
                              \ti{q}_{\ups_0;2}^*\nu_{t;\ups_1,\xi}\big\}
\big(u_{\ups_1,\xi}\!\circ\!\ti{q}_{\ups_0;2}\big) \\
& = \ti\pi_{\ups;\xi}^{t;+} \big(
\big(\!\big\{\bar\partial_{J_t}\!+\!\nu_{t;\ups_1,\xi}\big\}u_{\ups_1,\xi}\big)
\circ \partial\ti{q}_{\ups_0;2} \big)
+\ti\pi_{\ups;\xi}^{t;+}\bar\partial_{J_t} (u_{\ups,\xi}\!\circ\!\ti{q}_{\ups_0;2})\\
& = \ti\pi_{\ups;\xi}^{t;+}\bar\partial_{J_t} (u_{\ups,\xi}\!\circ\!\ti{q}_{\ups_0;2}).
\end{split}\end{equation}
The reason for the second equality is that the map $\ti{q}_{\ups_0;2}$
is holomorphic over the support of~$\nu_{t;\ups_1,\xi}$.
This last equality follows from the definition of~$\Ga_{t,\de}(\ups)$ above.
By our assumptions on $b'$ and~$\xi$, 
$$u_{\ups_1,\xi}\big|_{q_{\ups_1}^{-1}(\Si_{b';P})}=\const
\qquad\Lra\qquad 
du_{\ups_1,\xi}\!\circ\!\bar\partial\ti{q}_{\ups_0;2} \,
\big|_{\ti{q}_{\ups_0;2}^{-1}q_{\ups_1}^{-1}(\Si_{b';P})}=0.$$
By the construction of the map $\ti{q}_{\ups_0;2}$, 
the restriction of $\ti{q}_{\ups_0;2}$ to the complement of
$\ti{q}_{\ups_0;2}^{-1}q_{\ups_1}^{-1}(\Si_{b';P})$ in $\Si_{\ups}$
is holomorphic outside of the annuli 
$${\cal A}_{\ups,h}\equiv 
\ti{q}_{\ups_0;2}^{-1}q_{\ups_1}^{-1}\big( A_{b',h}(2|v_h|^2/\ep_b^2) \big),$$
with $h\!\in\!I_1$.
The~map $q_{\ups}$ maps such an annulus isomorphically onto the annulus of
radii $\ep_b/2$ and~$\ep_b$ around the point $x_{h(b')}\!\in\!\Si_{b';P}$.
The key advantage of using the map $\ti{q}_{\ups_0;2}$ instead of $q_{\ups_0;2}$
in~\e_ref{approxmap2_e} is that
\begin{equation}\label{torus_e7b}
\big\|d\ti{q}_{\ups_0;2}\big\|_{C^0({\cal A}_{\ups,h})}
\le C|v_h|\qquad\forall~h\!\in\!I_1,
\end{equation}
where $C^0$-norm of $d\ti{q}_{\ups_0;2}$ is computed with respect to the metrics
$g_{\ups_1}$  and $g_{\ups}$ on $\Si_{\ups_1}$ and $\Si_{\ups}$, respectively.\\

\noindent
Since $u_{\ups_1,\xi}|_{\ti{q}_{\ups_0;2}({\cal A}_{\ups,h})}$ is $J_t$-holomorphic,
\e_ref{derivexp_e1} and~\e_ref{derivexp_e2} give
\begin{equation}\label{torus_e7c}
\big\|du_{\ups_1,\xi}|_{\ti{q}_{\ups_0;2}({\cal A}_{\ups,h})}
\big\|_{\ups_1,p}
\le C|v_h|^{2/p} \!\!\!\!\! \sum_{i\in\chi(\T),h(i)=h} \!\!\!\!\!\!\!\!\!
\big|\ti{\rho}_i(\ups)\big| \qquad\forall~h\!\in\!I_1;
\end{equation}
see part (2b) of Corollary~3.8 in~\cite{g1comp}.
From \e_ref{torus_e7a}-\e_ref{torus_e7c}, 
we conclude that
\begin{equation}\label{torus_e8a}
\big\|\ti\pi_{\ups;\xi}^{t;+}\big\{\bar\partial_{J_t}\!+\!\nu_t\big\}u_{\ups,\xi}\big\|_{\ups,p}
\le C \! \sum_{i\in\chi(\T)} \!\!\! \big|\rho_i(\ups)\big| 
\qquad\forall~ t\!\in\![0,\de],~
\ups\!\in\!\wt{\cal F}_{\de}^{\eset}|_W, ~\xi\!\in\!\Ga_{t,\de}(\ups),
\end{equation}
if $W$ and $\de$ are sufficiently small.
This estimate is the analogue of the first estimate in part (3) of
Lemma~4.4 in~\cite{g1comp}.
Separately, by the construction of the homomorphism~$R_{\ups,\xi}$,
$$\hbox{supp}\,R_{\ups,\xi}\eta \cap \hbox{supp}\,\ti{\eta} = \eset
\qquad \forall~ R_{\ups,\xi}\eta\!\in\!\Ga_{t;-}^{0,1}(\ups;\xi),
~\ti\eta\!\in\!\ti\Ga_{t;-}^{0,1}(\ups;\xi).$$
Thus, by the $J_t$-holomorphicity of
$u_{\ups_1,\xi}|_{\ti{q}_{\ups_0;2}({\cal A}_{\ups,h})}$,
\e_ref{derivexp_e1}, \e_ref{derivexp_e2}, \e_ref{torus_e7a},
and integration by~parts, 
\begin{equation}\label{torus_e8b}\begin{split}
\Big|\bllrr{ \big\{\bar{\partial}_{J_t}\!+\!\nu_t\big\}u_{\ups,\xi},R_{\ups,\xi}\eta}_t
+2\pi\I\!\!\sum_{i\in\chi(\T)}\!\!\!\!
\blr{\cD_ib',\eta_{x_{h(i)}(b')}(\rho_i(\ups))}\Big| \qquad\qquad\qquad\qquad&\\
\le C \big(\|J_t\!-\!J_0\|_{C^1}\!+\!|\ups|^{1/p}+
\!|\ups|^{(p-2)/p}\!+\!\|\xi\|_{\ups,p,1}\big)  \|\eta\|
\!\!\!\sum_{i\in\chi({\T})}\!\!\!\!|\rho_i(\ups)|;&
\end{split}\end{equation}
see the proof of part (6) of Lemma~4.4 at the end of Subsection~4.2 of~\cite{g1comp}.
Here $\lr{\cdot,\cdot}$ is the Hermitian inner-product on $(T_{\ev_P(b')}X,J_0)$ 
defined via the metric $g$ on~$X$.\\

\noindent
We now finish the proof of Proposition~\ref{torus_prp}.
By \e_ref{torus_e4}, \e_ref{torus_e6a}, \e_ref{torus_e6b}, \e_ref{torus_e8a},
and the definition of~$\Ga_+^{0,1}(\ups;\xi)$,
\begin{equation}\label{torus_e9c}
\|\ze_r\|_{\ups_r,p,1}\le C \!\sum_{i\in\chi(\T)}\!\!\! \big|\rho_i(\ups_r)\big|,
\end{equation}
for all $r$ sufficiently large.
Combining this estimate with \e_ref{g1reg_lmm1e6},
\e_ref{torus_e5}, and \e_ref{torus_e6b}, we obtain
that for all $\eta\!\in\!\Ga_-^{0,1}(b_r';J_0)$,
\begin{equation}\label{torus_e9a}\begin{split}
\big| \bllrr{D_{\ups_r;\xi_r}^{t_r}\ze_r,R_{\ups_r,\xi_r}\eta}_{t_r}\big|
&\le C|\ups_r|^{1/2}\|\eta\|
\!\sum_{i\in\chi(\T)}\!\!\! \big|\rho_i(\ups_r)\big|;\\
\big| \bllrr{N_{\ups_r;\xi_r}^{t_r}\ze_r,R_{\ups_r,\xi_r}\eta}_{t_r}\big|
&\le C|\ups_r| \|\eta\|
\!\sum_{i\in\chi(\T)}\!\!\! \big|\rho_i(\ups_r)\big|.
\end{split}\end{equation}
Finally, by~\e_ref{torus_e1}, \e_ref{torus_e6a}, \e_ref{torus_e8b}, and~\e_ref{torus_e9a}, 
for a sequence $\ep_r$ converging to~$0$
\begin{equation}\label{torus_e10}
\Big|\!\sum_{i\in\chi(\T)}\!\!\!\!
\blr{\cD_ib_r',\eta_{x_{h(i)}(b_r')}(\rho_i(\ups_r))}\Big| 
\le \ep_r \|\eta\|\!\!\!\sum_{i\in\chi({\T})}\!\!\!\!|\rho_i(\ups_r)|
\qquad\forall~ \eta\!\in\!{\cal H}_{b_r'}^{0,1}\!\otimes_{J_0}\!T_{\ev_P(b')}X.
\end{equation}
Since $\ups_r\!\in\!\wt{\cal F}^{\eset}$ for all $r$,
$\rho_i(\ups_r)\!\neq\!0$ for all $i\!\in\!\chi(b)$.
Thus, \e_ref{torus_e10} implies the conclusion of Proposition~\ref{torus_prp},
since $\cD_ib_r'\!\lra\!\cD_ib$ as $r\!\lra\!\i$.

\section{Proof of Theorem~\ref{gwdiff_thm}}
\label{gwdiff_sec}

\subsection{Summary}
\label{gwdiffsumm_subs}

\noindent
Suppose $(X,\om)$, $A$, $k$, and $J$ are as in the statement of Theorem~\ref{gwdiff_thm}
and
$$\psi\in H^{\dim_{1,k}(X,A)}\big(\ov\M_{1,k}(X,A;J);\Q\big)$$
is a geometric cohomology class. By definition,
$$\psi=\prod_{l=1}^{l=k}\ev_l^*\mu_l
\qquad\hbox{for some}\qquad 
\mu_l\in H^{2n-d_l}(X;\Z),~~d_l<2n.$$
For each $l\!\in\![k]$, choose a pseudocycle representative 
$$f_l\!:\bar{Y}_l\lra X$$
for $\PD_X\mu_l$. 
In particular, $\bar{Y}_l$ is a disjoint union of smooth manifolds.
The dimension of one of them, $Y_{l;\mn}$, is~$d_l$; the dimensions of all others
are at most $d_l\!-\!2$.
The map $f_l$ is continuous, and its restriction to each smooth manifold is smooth;
see Chapter~7 in~\cite{McSa} or Section~1 in~\cite{RT}.
Let 
\begin{gather*}
\ev\!=\!\prod_{l=1}^{l=k}\ev_l\!:\X_{1,k}(X,A) \lra \prod_{l=1}^{l=k}X,\\ 
f\!=\!\prod_{l=1}^{l=k}f_l\!: \bar{Y}\!\equiv\!\prod_{l=1}^{l=k}\bar{Y}_l
\lra\prod_{l=1}^{l=k}X,
\qquad\hbox{and}\qquad  Y_{\mn}\!=\!\prod_{l=1}^{l=k}Y_{l;\mn}.
\end{gather*}\\

\noindent
With $(f_l)_{l\in[k]}$ as above, 
for any $\nu\!\in\!\G_{1,k}^{0,1}(X,A;J)$  and a bubble type $\T$ as 
in Subsection~\ref{notation1_subs}, let
\begin{alignat*}{1}
\X_{1,k}(X,A;f) &= \big\{(b,z)\!\in\!\X_{1,k}(X,A)\!\times\!\bar{Y}:
\ev(b)\!=\!f(z)\big\},\\
\ov\M_{1,k}(X,A;J,\nu;f) &= 
\big(\ov\M_{1,k}(X,A;J,\nu)\!\times\!\bar{Y}\big) \cap \X_{1,k}(X,A;f),
\qquad\hbox{and}\\
\U_{\T,\nu}(X;J;f) &= \big(\U_{\T,\nu}(X;J)\!\times\!\bar{Y}\big)\cap\X_{1,k}(X,A;f).
\end{alignat*}
If $\nu$ is sufficiently small, the space $\ov\M_{1,k}(X,A;J,\nu;f)$ is compact.\
Let
$$\De_X^k=\prod_{l=1}^{l=k}\De_X\subset\prod_{l=1}^{l=k} X\!\times\!X$$
be the $k$-fold product of the diagonals.
If $\nu$ and $f_l$ are chosen generically, then 
\begin{equation}\label{transver_e1a}
\ov\M_{1,k}(X,A;J,\nu;f) \subset  \M_{1,k}^0(X,A;J,\nu)\!\times\!Y_{\mn}
\end{equation}
and the map
$$\ev\!\times\!f\!:  \M_{1,k}^0(X,A;J,\nu)\!\times\!Y_{\mn} \lra (X^2)^k$$
is transverse to $\De_X^k$. 
Thus, $\ov\M_{1,k}(X,A;J,\nu;f)$ is a compact zero-dimensional orbifold.
By definition,
\begin{equation}\label{transver_e1b}
\GW_{1,k}^X(A;\psi)=~ ^{\pm}\!\big|\ov\M_{1,k}(X,A;J,\nu;f)\big|,
\end{equation}
if $\nu$ is sufficiently small and $\nu$ and $f_l$ are generic.\\

\noindent
Similarly to the previous paragraph, if $\nu_{\es}\!\in\!\G_{1,k}^{\es}(X,A;J)$, let
$$\ov\M_{1,k}^0(X,A;J,\nu_{\es};f) =
\big(\ov\M_{1,k}^0(X,A;J,\nu_{\es})\!\times\!\bar{Y}\big) 
\cap \ov\M_{1,k}(X,A;J,\nu_{\es};f).$$
For generic $\nu_{\es}\!\in\!\G_{1,k}^{\es}(X,A;J)$ and $f_l$,
\begin{equation}\label{transver_e2a}
\ov\M_{1,k}^0(X,A;J,\nu_{\es};f) \subset  \M_{1,k}^0(X,A;J,\nu_{\es})\!\times\!Y_{\mn}
\end{equation}
and the map
$$\ev\!\times\!f\!:  \M_{1,k}^0(X,A;J,\nu)\!\times\!Y_{\mn} \lra (X^2)^k$$
is transverse to $\De_X^k$, by the second half of Subsection~\ref{res_subs}.
Thus, $\ov\M_{1,k}^0(X,A;J,\nu_{\es};f)$ is a compact zero-dimensional 
orbifold.
By definition,
\begin{equation}\label{transver_e2b}
\GW_{1,k}^{0;X}(A;\psi)=~ ^{\pm}\!\big|\ov\M_{1,k}^0(X,A;J,\nu_{\es};f)\big|,
\end{equation}
if $\nu_{\es}$ is sufficiently small and $\nu_{\es}\!\in\!\G_{1,k}^{\es}(X,A;J)$ 
and $f_l$ are generic.\\

\noindent
If $\nu_{\es}\!\in\!\G_{1,k}^{\es}(X,A;J)$  and $f_l$ are generic,
\e_ref{transver_e1a} and \e_ref{transver_e1b} do not generally hold for 
$\nu\!=\!\nu_{\es}$ because the restriction of the bundle section
$\bar\partial\!+\!\nu_{\es}$ 
is not transverse to the zero set along some strata 
$$\X_{\T}(X,A)\subset \X_{1,k}(X,A)-\X_{1,k}^{\{0\}}(X,A).$$
Instead, we will apply \e_ref{transver_e1a} and \e_ref{transver_e1b}
with $\nu$ replaced by $\nu_{\es}\!+\!\nu$ for a generic $\nu\!\in\!\G_{1,k}^{0,1}(X,A;J)$ 
which is sufficiently small relatively to~$\nu_{\es}$.
In such a case, the compact zero-dimensional orbifold 
\begin{equation}\label{numoduli_e}
\ov\M_{1,k}\big(X,A;J,\nu\!+\!\nu_{\es};f\big)
\end{equation}
will lie in a small neighborhood of 
$$\ov\M_{1,k}(X,A;J,\nu_{\es};f) \subset \X_{1,k}(X,A;f).$$
We will express the cardinality of this orbifold in terms of data intrinsic 
to $\ov\M_{1,k}(X,A;J,\nu_{\es};f)$.
From the transversality of the relevant maps, 
it is straightforward to see that there is a unique element of~\e_ref{numoduli_e}
close to each of the elements of $\ov\M_{1,k}^0(X,A;J,\nu_{\es},f)$.
Thus,
$$\GW_{1,k}^X(A;\psi)- \GW_{1,k}^{0;X}(A;\psi)$$
is the number of elements of~\e_ref{numoduli_e} that lie close to the closed subset 
\begin{equation}\label{boundspace_e}
\ov\M_{1,k}(X,A;J,\nu_{\es};f)-\ov\M_{1,k}^0(X,A;J,\nu_{\es})\!\times\!\bar{Y}
\end{equation}
of $\ov\M_{1,k}(X,A;J,\nu_{\es})\!\times\!\bar{Y}$.
The contribution of \e_ref{boundspace_e} to~\e_ref{transver_e2b} can in fact be split 
into contributions of the subspaces $\U_{\T,\nu_{\es}}(X;J;f)$ of~\e_ref{boundspace_e}.
By studying local obstructions similarly to~\cite{g2n2and3},
each of these contributions will be shown to be equal to the number of zeros of an affine
bundle map between two finite-rank vector bundles over $\U_{\T,\nu_{\es}}(X;J;f)$;
see Proposition~\ref{bdcontr_prp} below.
Such numbers can be determined using the procedure described in 
\cite[Subsections 3.2,3.3]{g2n2and3}.\\

\noindent
If $\dim_{\R}X\!=\!4$, only two strata of~\e_ref{boundspace_e} are nonempty for a good
choice of~$\nu_{\es}$. 
They are isomorphic~to
$$\cM_{1,1}\times\M_{0,\{0\}\cup[k]}^0(X,A;J,\nu_B;f) \qquad\hbox{and}\qquad
\big(\ov\cM_{1,1}\!-\!\cM_{1,1}\big)\times\M_{0,\{0\}\cup[k]}^0(X,A;J,\nu_B;f)$$
for some $\nu_B\!\in\!\G_{0,\{0\}\cup[k]}^{0,1}(X,A;J)$, where
$$\cM_{1,1} \subset \ov\cM_{1,1}$$
is the complement of the equivalence class of the singular elliptic curve.
The dimension of\linebreak 
$\M_{0,\{0\}\cup[k]}^0(X,A;J,\nu_B;f)$ is zero,
even though no constraint has been imposed on the zeroth marked point.
In other words, $\M_{0,\{0\}\cup[k]}^0(X,A;J,\nu_B;f)$ is ``virtually empty".
It is thus not too surprising that neither of these strata contributes 
to~$\GW_{1,k}^X(A;\psi)$.\\

\noindent
{\it Remark:} If $J$ is a genus-one $A$-regular almost complex structure 
in the sense of  \cite[Definition~1.3]{g1comp}, 
we can take $\nu_{\es}\!=\!0$. 
If $\dim_{\R}X\!=\!4$, we then find that the space~\e_ref{boundspace_e} is empty, 
since 
$$\M_{0,\{0\}\cup[k]}^0(X,A;J;f)\equiv\M_{0,\{0\}\cup[k]}^0(X,A;J,0;f)$$
cannot be zero-dimensional.
Thus, if $(X,\om)$ admits a genus-one $A$-regular almost complex structure,
the first case of Theorem~\ref{gwdiff_thm} is immediate from dimension-counting,
once it is known that $\GW_{1,k}^{0;X}(A;\psi)$ is well defined.\\

\noindent
If $\dim_{\R}X\!=\!6$, for a good choice of~$\nu_{\es}$ only a few
strata of~\e_ref{boundspace_e} are nonempty.
All, but two of them, are either virtually empty or $\bpar_J$-hollow, 
in the sense of \cite[Subsection~3.1]{g2n2and3}.
In either of these cases, $\U_{\T,\nu_{\es}}(X;J;f)$ does not contribute
to~\e_ref{transver_e1b}.
The two remaining strata are isomorphic~to
$$\cM_{1,1}\times\M_{0,\{0\}\cup[k]}^0(X,A;J,\nu_B;f) \qquad\hbox{and}\qquad
\cM_{1,2}\times\M_{0,k}^0(X,A;J,\nu_B;f)$$
for some $\nu_B\!\in\!\G_{0,\{0\}\cup[k]}^{0,1}(X,A;J)$ or
$\nu_B\!\in\!\G_{0,k}^{0,1}(X,A;J)$, respectively.
In the second case, one of the elements of $[k]$ corresponds to the attaching node
of the only bubble component of each element in $\U_{\T,\nu_{\es}}(X;J;f)$;
we will denote it by~$0$.
Let
$$L_1\lra \ov\cM_{1,1},\,\ov\cM_{1,2} \qquad\hbox{and}\qquad
L_0\lra\ov\M_{0,\{0\}\sqcup k}(X,A;J,\nu_B;f),\,\M_{0,k}^0(X,A;J,\nu_B;f)$$
be the universal tangent line bundles for the marked points labeled by $1$ and $0$,
respectively.
Let
$$\cD_0\!: L_0 \lra \ev_0^*TX$$
be the bundle homomorphism over $\ov\M_{0,\{0\}\sqcup k}(X,A;J,\nu_B;f)$ or
$\M_{0,k}^0(X,A;J,\nu_B;f)$ given~by
$$\cD_0\big([b,v]\big)=du_b|_{y_0(b)}v.$$
We denote~by
\begin{alignat*}{1}
&\pi_P,\pi_B\!: \ov\cM_{1,1}\!\times\!\ov\M_{0,\{0\}\cup[k]}(X,A;J,\nu_B;f)
\lra \ov\cM_{1,1},\,\ov\M_{0,\{0\}\cup[k]}(X,A;J,\nu_B;f)
\qquad\hbox{and}\\
&\pi_P,\pi_B\!: \ov\cM_{1,2}\!\times\!\M_{0,k}^0(X,A;J,\nu_B;f)
\lra \ov\cM_{1,2},\,\M_{0,k}^0(X,A;J,\nu_B;f)
\end{alignat*}
the two projection maps.
In both cases, the linear part of the affine map determining the contribution of
the stratum $\U_{\T,\nu_{\es}}(X;J;f)$ to~\e_ref{transver_e1b}~is
\begin{gather*}
\cD_{\T}\!:\pi_P^*L_1\!\otimes\!\pi_B^*L_0\lra \pi_P^*\E^*\!\otimes\!\pi_B^*\ev_0^*TX,\\
\big\{\cD_{\T}[b_P,b_B,v_P\!\otimes\!v_B]\big\}\big([b_B,\psi])=
\psi_{x_1(b_P)}(v_P)\cdot_J\!\cD_0\big([b_B,v_B]\big),\\
\hbox{if}\qquad [b_P,b_B,v_P\!\otimes\!v_B]\!\in\!\pi_P^*L_1\!\otimes\!\pi_B^*L_0,\,
[b_B,\psi]\!\in\!\pi_P^*\E.
\end{gather*}
The constant term $\bar\nu$ of each of the affine maps is generic.
In the second case, $\M_{0,k}^0(X,A;J,\nu_B;f)$ is a finite collection of points.
It is then straightforward to see that for a generic $\nu$, the affine bundle map
$\cD_{\T}\!+\!\bar\nu$ does not vanish.
Thus, the corresponding stratum $\U_{\T,\nu_{\es}}(X;J;f)$ does not contribute
to~\e_ref{transver_e1b}.
We will show in Subsection~\ref{affinemap_subs} that 
the number $N(\cD_{\T})$ of zeros of $\cD_{\T}\!+\!\bar\nu$ in the first case~is
$$N(\cD_{\T})=\frac{2\!-\!\lr{c_1(TX),A}}{24} \, \GW_{0,k}^X(A;\psi),$$
proving the second case of Theorem~\ref{gwdiff_thm}.\\

\noindent
{\it Remark:} If $J$ is a genus-one $A$-regular almost complex structure and 
$\dim_{\R}X\!=\!6$, the space \e_ref{boundspace_e} is the union of the spaces
\begin{equation}\label{rigidcurve_e}
\ov\M_{1,k}\big(\ka,1;J|_{\ka};\under{f}|_{\ka}\big)
\end{equation}
taken over all degree-$A$ genus-zero curves $\ka$ in $X$ that intersect $f_l(Y_{l;\mn})$
for every $l\!\in\![k]$.
Based on~\cite{P}, one would expect that each of the spaces~\e_ref{rigidcurve_e}
contributes $\frac{(2-\lr{c_1(TX),A}}{24}$ to~$\GW_{1,k}^X(A;\psi)$.
The total number of the curves~$\ka$ is~$\GW_{0,k}^X(A;\psi)$.
In particular, the second case of Theorem~\ref{gwdiff_thm}, just like the first,
is consistent with geometric expectations.

\subsection{Analytic Setup}
\label{setup_subs}

\noindent
Let $\T\!=\!([k],I,\aleph;j,\under{A})$ be a bubble type such that
\begin{equation}\label{bubbtype_e}
A_i=0~~\forall\,i\!\in\!I_0 \qquad\hbox{and}\qquad
\sum_{i\in I}\!A_i=A.
\end{equation}
For each $i\!\in\!I\!-\!I_0$, let
$$H_i\T=\big\{h\!\in\!\hat{I}\!:\io_h\!=\!i\big\} \qquad\hbox{and}\qquad
M_i\T=\big\{l\!\in\![k]\!: j_l\!=\!i\big\}.$$
We denote by
$$\pi_{\T;i}\!: \X_{\T}(X) \lra \X_{0,\{0\}\sqcup H_i\T\sqcup M_i\T}(X,A_i)$$
the map sending each element $[b]$ of $\X_{\T}(X)$ to its restriction to $\Si_{b,i}$:
$$\big[(S,[k],I,\aleph;x,(j,y),u)\big] \lra 
\big[(H_i\T\!\sqcup\!M_i\T,\{i\};,(\io,x)|_{H_i\T}\!\sqcup\!(j,y)|_{M_i\T},u_i)\big].$$
Let $\G_{1,k}^{\gd}(X,A;J)$ be the subspace of elements $\nu$ in 
$\G_{1,k}^{\es}(X,A;J)$ such that for every bubble type $\T$ as above, 
$$\nu|_{\X_{\T}(X)}=\sum_{i\in I-I_0}\pi_{\T;i}^*\nu_{\T;i} \qquad\hbox{for some}\qquad
\nu_{\T;i}\in\G_{0,\{0\}\sqcup H_i\T\sqcup M_i\T}^{0,1}(X,A_i;J)$$
and for every $[b]\!\in\!\M_{0,\{0\}\sqcup H_i\T\sqcup M_i\T}^0(X,A_i;J,\nu_{\T;i})$ 
the linearization
\begin{equation}\label{setup_e1b}
D_{J,\nu_{\T;i};b}\!: \big\{\xi\!\in\!\Ga(\Si_b;u_b^*TX)\!: \xi(y_0(b))\!=\!0\big\} 
\lra \Ga\big(\Si_b;\La^{0,1}_{J,j}T^*\Si_b\!\otimes\!u_b^*TX\big)
\end{equation}
of the section $\bpar_J\!+\!\nu_{\T;i}$ at $b$ is surjective.\\

\noindent
For a generic element $\nu_{\T;i}\!\in\!\G_{0,\{0\}\sqcup H_i\T\sqcup M_i\T}^{\es}(X,A_i;J)$,
the operator~\e_ref{setup_e1b} is surjective for every 
$$[b]\in\ov\M_{0,\{0\}\sqcup H_i\T\sqcup M_i\T}(X,A_i;J,\nu_{\T;i}).$$
This implies that the closure of $\G_{1,k}^{\gd}(X,A;J)$ in $\G_{1,k}^{0,1}(X,A;J)$
contains the zero section, since
we can construct an element $\nu$ of $\G_{1,k}^{\gd}(X,A;J)$ inductively 
starting from the highest-codimension strata of $\X_{1,k}(X,A)$.
If $\T$ is a bubble type as above and $\nu$ has been defined on 
$$\ov\X_{\T}(X)-\X_{\T}(X)\subset \X_{1,k}(X,A)$$
subject to the above restriction and regularity conditions, 
then~$\nu$ induces a multisection $\nu_{\T;i}$~of 
$$\Ga_{0,\{0\}\sqcup H_i\T\sqcup M_i\T}^{0,1}(X,A_i;J) \lra 
\X_{0,\{0\}\sqcup H_i\T\sqcup M_i\T}(X,A_i)-\X_{0,\{0\}\sqcup H_i\T\sqcup M_i\T}^0(X,A_i).$$
It extends continuously to an effectively supported multisection over all of 
$\X_{0,\{0\}\sqcup H_i\T\sqcup M_i\T}(X,A_i)$.
By perturbing this extension outside of 
$$\X_{0,\{0\}\sqcup H_i\T\sqcup M_i\T}(X,A_i)
 -\X_{0,\{0\}\sqcup H_i\T\sqcup M_i\T}^0(X,A_i),$$
we obtain an element $\nu_{\T;i}$ of $\G_{0,\{0\}\sqcup H_i\T\sqcup M_i\T}^{\es}(X,A_i;J)$
such that the operator~\e_ref{setup_e1b} is surjective for every 
$$[b]\in\ov\M_{0,\{0\}\sqcup H_i\T\sqcup M_i\T}(X,A_i;J,\nu_{\T;i}).$$\\

\noindent
We fix small generic elements 
$$\nu_{\es}\in\G_{1,k}^{\gd}(X,A;J) \qquad\hbox{and}\qquad
\nu\in\G_{1,k}^{0,1}(X,A;J)$$
such that for all $t\!\in\!\R^+$ sufficiently small the section 
$$\big\{\bpar_J\!+\!\nu_{\es}\!+\!t\nu\big\}\big|_{\X_{\T}(X)}$$
is transverse to the zero set in $\Ga_{1,k}^{0,1}(X,A;J)|_{\X_{\T}(X)}$
for every stratum $\X_{\T}(X)$ of~$\X_{1,k}(X,A)$.
Let $\{Y_{\la}\}_{\la\in\A}$ be the strata of $\bar{Y}$ induced by 
the partitions of each $\bar{Y}_l$ into smooth manifolds.
By our assumptions,
$$\dim_{\R}\!Y_{\mn}=2nk-\dim_{1,k}(X,A) \qquad\hbox{and}\qquad
\dim_{\R}\!Y_{\la}\le \dim_{\R}\!Y_{\mn}\!-\!2~~~\forall\,\la\!\in\!\A\!-\!\{0\}.$$
For each bubble type $\T$ as above and $\la\!\in\!\A$, let
$$\U_{\T,\nu_{\es}}(X;J;f_{\la})=
\big(\U_{\T,\nu_{\es}}(X;J)\!\times\!Y_{\la}\big)\cap
\ov\M_{1,k}(X,A;J,\nu_{\es};f).$$
We will call a bubble type $\T$ as above \sf{simple} if
$$\aleph=\eset \qquad\hbox{and}\qquad \hat{I}=\chi(\T).$$
In other words, $\T$ is simple if and only if for every element 
$[b]\!\in\!\U_{\T,\nu_{\es}}(X,A;J)$ the domain $\Si_b$ consists of 
a smooth principal component~$\Si_{b;P}$, 
on which the map $u_b$ is necessarily constant, and $|\hat{I}|$ 
bubble components, all of which are attached directly to~$\Si_{b;P}$
and on which the map $u_b$ is not constant. \\

\noindent
Suppose $\ov\cM$ is a compact topological space which is a disjoint union of 
smooth orbifolds, one of which, $\cM$, is a dense open subset of $\ov\cM$,
and the dimensions of all others do not exceed $\dim\cM\!-\!2$.
Let 
$$E,\O\lra \ov\cM$$
be vector bundles such that the restrictions of $E$ and $\O$ to every stratum of $\ov\cM$
is smooth~and
$$\rk\,\O-\rk\,E=\frac{1}{2}\dim_{\R}\!\cM.$$
If 
$$\al\in\Ga\big(\ov\cM;\Hom(E,\O)\big)$$
is a regular section in the sense of  \cite[Definition~3.9]{g2n2and3},
then the cardinality of the zero set of the affine bundle map
$$\psi_{\al,\bar\nu}\!\!\equiv\!\al\!+\!\bar\nu\!: E\lra \O$$
is finite and independent of a generic choice of $\bar\nu\!\in\!\Ga(\ov\cM;\O)$,
by  \cite[Lemma~3.14]{g2n2and3}.
We denote it by $N(\al)$.
A key step in our proof of Theorem~\ref{gwdiff_thm} is the following proposition.

\begin{prp}
\label{bdcontr_prp}
Suppose $(X,\om,J)$, $A$, $k$, $\psi$, $f_l$, $\nu_{\es}$, $\nu$, and  
$\{Y_{\la}\}_{\la\in\A}$ are as above.
If $\T$ is a bubble as above and $\la\!\in\!\A$, 
there exist 
$$\cC_{\U_{\T,\nu_{\es}}(X;J;f_{\la})}(\bpar_J)\in\Q, \qquad
\ep_{\nu}\in\R^+,$$
and a compact subset $K_{\nu}$ of $\U_{\T,\nu_{\es}}(X;J;f_{\la})$
with the following property.
For every compact subset $K$ of $\U_{\T,\nu_{\es}}(X;J;f_{\la})$ and 
open subset $U$ of $\X_{1,k}(X,A;f)$,
there exist an open neighborhood $U_{\nu}(K)$ of $K$ in $\X_{1,k}(X,A;f)$
and $\ep_{\nu}(U)\!\in\!(0,\ep_{\nu})$, respectively, such that
$$^{\pm}\big|\ov\M_{1,k}\big(X,A;J,\nu_{\es}\!+\!t\nu;f\big)\!\cap\!U\big|
=\cC_{\U_{\T,\nu_{\es}}(X;J;f_{\la})}(\bpar_J)
\quad\hbox{if}~~
t\!\in\!(0,\ep_{\nu}(U)),~K_{\nu}\!\subset\! K\!\subset\! U\!\subset\! U_{\nu}(K).$$
Furthermore, if $\T$ is simple and $\la\!=\!\mn$,
$$\cC_{\U_{\T,\nu_{\es}}(X;J;f_{\la})}(\bpar_J)=N(\cD_{\T})$$
for some regular vector bundle homomorphism $\cD_T$ over $\U_{\T,\nu_{\es}}(X;J;f)$.
Otherwise, 
$$\cC_{\U_{\T,\nu_{\es}}(X;J;f_{\la})}(\bpar_J)=0.$$\\
\end{prp}

\noindent
This proposition is proved in the next subsection.
In the previous subsection we described the homomorphism $\cD_{\T}$ for a simple
bubble type $\T$ such that $|\chi(\T)|\!=\!1$.
Below we describe this homomorphism for an arbitrary bubble $\T$ 
satisfying~\e_ref{bubbtype_e}.\\

\noindent
For each $i\!\in\!\hat{I}$, let
$$\ti{H}_i\T=\big\{h\!\in\!\hat{I}\!: \io_h\!\ge\!i\big\}, \qquad
\ti{M}_i\T=\big\{l\!\in\![k]\!: j_l\!\ge\!i\big\}, \quad\hbox{and}\quad
\ti{A}_i=\sum_{h\ge i}A_h.$$
We denote by
$$\ti\pi_{\T;i}\!: \X_{\T}(X) \lra \X_{0,\{0\}\sqcup\ti{M}_i\T}(X,\ti{A}_i)$$
the map sending each element $[b]$ of $\X_{\T}(X)$ to its restriction to 
the tree of bubble components beginning with~$\Si_{b,i}$:
$$\big[(S,[k],I,\aleph;x,(j,y),u)\big] \lra 
\big[(\ti{M}_i\T,\{i\}\!\cup\!\ti{H}_i\T;(\io,x)|_{\ti{H}_i\T},(j,y)|_{\ti{M}_i\T},
u|_{\{i\}\cup\ti{H}_i\T})\big].$$
By our assumptions on $\nu_{\es}$,
$$\ti\pi_{\T;i}\!:  \U_{\T,\nu_{\es}}(X;J)\lra 
\ov\M_{0,\{0\}\sqcup\ti{M}_i\T}\big(X,\ti{A}_i;J,\ti\nu_{B;i}\big)$$
for some $\ti\nu_{B;i}\!\in\!\G_{0,\{0\}\sqcup\ti{M}_i\T}^{\es}(X,\ti{A}_i;J)$.
Let 
$$\F\T=\bigg(\bigoplus_{i\in\chi(\T)}\!\!\pi_P^*L_{h(i)}\!\otimes\!\ti\pi_{\T;i}^*L_0
\bigg)\Big/\Aut(\T) \lra \U_{\T,\nu_{\es}}(X;J),\,\U_{\T,\nu_{\es}}(X;J;f).$$\\

\noindent
If $M_0$ is a finite set and $h\!\in\!M_0$, let 
$$s_h\in\Ga\big(\ov\cM_{1,M_0};\Hom(L_h,\E^*)\big)$$
be the the section given by
$$\big\{s_h(b;v)\big\}(\psi)=\psi_{x_h(b)}v\in\C \qquad\hbox{if}\quad
v\in L_h|_b,~\psi\in\E_b,$$
where $x_h(b)\!\in\!\Si_b$ is the $h$th marked point.
We define the homomorphism
$$\cD_{\T}\!:\F\T\lra \big(\pi_P^*\E^*\!\otimes\!\ev_P^*TX\big)\big/\Aut(\T)$$
over $\U_{\T,\nu_{\es}}(X;J)$ or $\U_{\T,\nu_{\es}}(X;J;f)$ by
$$\cD_{\T}=\sum_{i\in\chi(\T)}\!\pi_P^*s_{h(i)}\!\otimes\!\ti\pi_{\T;i}^*\cD_0.$$

\subsection{Proof of Proposition~\ref{bdcontr_prp}}
\label{bdcontrprp_subs}

\noindent
We continue with the notation of Subsections~\ref{g0prp_subs} and~\ref{setup_subs}.
By our assumptions on~$\nu_{\es}$, the operators $D_{J,\nu_{\es};b,i}$ are surjective 
for all $[b]\!\in\!\U_{\T,\nu_{\es}}(X;J)$ and $i\!\in\!\hat{I}$.
Thus, we can take
$$\ti\Ga_-^{0,1}(b;i) = \{0\} \qquad\forall\,
[b]\!\in\!\U_{\T,\nu_{\es}}(X;J),\, i\!\in\!I.$$
The corresponding space $\ti\U$ of Subsection~2.4 
is a smooth manifold of $(J,\nu_{\es})$-holomorphic maps.
In Subsections~2.5 of~\cite{gluing} and~4.2 of~\cite{g1comp}, 
we describe a space $\U_{\T}^{(0)}(X;J)$ of balanced $J$-holomorphic maps, 
not of equivalence classes of such maps.
If $\nu_{\es}$ is sufficiently small, the same definitions can be used to describe
a submanifold $\U_{\T,\nu_{\es}}^{(0)}(X;J)$ of~$\ti\U$.
In particular,
$$ \U_{\T,\nu_{\es}}(X;J)=\U_{\T,\nu_{\es}}^{(0)}(X;J)\big/
\Aut(\T)\!\propto\!(S^1)^{\hat{I}}$$
for a natural action of $\Aut(\T)$ on $(S^1)^{\hat{I}}$ and of 
$\Aut(\T)\!\propto\!(S^1)^{\hat{I}}$ on~$\ti\U$.
Let
$$\ti\cF\T=\wt\cF|_{\U_{\T,\nu_{\es}}^{(0)}(X;J)},$$
where $\wt\cF\!\lra\!\ti\U$ is the vector bundle defined in Subsection~\ref{comp1prp_subs}.
The above group action on $\ti\U$ lifts to an action on $\wt\cF$ so~that
$$\cF\T \equiv \ti\cF\T \big/ \Aut(\T)\!\propto\!(S^1)^{\hat{I}}$$
is the bundle of gluing parameters for $\U_{\T,\nu_{\es}}(X;J)$.\\

\noindent
We will apply the construction of Subsection~\ref{comp1prp_subs},
with some refinements, to the entire space $\U_{\T,\nu_{\es}}^{(0)}(X;J)$,
instead of a small open subset of~$\ti\U$.
We will view $\R$-valued functions on $\U_{\T,\nu_{\es}}(X;J)$
as functions on $\U_{\T,\nu_{\es}}^{(0)}(X;J)$ via the quotient projection~map
$$\U_{\T,\nu_{\es}}^{(0)}(X;J) \lra \U_{\T,\nu_{\es}}(X;J).$$
Fix small 
$$\de,\ep \in C^{\i}(\U_{\T,\nu_{\es}}(X;J);\R^+)$$ 
such that the basic gluing map
$$q_{\ups}\!: \Si_{\ups}\lra \Si_b$$
of  \cite[Subsection~2.2]{gluing} and the modified gluing map
$$\ti{q}_{\ups_0;2}\!: \Si_{\ups}\lra \Si_{\ups_1}$$
of  \cite[Subsection~4.2]{g1comp} with the collapsing radius $\ep(b)$
are defined for all $\ups\!\equiv\!(b,v)\!\in\!\ti\cF\T_{\de}$.
For all 
$$\xi\in\Ga(\ups_1) \qquad\hbox{s.t.}\qquad \|\xi\|_{\ups_1,p,1}\le\de(b),$$
let $u_{\ups_1,\xi}$, $u_{\ups,\xi}$, $\Ga_+(\ups;\xi)$,
$$D_{J,\nu_{\es};\ups;\xi}\!: \Ga(\ups;\xi) \lra \Ga^{0,1}(\ups;\xi), \qquad
R_{\ups,\xi}\!: \Ga_-^{0,1}(b;J)\lra\Ga^{0,1}(\ups;\xi),$$
and $\nu_{\es;\ups_1,\xi}$ be as in Subsection~\ref{comp1prp_subs}, 
with $J_t$, $\nu_t$, and  $\na^t$ replaced
by~$J$, $\nu_{\es}$, and a $J$-compatible connection~$\na$, respectively.
The estimates \e_ref{g1reg_lmm1e6}, \e_ref{torus_e5}, \e_ref{torus_e6b},
and \e_ref{torus_e7b}-\e_ref{torus_e8b} continue to hold if $C\!\in\!\R^+$, 
$D_{\ups;\xi}^t$, $J_t$,  and $\nu_t$ 
are replaced by 
$$C\in C^{\i}(\U_{\T,\nu_{\es}}(X;J);\R^+),$$
$D_{J,\nu_{\es};\ups;\xi}$, $J$, and $\nu_{\es}$, respectively.
In \e_ref{torus_e8a} and~\e_ref{torus_e8b}, $\ti\pi_{\ups;\xi}^{t;+}\!=\!\id$
and  $J_0\!=\!J$.\\

\noindent
With notation as in Subsection~\ref{comp1prp_subs}, for each 
$\ups\!\in\!\ti\cF\T_{\de}^{\eset}$ let
$$\Ga_{B;-}(\ups_1)=\big\{\xi\!\circ\!q_{\ups_1}\!: \xi\!\in\!\ker D_{J,\nu_{\es};b},\,
\xi|_{\Si_{b;P}}\!=\!0\big\}.$$
We denote by $\Ga_{B;+}(\ups_1)$ the $L^2$-orthogonal complement of $\Ga_{B;-}(\ups_1)$
in~$\Ga_B(\ups_1)$.
Let 
$$D_{J,\nu_{\es};\ups_1}^B\!:\Ga_B(\ups_1) \lra \Ga_B^{0,1}(\ups_1)$$
be the linearization of section $\bpar_J\!+\!\nu_{\es}$ at 
$(\Si_{\ups_1};u_{\ups_1})$ defined via the connection~$\na$.
Similarly to~\e_ref{torus_e4},
\begin{equation}\label{gwdiff_e3}
C(b)^{-1}\|\xi\|_{\ups_1,p,1}\le \big\|D_{J,\nu_{\es};\ups_1}^B\xi\big\|_{\ups_1,p}
\le C(b)\|\xi\|_{\ups_1,p,1} 
\quad \forall~ \ups\!=\!(b,v)\!\in\!\ti\cF\T_{\de}^{\eset}, ~
~\xi\!\in\!\Ga_{B;+}(\ups_1),
\end{equation}
for some $C\!\in\!C^{\i}(\U_{\T,\nu_{\es}}(X;J);\R^+)$, provided that
$\de\!\in\!C^{\i}(\U_{\T,\nu_{\es}}(X;J);\R^+)$ is sufficiently small.
In particular, the operator 
$$D_{J,\nu_{\es};\ups_1}^B\!:\Ga_{B;+}(\ups_1) \lra \Ga_B^{0,1}(\ups_1)$$
is an isomorphism. Its norm and the norm of its inverse are dependent
only on $[b]\!\in\!\U_{\T;\nu_{\es}}(X;J)$.
Thus, by the Contraction Principle, for each $\ups\!\in\!\ti\cF\T_{\de}^{\eset}$, 
the equation 
$$\big\{\bpar_J\!+\!\nu_{\es;\ups_1,\xi}\big\}u_{\ups_1,\xi}=0,
\qquad\xi\in\Ga_{B;+}(\ups_1),$$
has a unique small solution $\xi_{\nu_{\es}}(\ups_1)$.\\

\noindent
{\it Remark:} Since $\nu_{\es}$ is a multisection, the uniqueness statement above,
as well as similar statements below, should be interpreted in terms of local 
branches of~$\nu_{\es}$ as defined in \cite[Section~3]{FuO}.

\begin{lmm}
\label{gwdiff_lmm1}
If $\T$ is a bubble type as in~\e_ref{bubbtype_e} and 
$\nu_{\es}\!\in\!\G_{1,k}^{\gd}(X,A;J)$ is a sufficiently small generic perturbation,
there exist 
$$\de\in C^{\i}\big(\U_{\T,\nu_{\es}}(X;J);\R^+\big)$$
and an open neighborhood $U_{\T}$ of $\U_{\T,\nu_{\es}}(X;J)$ in $\X_{1,k}(X,A)$ 
such that the~map
\begin{gather*}
\begin{split}
\big\{(\ups,\ze)\!: \ups\!=\!(b,v)\!\in\!\ti\cF\T_{\de}^{\eset};
\ze\!\in\!\Ga_+(\ups;\xi_{\nu_{\es}}(\ups_1)),\,
\|\ze\|_{\ups,p,1}\!<\!\de(b)\big\}\big/ \Aut(\T)\!\propto\!(S^1)^{|\hat{I}|}  \qquad\qquad&\\
\lra \X_{1,k}^0(X,A)\!\cap\!U_{\T}, &
\end{split} \\
\big[(\ups,\ze)\big]\lra \big[(\Si_{\ups},\exp_{u_{\ups,\xi_{\nu_{\es}}(\ups_1)}}\!\ze)\big],
\end{gather*}
is a diffeomorphism.
\end{lmm}

\noindent 
It is immediate from the construction that the map
$$\big(\ups,\ze\big) \lra 
\big[(\Si_{\ups},\exp_{u_{\ups,\xi_{\nu_{\es}}(\ups_1)}}\!\ze)\big]$$
is $\Aut(\T)\!\propto\!(S^1)^{|\hat{I}|}$-invariant and smooth.
The injectivity and surjectivity of the induced map on the quotient are proved
by arguments similar to Subsections~4.2 and 4.3-4.5 in~\cite{gluing}, 
respectively; see also the paragraph following Lemma~4.4
in~\cite{g1comp}.\\

\noindent
For each $\ups\!=\!(b,v)\!\in\!\ti\cF\T_{\de}^{\eset}$, we define the homomorphism 
$$\pi_{\ups;-}^{0,1}\!: \Ga^{0,1}\big(\ups;\xi_{\nu_{\es}}(\ups_1)\big) \lra 
\Ga_-^{0,1}(b;J)\!\approx\!\E_{\pi_P(b)}^*\!\otimes\!T_{\ev_P(b)}X$$
as follows. 
If $\{\eta_r\}_{r\in[n]}$ is an orthonormal basis for $\Ga_-^{0,1}(b;J)$,
we~put
$$\pi_{\ups;-}^{0,1}\eta'=\sum_{r=1}^{r=n}
\bllrr{\eta',R_{\ups,\xi_{\nu_{\es}}(\ups_1)}\eta_r}\eta_r
\qquad\forall\,\eta'\!\in\!\Ga^{0,1}\big(\ups;\xi_{\nu_{\es}}(\ups_1)\big).$$
This map is well defined.\\

\noindent
For each $(\ups,\ze)$ as above, let
$$\Pi_{\ze}\!: L^p(\Si_{\ups};\La^{0,1}_{J,j}T^*\Si_{\ups}\!\otimes\!
u_{\ups,\xi_{\nu_{\es}}(\ups_1)}^*TX)\lra
L^p(\Si_{\ups};\La^{0,1}_{J,j}T^*\Si_{\ups}\!\otimes\!
\{\exp_{u_{\ups,\xi_{\nu_{\es}}(\ups_1)}}\!\ze\}^*TX)$$
be the isomorphism induced by the $\na$-parallel transport along the $\na$-geodesics
$$\tau\lra\exp_{u_{\ups,\xi_{\nu_{\es}}(\ups_1)}}\!\tau\ze, \qquad \tau\!\in\![0,1].$$
Similarly to~\e_ref{torus_e6a}, 
\begin{equation}\label{gwdiff_e6a}\begin{split}
\Phi_t(\ups,\ze) &\equiv\Pi_{\ze}^{-1}
\{\bpar_J\!+\!\nu_{\es}\!+\!t\nu\}\exp_{u_{\ups,\xi_{\nu_{\es}}(\ups_1)}}\!\ze\\
&=\{\bpar_J\!+\!\nu_{\es}\}u_{\ups}+D_{J,\nu_{\es};\ups;\xi_{\nu_{\es}}(\ups_1)}\ze
  +t\,\nu\big(\Si_{\ups},u_{\ups,\xi_{\nu_{\es}}(\ups_1)}\big)
  +N_{\ups;t}\xi,
\end{split}\end{equation}
where the quadratic term $N_{\ups;t}$ satisfies
\begin{gather}\label{gwdiff_e6b}
N_{\ups;t}0\!=\!0, \qquad
\big\|N_{\ups;t}\ze-N_{\ups;t}\ze'\big\|_{\ups,p} 
\le C(b)\big(t\!+\!\|\ze\|_{\ups,p,1}\!+\!\|\ze'\|_{\ups,p,1}\big)
\big\|\ze\!-\!\ze'\big\|_{\ups,p,1}\\
\forall\quad t\!\in\![0,\de(b)],~
\ups\!=\!(b,v)\!\in\!\ti\cF\T_{\de}^{\eset}, ~
\ze,\ze'\in\!\Ga\big(\ups;\xi_{\nu_{\es}}(\ups_1)\big) ~\st~ 
\|\ze\|_{\ups,p,1},\|\ze'\|_{\ups,p,1}\le\de(b).\notag
\end{gather}
Thus, by Lemma~\ref{gwdiff_lmm1}, the analogues of~\e_ref{torus_e4} 
and~\e_ref{torus_e8a} mentioned above, 
for every precompact open subset $K$ of $\U_{\T,\nu}(X;A)$ there exist 
$\de_K,C_K\!\in\!\R^+$ and a neighborhood of $U_K$ of $K$ in $\X_{1,k}(X,A)$ such that 
for all $t\!\in\![0,\de_K]$
\begin{gather}\label{gwdiff_e5a}
\M_{1,k}^0(X,A;J,\nu_{\es}\!+\!t\nu) \!\cap\! U_K
\approx \big\{[(\ups,\ze)]\!\in\!\Om_K(t)\!:\Phi_t(\ups,\ze)\!=\!0\big\},\\
\hbox{where}\quad 
\Om_K(t)=\big\{[(\ups,\ze)]\!: \ups\!=\!(b,v)\!\in\!\ti\cF\T_{\de_K}^{\eset}|_K;~
\ze\!\in\!\Ga_+(\ups;\xi_{\nu_{\es}}(\ups_1)),\,
 \|\ze\|_{\ups,p,1}\!\le\!C_K\big(|\rho(\ups)|\!+\!t\big)\big\}.\notag
\end{gather}\\

\noindent
For each $b\!\in\!\U_{\T,\nu_{\es}}^{(0)}(X;J)$, let
$$\bar\nu_b \in \Ga_-^{0,1}(b;J)$$
be the $L^2$-projection of $\nu(b)$.
We note that the~map
$$\G_{1,k}^{0,1}(X,A;J)\lra 
 \Ga\big(\U_{\T,\nu_{\es}}(X;J);\pi_P^*\E^*\!\otimes\!\ev_P^*TX\big), \qquad
 \nu\lra\bar\nu,$$
is surjective for every bubble type~$\T$.
By \e_ref{gwdiff_e6a}, \e_ref{gwdiff_e6b}, and 
the analogues of~\e_ref{g1reg_lmm1e6}, \e_ref{torus_e5}, and~\e_ref{torus_e8b},
\begin{gather}\label{gwdiff_e5b}
\big\|\pi^{0,1}_{\ups;-}\Phi_t(\ups,\ze)-(\cD_{\T}\rho(\ups)\!+\!t\bar\nu_b)\big\|
\le \ve(\ups)\big(|\rho(\ups)|\!+\!t) \\
\forall\quad t\!\in\![0,\de_K],\, t\!\in\![0,\de_K],\,
\ups\!=\!(b,v)\!\in\!\ti\cF\T_{\de_K}^{\eset}|_K,\,
[(\ups,\ze)]\!\in\!\Om_K(t),\notag
\end{gather}
for some function
$$\ve\!:\ti\cF\T^{\eset}\lra\R^+ \qquad\hbox{s.t.}\qquad
\lim_{|\ups|\lra0}\ve(\ups)=0.$$\\

\noindent
We denote by
$$\pi_{\bP\F\T}\!: \bP\F\T\lra \U_{\T,\nu_{\es}}(X;J)
\qquad\hbox{and}\qquad
\ga_{\bP\F\T}\!\lra \bP\F\T$$
the bundle projection map and the tautological line bundle.
If $\nu_{\es}\!\in\!\G_{1,k}^{\gd}(X,A;J)$ is generic, the section
$$\ti\cD_{\T}\in\Ga\big(\U_{\T,\nu_{\es}}(X;J);
\Hom(\ga_{\bP\F\T};\pi_{\bP\F\T}^*(\pi_P^*\E^*\!\otimes\!\ev_P^*TX)\big)$$
induced by $\cD_{\T}$ is transverse to the zero set.
For a generic choice of the pseudocycles~$f_l$, this is also the case for 
the restriction of~$\ti\cD_{\T}$ to $\bP\F\T\big|_{\U_{\T,\nu_{\es}}(X;J;f_{\la})}$
for every $\la\!\in\!\A$.
On the other hand,
\begin{equation}\label{dimcount_e}\begin{split}
\dim\,\U_{\T,\nu_{\es}}(X;J;f_{\la})
 &=\dim\,\U_{\T,\nu_{\es}}(X;J)+\big(\!\dim Y_{\la}\!-\!nk\big)\\
&\le\big(\!\dim_{1,k}(X,A)+2(n\!-\!|\hat{I}|\!-\!|\aleph|)\big)-\dim_{1,k}(X,A)\\
&=2(n\!-\!|\chi(\T)|)-2\big(|\hat{I}\!-\!\chi(\T)|\!+\!|\aleph|\big).
\end{split}\end{equation}
The middle inequality is an equality if and only if $\la\!=\!\mn$.
Thus, the section $\ti\cD_{\T}$ does not vanish over 
$\bP\F\T\big|_{\U_{\T,\nu_{\es}}(X;J;f_{\la})}$.
This is equivalent to saying that the bundle homomorphism
$$\cD_{\T}\!: \F\T\lra\pi_P^*\E^*\!\otimes\!\ev_P^*TX$$
is \sf{nondegenerate} over $\U_{\T,\nu_{\es}}(X;J;f_{\la})$,
i.e.~is injective on every fiber over $\U_{\T,\nu_{\es}}(X;J;f_{\la})$.\\

\noindent
Suppose $\T$ is not a simple bubble type or $\la\!\neq\!\mn$. 
By~\e_ref{dimcount_e},
$$\frac{1}{2}\dim\,\U_{\T,\nu_{\es}}(X;J;f_{\la}) +\rk\,\F\T 
<  \rk\,\pi_P^*\E^*\!\otimes\!\ev_P^*TX.$$
Thus, for a generic $\nu\!\in\!\G_{1,k}^{\gd}(X,A;J)$, the affine bundle map
$$\F\T\lra\pi_P^*\E^*\!\otimes\!\ev_P^*TX, \qquad
\ups\!=\![b,v]\lra \cD_{\T}\ups\!+\!\bar\nu_b,$$
has no zeros over $\U_{\T,\nu_{\es}}(X;J;f_{\la})$.
Since $\cD_{\T}$ is nondegenerate over $\U_{\T,\nu_{\es}}(X;J;f_{\la})$,
\e_ref{gwdiff_e5b} and
the proof of Lemma~3.2 in~\cite{g2n2and3} then imply that 
for every compact subset $K$ of $\U_{\T,\nu_{\es}}(X;J;f_{\la})$ there exist 
$\de_K\!\in\!\R^+$ and a neighborhood $U_K'$ of $K$ in 
$\U_{\T,\nu_{\es}}(X;J)\!\times\!\bar{Y}$ such~that
$$\big\{[(\ups,\ze),z]\!\in\!\Om_K(t)|_{U_K'}\!: 
\pi^{0,1}_{\ups;-}\Phi_t(\ups,\ze)\!=\!0\big\}
=\eset \qquad\forall~t\!\in\!(0,\de_K).$$
Thus, there exists a neighborhood $U_K$ of $K$ in $\X_{1,k}(X,A)\!\times\!\bar{Y}$ such that 
$$\big(\M_{1,k}^0(X,A;J,\nu_{\es}\!+\!t\nu)\!\times\!\bar{Y}\big) 
\cap U_K=\eset.$$
The proof of Proposition~\ref{bdcontr_prp} in the case 
$\T$ is not a simple bubble type or $\la\!\neq\!\mn$ is now complete.\\

\noindent
For remainder of this subsection we assume that
$\T$ is a simple bubble type and $\la\!=\!\mn$.
If $\ups\!=\!(b,v)\!\in\!\ti\cF\T_{\de}^{\eset}$, we denote by 
$$\Ga_-^{0,1}(\ups)\subset\Ga^{0,1}\big(\ups;\xi_{\nu_{\es}}(\ups_1)\big)$$
the image of $\Ga_-^{0,1}(b;J)$ under $R_{\ups}$ and by $\Ga_+^{0,1}(\ups)$
the $L^2$-orthogonal complement of $\Ga_-^{0,1}(\ups)$ in
$\Ga^{0,1}\big(\ups;\xi_{\nu_{\es}}(\ups_1)\big)$.
Let 
$$\pi^{0,1}_{\ups;+}\!: \Ga^{0,1}\big(\ups;\xi_{\nu_{\es}}(\ups_1)\big)
\lra \Ga_+^{0,1}(\ups)$$
the $L^2$-projection map.
Since $\aleph\!=\!\eset$,
$$\|R_{\ups,\xi_{\nu_{\es}}(\ups_1)}\eta\|_{\ups,p} \le 
C(b)\|R_{\ups,\xi_{\nu_{\es}}(\ups_1)}\eta\|_{\ups,2}
\qquad\forall~ \ups\!=\!(b,v)\!\in\!\ti\cF\T_{\de}^{\eset},\,
\eta\!\in\!\Ga_-^{0,1}(b;J),$$
for some $C\!\in\!C^{\i}(\U_{\T,\nu_{\es}}(X;J);\R^+)$.
It follows that
$$\big\|\pi^{0,1}_+\eta'\big\|_{\ups,p} \le C(b)\|\eta'\|_{\ups,p}
\qquad\forall~ \ups\!=\!(b,v)\!\in\!\ti\cF\T_{\de}^{\eset},\,
\eta'\!\in\!\Ga^{0,1}\big(\ups;\xi_{\nu_{\es}}(\ups_1)\big).$$
Thus, by the analogues of \e_ref{torus_e4} and~\e_ref{torus_e5},
\begin{gather}\label{gwdiff_e7}
C(b)^{-1}\|\ze\|_{\ups,p,1} \le 
\big\|\pi^{0,1}_+D_{J,\nu_{\es};\ups;\xi_{\nu_{\es}}(\ups_1)}\ze\big\|_{\ups,p}
\le C(b)\|\ze\|_{\ups,p,1} \\ 
\forall\quad  \ups\!=\!(b,v)\in\!\ti\cF\T_{\de}^{\eset}, 
~\ze\!\in\!\Ga_+\big(\ups;\xi_{\nu_{\es}}(\ups_1)\big).\notag
\end{gather}
In particular, the operator
$$\pi^{0,1}_+D_{J,\nu_{\es};\ups;\xi_{\nu_{\es}}(\ups_1)}\!:
\Ga_+\big(\ups;\xi_{\nu_{\es}}(\ups_1)\big) \lra \Ga_+^{0,1}(\ups)$$
is an isomorphism.
Its norm and the norm of its inverse are dependent only on 
$[b]\!\in\!\U_{\T,\nu_{\es}}(X;J)$.
By \e_ref{gwdiff_e6a}, \e_ref{gwdiff_e6b}, \e_ref{gwdiff_e7}, 
the analogue of~\e_ref{torus_e8a}, and the Contraction Principle,
for every compact subset $K$ of $\U_{\T,\nu_{\es}}(X;J)$ there exists 
$\de_K\!\in\!\R^+$ such that for all 
$$\ups\in\ti\cF\T_{\de_K}^{\eset}|_K \qquad\hbox{and}\qquad
t\in[0,\de_K]$$
the equation
$$\pi^{0,1}_+\Phi_t(\ups,\ze)=0, \qquad \ze\!\in\!
\Ga_+\big(\ups;\xi_{\nu_{\es}}(\ups_1)\big),$$
has a unique small solution $\ze_t(\ups)$.
Furthermore, $\ze_t(\ups)\!\in\!\Om_K(t)$.\\

\noindent
By the above, for every compact subset $K$ of $\U_{\T,\nu_{\es}}(X;J;f_{\mn})$
there exist $\de_K\!\in\!\R^+$ and small neighborhoods $U_K'$ and $U_K$ 
of $K$ in $\U_{\T,\nu_{\es}}(X;J)\!\times\!Y_{\mn}$ and $\X_{1,k}(X,A)\!\times\!\bar{Y}$, 
respectively,
such that
\begin{equation*}\begin{split}
\M_{1,k}^0(X,A;J,\nu_{\es}\!+\!t\nu;f_{\mn}) \!\cap\! U_K \approx
\big\{&([\ups],z)\!\in\!\cF\T_{\de_K}^{\eset}|_K\!\times\!Y_{\mn}\!: \pi^{0,1}_{\ups;-}\Phi_t(\ups,\ze_t(\ups))\!=\!0;\\
&\qquad\qquad\big\{\ev\!\times\!f\big\}\big(
(\Si_{\ups},\exp_{u_{\ups,\xi_{\nu_{\es}}(\ups_1)}}\!\ze_t(\ups)),z)\!\in\!\De_X^k\big\}
\end{split}\end{equation*}
for all $t\!\in\![0,\de_K]$.
On the other hand, the bundle homomorphism $\cD_{\T}$ is regular over 
$\bar\U_{\T,\nu_{\es}}(X;J;f)$ by the $m\!=\!1$ case of \e_ref{derivexp_e1} 
and~\e_ref{derivexp_e2}, i.e.~$\cD_{\T}$ can be approximated by a polynomial
on the normal bundle near every stratum of $\bar\U_{\T,\nu_{\es}}(X;J;f)$;
see  \cite[Definition~3.9]{g2n2and3}.
Since $\nu$ is generic, 
$\cD_{\T}$ is nondegenerate over $\U_{\T,\nu_{\es}}(X;J;f_{\mn})$,
and $\ev\!\times\!f$ is transverse to~$\De_X^k$ in $(X^2)^k$, \e_ref{gwdiff_e5a} and 
the proof of Lemma~3.5 in~\cite{g2n2and3} then imply that there exists a compact subset
$K_{\nu}$ of $\U_{\T,\nu_{\es}}(X;J;f_{\mn})$ with the following property.
For every compact subset $K$ of $\U_{\T,\nu_{\es}}(X;J;f_{\mn})$ 
that contains~$K_{\nu}$ there exist $\de_K\!\in\!\R^+$ 
and a neighborhood $U_K$ of $K$ in $\X_{1,k}(X,A)\!\times\!\bar{Y}$ such that 
$$^{\pm}\!\big|\M_{1,k}^0(X,A;J,\nu_{\es}\!+\!t\nu;f_{\mn})\!\cap\!U_K\big|
=N(\cD_{\T}) \qquad\forall\,t\!\in\![0,\de_K].$$
We note that $K_{\nu}$ can taken to be any compact subset of 
$\U_{\T,\nu_{\es}}(X;J;f_{\mn})$ such that all of the finitely many zeros of the affine map
$$\cD_{\T}\!+\!\bar\nu\!: \F\T\lra\pi_P^*\E^*\!\otimes\!\ev_P^*TX$$
over $\U_{\T,\nu_{\es}}(X;J;f_{\mn})$ lie in~$\F\T|_{K_{\nu}}$.

\subsection{Counting Zeros of Affine Bundle Maps}
\label{affinemap_subs}

\noindent
In this subsection we conclude the proof of Theorem~\ref{gwdiff_thm} 
by computing the numbers $N(\cD_{\T})$ when the dimension of $X$ is $4$ or~$6$.
Using the method of \cite[Section~3]{g2n2and3}, $N(\cD_{\T})$ can be computed for 
arbitrary-dimensional symplectic manifolds~$X$ and more general cohomology classes~$\psi$;
this is the main subject of~\cite{g1diff}.
In order to avoid introducing quite a bit of additional notation in this paper 
we restrict the computation to the special cases of Theorem~\ref{gwdiff_thm}.\\

\noindent
If $\T$ is a bubble type as in~\e_ref{bubbtype_e}, let
$$M_P\T=\big\{l\!\in\![k]\!:j_l\!\not\in\!\hat{I}\big\} \quad
\hbox{and}\quad
\Aut^*(\T)=\Aut(\T)/\{g\!\in\!\Aut(\T)\!:g\cdot h\!=\!h~\forall\, h\!\in\!I_1\}.$$
For each $i\!\in\!I_1$, let
$$\ti\T_i=\big(\ti{M}_i\T,\{i\}\!\cup\!\ti{H}_i\T;
j|_{\ti{M}_i\T},\under{A}|_{\{i\}\cup\ti{H}_i\T}\big).$$
If $\nu_{\es}\!\in\!\G_{1,k}^{\gd}(X,A;J)$ and $\aleph\!=\!\eset$, then
\begin{equation}\label{g1gendecomp_e3}
\U_{\T,\nu_{\es}}(X;J)\approx \big(\cM_{1,M_P\T\sqcup I_1}
\!\times\!\U_{\bar\T,\nu_{\es}}(X;J)\big)\big/\Aut^*(\T),
\end{equation}
where $\cM_{1,M_P\T\sqcup I_1}$ is the subspace of the moduli space 
$\ov\cM_{1,M_P\T\sqcup I_1}$ consisting of smooth curves and
\begin{equation}\label{g1gendecomp_e4}\begin{split}
\U_{\bar\T,\nu_{\es}}(X;J)=
\big\{(b_i)_{i\in I_1}\!\in\!\prod_{i\in I_1}\!\U_{\T_i,\nu_{B;i}}\!:
\ev_0(b_{i_1})\!=\!\ev_0(b_{i_2})\, \forall i_1,i_2\!\in\!I_1\big\}\qquad&\\
\subset \prod_{i\in I_1}\!\ov\M_{0,\{0\}\cup\ti{M}_i\T}(X,\ti{A}_i;J,\nu_{B;i}),&
\end{split}\end{equation}
for some $\nu_{B;i}\!\in\!\G_{0,\{0\}\cup\ti{M}_i\T}^{\es}(X,\ti{A}_i;J)$.
This decomposition is illustrated in Figure~\ref{g1gendecomp_fig}.
In this figure, we represent an entire stratum $\U_{\T,\nu_{\es}}(X;J)$
of bubble maps by the domain of the stable maps in $\U_{\T,\nu_{\es}}(X;J)$.
The right-hand side of Figure~\ref{g1gendecomp_fig} 
represents the subset of the cartesian product of the three spaces
of bubble maps, corresponding to the three drawings,
on which the appropriate evaluation maps agree,
as indicated by the solid line and defined in~\e_ref{g1gendecomp_e4}.\\

\begin{figure}
\begin{pspicture}(-1.1,-2)(10,1.25)
\psset{unit=.4cm}
\rput{45}(0,-4){\psellipse(5,-1.5)(2.5,1.5)
\psarc[linewidth=.05](5,-3.3){2}{60}{120}\psarc[linewidth=.05](5,0.3){2}{240}{300}
\pscircle[fillstyle=solid,fillcolor=gray](5,-4){1}\pscircle*(5,-3){.2}
\pscircle[fillstyle=solid,fillcolor=gray](6.83,.65){1}\pscircle*(6.44,-.28){.2}
\pscircle[fillstyle=solid,fillcolor=gray](3.17,.65){1}\pscircle*(3.56,-.28){.2}
\pscircle*(2.5,-1.5){.2}\pscircle*(7.5,-1.5){.2}}
\rput(7.1,.6){$y_1$}\rput(3,-4){$y_2$}
\rput(19,-1.5){$\approx\qquad\qquad\cM_{1,5}\qquad\times$}
\pscircle[fillstyle=solid,fillcolor=gray](29,2){1}\pscircle*(28,2){.2}
\pscircle[fillstyle=solid,fillcolor=gray](29,-1){1}\pscircle*(28,-1){.2}
\pscircle[fillstyle=solid,fillcolor=gray](29,-4){1}\pscircle*(28,-4){.2}
\psline[linewidth=.03](28,2)(28,-4)
\end{pspicture}
\caption{An example of the decomposition~\e_ref{g1gendecomp_e3}}
\label{g1gendecomp_fig}
\end{figure}

\noindent
We next define a similar splitting for the space $\U_{\T,\nu_{\es}}(X;J;f_{\mn})$.
We can assume that the maps $\{f_{l;\mn}\}_{l\in M_P\T}$ intersect transversally.
Let
$$f_0^{\T}\!: \bar{Y}_0^{\T} \lra \bigcap_{l\in M_P\T}\!\!\!\!f_l(\bar{Y}_l)
~\subset X$$
be a pseudocycle representative for 
$$ \bigcap_{l\in M_P\T}\!\!\!\!\PD_X\psi_l ~\in H_*(X;\Z)$$
such that
$$f_{0;\mn}^{\T}\!: Y_{0;\mn}^{\T} \lra 
\bigcap_{l\in M_P\T}\!\!\!\!f_l(\bar{Y}_{l;\mn})$$
is one-to-one.
We put
\begin{gather*}
\ev^{\T}\!=\ev_0\!\times\!\!\!\prod_{l\in[k]-M_P\T}\!\!\!\!\!\!\!\!\!\ev_l:
\U_{\bar\T,\nu_{\es}}(X;J)\lra X\!\times\!\prod_{l\in[k]-M_P\T}\!\!\!\!\!\!\!\!\!X~,\\ 
f^{\T}\!=f_0^{\T}\!\times\!\prod_{l\in[k]-M_P\T}\!\!\!\!\!\!\!\!\!f_l: 
\bar{Y}^{\T}\!\equiv\bar{Y}_0^{\T}\!\times\!\!\!
\prod_{l\in[k]-M_P\T}\!\!\!\!\!\!\!\!\!\bar{Y}_l
~\lra X\!\times\!\prod_{l\in[k]-M_P\T}\!\!\!\!\!\!\!\!\!X~,
\quad\hbox{and}\quad  Y_{\mn}^{\T}\!=
Y_{0;\mn}^{\T}\!\times\!\! \prod_{l\in[k]-M_P\T}\!\!\!\!\!\!\!\!\!Y_{l;\mn}.
\end{gather*}\\
Similarly to Subsection~\ref{gwdiffsumm_subs}, let
\begin{alignat*}{1}
\bar\U_{\bar\T,\nu}(X;J;f)&= \big\{(b,z)\!\in\!
\bar\U_{\bar\T,\nu}(X;J)\!\times\!\bar{Y}^{\T}\!:
\ev^{\T}(b)\!=\!f^{\T}(z)\big\}, \qquad\hbox{and}\\
\U_{\bar\T,\nu}(X;J;f_{\mn}) &= \big(\U_{\bar\T,\nu}(X;J)\!\times\!Y_{\mn}^{\T}\big)
\cap\bar\U_{\bar\T,\nu}(X;J;f).
\end{alignat*}
We have
\begin{equation}\label{g1gendecomp_e5}
\U_{\T,\nu_{\es}}(X;J;f_{\mn})\approx \big(\cM_{1,M_P\T\sqcup I_1}
\!\times\!\U_{\bar\T,\nu_{\es}}(X;J;f_{\mn})\big)\big/\Aut^*(\T).
\end{equation}\\

\noindent
If $\T$ is a simple bubble type, we define the homomorphism
$$\wt\cD_{\T}\!: \wt\F\T\lra \pi_P^*\E^*\!\otimes\!\ev_0^*TX$$
over $\ov\cM_{1,M_P\T\sqcup I_1}\!\times\!\bar\U_{\bar\T,\nu_{\es}}(X;J;f)$
similarly to the homomorphism $\cD_{\T}$.
By~\e_ref{g1gendecomp_e5},
\begin{equation}\label{gwdiff_e11}
N(\cD_{\T})=N(\wt\cD_{\T})\big/ \big|\Aut^*(\T)\big|.
\end{equation}
By Proposition~\ref{bdcontr_prp} and \e_ref{gwdiff_e11}, the difference
between the standard and reduced genus-one GW-invariants of Theorem~\ref{gwdiff_thm}
is determined by the numbers $N(\wt\cD_{\T})$, where $\T$ is a simple bubble type
as in~\e_ref{bubbtype_e}.
We will compute these numbers in the two special cases of Theorem~\ref{gwdiff_thm}.\\

\noindent
First, we note~that if $\T$ is a bubble type as in~\e_ref{bubbtype_e},
not necessarily simple, and $\aleph\!=\!\eset$,
\begin{equation}\label{dimcount_e2}\begin{split}
\dim\,\bar\U_{\bar\T,\nu_{\es}}(X;J;f)
&=\dim\,\bar\U_{\T,\nu_{\es}}(X;J;f)-\dim\,\ov\cM_{1,M_P\T\sqcup I_1}\\
&=2(n\!-\!|\chi(\T)|\!-\!|I_1|)-2\big(|\hat{I}\!-\!\chi(\T)|\!+\!|M_P\T|\big),
\end{split}\end{equation}
by~\e_ref{dimcount_e}.
In particular, if $n\!=\!2$, then
$$\bar\U_{\bar\T,\nu_{\es}}(X;J;f)\neq\eset \qquad\Lra\qquad
|\chi(\T)|=1,~~   \chi(\T)=I_1=\hat{I}, ~~M_P\T=\eset.$$
Furthermore, if $\bar\U_{\bar\T,\nu_{\es}}(X;J;f)$ is nonempty,
it is a finite collection of points. 
In this case, $\wt\cD_{\T}$ is the bundle homomorphism
$$\pi_P^*\big(s_1\!\oplus\!0\big) \!: 
\pi_P^*L_1\lra \pi_P^*\big(\E^*\!\oplus\!\E^*\big)$$
over $\ov\cM_{1,1}\!\times\!\bar\U_{\bar\T,\nu_{\es}}(X;J;f)$.
Thus,
$$N(\wt\cD_{\T})=-\frac{1}{24}\,^{\pm}\!\big|\bar\U_{\bar\T,\nu_{\es}}(X;J;f)\big|.$$
Since $^{\pm}\!\big|\bar\U_{\bar\T,\nu_{\es}}(X;J;f)\big|
=^{\pm}\!\big|\bar\U_{\bar\T,\nu'}(X;J;f)\big|$
for any sufficiently small $\nu'$ such that the restriction of the bundle section $\bar\partial\!+\!\nu'$ to every stratum of $\X_{0,\{0\}\sqcup[k]}(X,A)$ 
is transverse to the zero set, we can take $\nu'\!=\!\pi^*\nu_0$, where
$$\pi\!:\X_{0,\{0\}\sqcup[k]}(X,A)\lra \X_{0,[k]}(X,A)$$
is the forgetful map and $\nu_0$ is a small generic 
deformation of $\bar\partial$ on $\X_{0,[k]}(X,A)$.
Since 
$$\ov\M_{0,k}(X,A;J,\nu_0;f)=\eset$$
for a generic $\nu_0$ for dimensional reasons, it follows 
$\bar\U_{\bar\T,\nu'}(X;J;f)\!=\!\eset$ as well and thus 
$N(\wt\cD_{\T})=0$, as claimed by  $n\!=\!2$ case
of Theorem~\ref{gwdiff_thm}.\\

\noindent
Suppose $n\!=\!3$. By~\e_ref{dimcount_e2},
$$\bar\U_{\bar\T,\nu_{\es}}(X;J;f)\neq\eset \qquad\Lra\qquad
|\chi(\T)|=1,~~   
\big(|I_1|\!-\!|\chi(\T)|\big)+\big|\hat{I}\!-\!\chi(\T)\big|
+|M_P\T|\in\{0,1\}.$$
If $\T$ is a simple bubble type and $\bar\U_{\bar\T,\nu_{\es}}(X;J;f)$ is nonempty,
it follows that
$$|\chi(\T)|=1 \qquad\hbox{and}\qquad |M_P\T|\in\{0,1\}.$$
If $|M_P\T|\!=\!1$, $\bar\U_{\bar\T,\nu_{\es}}(X;J;f)$ is again a finite collection of points. 
In this case, $\wt\cD_{\T}$ is the bundle homomorphism
$$\pi_P^*\big(s_1\!\oplus\!0\!\oplus\!0\big) \!: 
\pi_P^*L_1\lra \pi_P^*\big(\E^*\!\oplus\!\E^*\!\oplus\!\E^*\big)$$
over $\ov\cM_{1,2}\!\times\!\bar\U_{\bar\T,\nu_{\es}}(X;J;f)$.
Thus, as in the $n\!=\!2$ case above, 
$N(\wt\cD_{\T})\!=\!0$.\footnote{As in the $n\!=\!2$ case, $^{\pm}\!\big|\bar\U_{\bar\T,\nu_{\es}}(X;J;f)\big|\!=\!0$.
Furthermore, since $c_1(\E^*)^2\!=\!0$ on $\ov\cM_{1,2}$,
we can choose $\nu$ so that~$\bar\nu$ does not vanish and thus $N(\wt\cD_{\T})\!=\!0$
by definition.}\\

\noindent
Finally, suppose $n\!=\!3$, $\T$ is a simple bubble type, and $|M_P\T|\!=\!0$.
If $i$ is the unique element of~$\chi(\T)$, 
$\wt\cD_{\T}$ is the bundle homomorphism
$$\pi_P^*s_1\!\otimes\!\pi_i^*\cD_0\!: \pi_P^*L_1\!\otimes\!\pi_i^*L_0
\lra \pi_P^*\E^*\!\otimes\!\pi_i^*\ev_0^*TX$$
over 
$$\ov\cM_{1,1}\times \bar\U_{\bar\T,\nu_{\es}}(X;J;f)
=\ov\cM_{1,1}\times \ov\M_{0,\{0\}\cup[k]}(X,A;J,\nu_{B;i};f).$$
It can be assumed that $\nu_{B,i}\!=\!\pi^*\nu_B$ for a generic element
$\nu_B\!\in\!\G_{0,\{0\}}^{\es}(X,A;J)$, where 
$$\pi\!: \X_{0,\{0\}\sqcup[k]}(X,A)\lra  \X_{0,\{0\}}(X,A)$$
is the forgetful map.
Then,
\begin{equation}\label{piim_e}
\pi\!: \ov\M_{0,\{0\}\cup[k]}(X,A;J,\nu_{B;i};f)
\lra\ov\M_{0,\{0\}}(X,A;J,\nu_B)
\end{equation}
and there are identifications $L_0\!=\!\pi^*L_0$ and $\cD_0\!=\!\pi^*\cD_0$.
Thus,
\begin{equation}\label{Ntwist_e}
N(\wt\cD_{\T})=N\big(\pi_P^*s_1\!\otimes\!\pi_i^*\pi^*\cD_0\big).
\end{equation}\\

\noindent
On the other hand, if $\nu_B$ is generic, the image of the projection $\pi$
in \e_ref{piim_e} is contained in $\M_{0,\{0\}}^{\{0\}}(X,A;J,\nu_B)$,
while the restriction of $\cD_0$ to every stratum of $\M_{0,\{0\}}^{\{0\}}(X,A;J,\nu_B)$
is transverse to the zero set.
Thus, $\pi^*\cD_0$ does not vanish over 
$\ov\M_{0,\{0\}\cup[k]}(X,A;J,\nu_{B;i};f)$
for dimensional reasons if $\nu_B$ is generic.
Since $s_1$ does not vanish over $\ov\cM_{1,1}$, the bundle homomorphism
$$\pi_P^*s_1\!\otimes\!\pi_i^*\pi^*\cD_0\!: \pi_P^*L_1\!\otimes\!\pi_i^*\pi^*L_0
\lra \pi_P^*\E^*\!\otimes\!\pi_i^*\pi^*\ev_0^*TX$$
does not vanish over 
$\ov\cM_{1,1}\times \ov\M_{0,\{0\}\cup[k]}(X,A;J,\nu_{B;i};f)$.
Thus, by Lemma~3.14 in~\cite{g2n2and3},
\begin{equation}\label{bdcontr_e1}\begin{split}
N(\wt\cD_{\T})&=\blr{c(\pi_P^*\E^*\!\otimes\!\pi_i^*\ev_0^*TX)
c(\pi_P^*L_1\!\otimes\!\pi_i^*\pi^*L_0)^{-1},
\ov\cM_{1,1}\!\times\!\ov\M_{0,\{0\}\cup[k]}(X,A;J,\nu_{B;i};f)}\\
&=-\frac{1}{24}\blr{c_1(TX)+\pi^*c_1(L_0^*),
\ov\M_{0,\{0\}\cup[k]}(X,A;J,\nu_{B;i};f)}.
\end{split}\end{equation}
By the divisor and dilaton relations for GW-invariants, 
\begin{equation}\label{bdcontr_e8}\begin{split}
\blr{c_1(TX),\ov\M_{0,\{0\}\cup[k]}(X,A;J,\nu_{B;i};f)}
&=\lr{c_1(TX),A}\cdot\GW_{0,k}^X(A;\psi),\\
\blr{\pi^*c_1(L_0^*),\ov\M_{0,\{0\}\cup[k]}(X,A;J,\nu_{B;i};f)}
&=-2\cdot\GW_{0,k}^X(A;\psi);
\end{split}\end{equation}
see Section~26.3 in~\cite{MirSym}, for example.
Combining \e_ref{Ntwist_e} with \e_ref{bdcontr_e1} and \e_ref{bdcontr_e8}, 
we conclude that
$$ N(\wt\cD_{\T})= \frac{2\!-\!\lr{c_1(TX),A}}{24}\GW_{0,k}^X(A;\psi).$$
The proof of Theorem~\ref{gwdiff_thm} is now complete.

\vspace{.2in}

{\it 
\begin{tabbing}
${}\qquad$
\= Department of Mathematics, Stanford University, Stanford,
CA 94305-2125\\ 
\> \textnormal{Current Address:} Department of Mathematics, SUNY, Stony Brook, NY 11794-3651\\
\> azinger@math.sunysb.edu
\end{tabbing}}

\end{document}